\documentclass[11pt]{amsart}
\usepackage{amsmath,amsfonts,amssymb,amsxtra,latexsym,amscd,enumerate,amsthm}

\usepackage{latexsym}
\usepackage{epsfig}

\def\z{\zeta}
\def\t{\theta}
\def\r{\rho}
\def\g{\gamma}
\def\a{\alpha}
\def\q{\frac{1}{2}}

\def\b{\beta}
\def\e{\varepsilon}

\def\D{\Delta}
\def\Ra{\Rightarrow}

\def\R{\mathbb{R}}
\def\C{\mathbb{C}}
\def\Z{\mathbb{Z}}
\def\N{\mathbb{N}}

\def\beq{\begin{equation}}
\def\eeq{\end{equation}}

\def\beq{\begin{equation}}
\def\eeq{\end{equation}}

\newtheorem{t1}{Theorem}
\newtheorem{l1}{Lemma}
\newtheorem{p2}{Proposition}

\newtheorem{c1}{Corollary}

\newtheorem{c3}{Corollary}
\newtheorem{r1}{Remark}

\newtheorem{r3}{Remark}

\begin{document}

\title[Quadratic Nonlinear Derivative Schr\"odinger Equations]{Quadratic Nonlinear Derivative Schr\"odinger Equations - Part 1}

\author{Ioan Bejenaru}
\address{Department of Mathematics, UCLA, Los Angeles CA 90095-1555}
\email{bejenaru@math.ucla.edu}

\begin{abstract} 

In this paper we consider the local well-posedness theory for the quadratic nonlinear Schr\"odinger equation with low regularity initial data in the case when the nonlinearity contains derivatives. We work in $2+1$ dimensions and prove a local well-posedness result up to the scaling for small initial data with some spherical symmetry structure.

\end{abstract}

\maketitle

\section{Introduction}

This work is concerned with  the initial value problem for the nonlinear Schr\"odinger equations which generically have the form:

\begin{equation} \label{E:s}
\begin{cases}
\begin{aligned} 
iu_{t}-\Delta u &= P(u,\bar{u},\nabla u, \nabla \bar{u}), \   t\in \mathbb{R}, x\in \mathbb{R}^{n} \\
u(x,0) &=u_{0}(x)
\end{aligned}
\end{cases}
\end{equation}

\noindent
where $u: \mathbb{R}^{n} \times \mathbb{R} \rightarrow \mathbb{C}$ and $P: \mathbb{C}^{2n+2} \rightarrow \mathbb{C}$ is a polynomial. 

We are interested in the theory of local well-posedness for this problem in Sobolev spaces. It is natural to discuss this in terms of the Taylor expansion of $P$ around $0$. 

Constant functions are not in $H^{s}$, hence the first natural condition to impose on $P$ is $P(0)=0$. The next step is to consider the linear terms. Those without derivatives are harmless for the local well-posedness theory because they have Lipschitz contribution in all Sobolev spaces. 

The problem becomes nontrivial when we have linear terms with derivatives. This is made clear by the following result due to Mizohata \cite{m}  which proves that for the problem:

\begin{equation}
\begin{cases}
\begin{aligned} 
iu_{t}-\Delta u &= b_{1}(x) \nabla u , \   t\in \mathbb{R}, x\in \mathbb{R}^{n} \\u(x,0) &=u_{0}(x)
\end{aligned}
\end{cases}
\end{equation}

the following condition on $b_{1}$ is necessary for the $L^{2}$ well-posedness theory:

\beq \label{co1}
\sup_{x \in \R^{n}, \omega \in \mathbb{S}^{n-1}, R > 0}  |Re \int_{0}^{R} b_{1}(x+r \omega) \cdot \omega dr| < \infty
\eeq

The idea behind this condition is that $Re \ b_{1}$ contributes to exponential growth of the solution along the flow.

One of the consequences for our problem is that we cannot have in the nonlinearity terms like $b_{1} \nabla u$ with $b_{1}$ real. On the other hand if the inhomogeneity is of type $b_{1} \nabla u$ with $b_{1}$ imaginary or $b_{2} \cdot \nabla \bar{u}$, the problem is  well-posed in any $H^{s}$, see \cite{k3} and \cite{k4}. 

The problem becomes more complicated once we have to deal with quadratic and higher order terms. There is a very general result due to Kenig, Ponce and Vega, see \cite{k2} which we summarize in what follows. 

\begin{t1}

Assume that $P$ has no constant or linear terms. Then there exist $s=s(n,P) > 0$ and $m=m(n,P) >0$ such that $\forall u_{0} \in H^{s}(\R^{n}) \cap L^{2}(\R^{n}:|x|^{2m}dx)$ the problem (\ref{E:s}) has a unique solution in $C([0,T]: H^{s} \cap L^{2}(\R^{n}:|x|^{2m}dx)$ where $T=T(||u||_{H^{s} \cap L^{2}(\R^{n}:|x|^{2m}dx) })$.
\end{t1}

 If $P$ does not contain quadratic terms, then the above authors also obtain a similar result without involving any decay, see  \cite{k2}.

 When we have terms with derivatives in the nonlinearity there is a loss of a derivative in the right-hand side of the equation which should be recovered. Hence one of the main ingredients in dealing with the problem is the local smoothing effect for Schr\"odinger equation which we describe bellow.

Let $(Q_{\a})_{\a \in \Z^{n}}$ be a system of disjoint cubes of size $R$ such that $\R^{n}= \cup_{\a \in \Z^{n}} Q_{\a}$. The homogeneous Schr\"odinger equation has a local gain of a half-derivative:

$$\sup_{\a \in \Z^{n}} ||D_{x}^{\q} e^{it\D}u_{0}||_{L^{2}(Q_{\a} \times \R)} \leq c R ||u_{0}||_{L^{2}(\R^{n})}$$  

The solution of the inhomogeneous Schr\"odinger equation

\begin{equation}
\begin{cases}
\begin{aligned} 
iu_{t}-\Delta u &= f , \   t\in \mathbb{R}, x\in \mathbb{R}^{n} \\u(x,0) &=0
\end{aligned}
\end{cases}
\end{equation}

\noindent
has a local gain of a derivative:

$$\sup_{\a \in \Z^{n}} ||\nabla_{x} u||_{L^{2}(Q_{\a} \times \R)} \leq c R\sum_{\a \in \Z^{n}} ||f||_{L^{2}(Q_{\a} \times \R)}$$  

These estimates can be found in \cite{k1} and there are even more refined versions in \cite{k2}. 

Let us briefly justify the need of decay. If we look at the problem:

\beq \label{ch}
iu_{t}-\Delta u = u \nabla{u}
\eeq

\noindent
then we need to recover an $l^{1}_{\a}$ structure for the term $u \nabla{u}$, while $\nabla{u}$ comes only with a $l^{\infty}_{\a}$ structure. Hence we need an $l^{1}_{\a}$ structure for $u$ and this can be achieved by a decay condition on $u$.

Another way is to think of $u \nabla{u}$ as being $b_{1}(t,x) \nabla{u}$. The result in (\ref{co1}) tells us that we need an integrability condition for $b_{1}=u$ which can be fulfilled  via a decay condition at infinity. 

In \cite{chr}, Christ gives a complete proof of ill-posedness of (\ref{ch}) in one-dimension with $u(0)=u_{0} \in H^{s}$ no matter how large we take $s$. This is because Sobolev regularity cannot be traded for decay.
  
Bringing in decay in Schr\"odinger equation forces some kind of trade-off. The linear equation conserves the $H^{s}$ structure of the initial data, but it does not conserve the decay structure. A simple way to see this is to consider the homogeneous equation with initial data $u_{0} e^{ix \xi_{0}}$, where $u_{0}$ is a smooth approximation of the characteristic function of the unit ball in $\R^{n}$. A standard approximation of the solution is $u_{0}(x-2t\xi_{0})  e^{ix \xi_{0}} e^{it \xi^{2}_{0}}$. At time $t=1$ this solution is supported at $|x| \approx |\xi_{0}|$, therefore we have to recover a decay of type $|\xi_{0}|^{m}$, while $u_{0}$ comes with no decay since it was supported around the origin. Therefore we have to spend derivatives in order to conserve the decay. Thus we should increase the regularity of the initial data more. 

This last remark brings us closer to the goal of this paper. The general result of Kenig, Ponce and Vega is not concerned with is the following question: what is the lowest Sobolev regularity the initial data can have so that we have well-posedness? When asking this question, one should be more specific about the type of the equation and the dimension of the space. 

The quadratic terms in $P$ are the ones we want to understand. If the nonlinearity contains only terms without derivatives then the problem is called semilinear Schr\"odinger equation (NLS). The generic quadratic (NLS) are those with the nonlinearity of type:  

$$u^{2}, |u|^{2}=u \cdot \bar{u}, (\bar{u})^{2}$$

If the nonlinearity contains terms with derivatives then the problem is called derivative non-linear Schr\"odinger equation (D-NLS). The generic quadratic (D-NLS) are those with the nonlinearity of type:  

$$u \cdot Du, \bar{u} \cdot Du, u \cdot D\bar{u}, \bar{u} \cdot D\bar{u}, Du \cdot Du, D\bar{u} \cdot Du, D\bar{u} \cdot D\bar{u} $$

To obtain a sharp result for these kind of problems means to obtain for each equation an index $s_{0}$ such that if $s > s_{0}$ and $u_{0} \in H^{s}$ (maybe with some additional structure) we have a well-posedness result and if  $s < s_{0}$ and $u_{0} \in H^{s}$ we have ill-posedness. What exactly is meant by ill-posedness is in itself a delicate issue. The case $s=s_{0}$ is the hardest one and it depends on the specific equation whether there is a positive or negative result.

An important concept for these problems is the scaling exponent, $s_{c}$. This is the exponent of the Sobolev space which scales the same way as the equation. Heuristically one would expect that $s_{0}=s_{c}$, but many times this is not the case. Let's become more specific. 

For the case of quadratic (NLS) we have $s_{c}=\frac{n}{2}-2$. There has been considerable progress in the study of quadratic (NLS) in two dimensions (here $s_{c}=-1$), see \cite{c1}, which led to the following results: if the nonlinearity is of type $u^{2}$ or $(\bar{u})^{2}$ then $s_{0}=-\frac{3}{4}$ and if the nonlinearity is of type $|u|^{2}$ then $s_{0}=-\frac{1}{4}$. The problem becomes ill-posed if $s < s_{0}$ in the following sense: the bilinear estimates fail to hold true.

The results for quadratic (D-NLS) did not yet reach this level of precision, the main difficulty being generated by the loss of one-derivative in the nonlinearity. What are the expectations for this problem? If only one of the terms contains derivatives (for instance $u Du$) then $s_{c}=\frac{n}{2}-1$. If both terms contain contain derivatives (for instance $Du Du$) then $s_{c}=\frac{n}{2}$. One new aspect in these cases is the use of decay: one needs a condition of the form $ \langle x \rangle^{m} u_{0} \in L^{2}$ in addition to $u_{0} \in H^{s}$. The best result up to date we know is of the form, see \cite{ch1}: if $m=\frac{n}{2}+2$ and $s=\frac{n}{2}+4$ then the quadratic (D-NLS) is locally well-posed.  This is a bit too far from the scaling exponent and, as we will see later on, the decay is too strong also.

There is one exception. For the terms $\bar{u} \cdot D \bar{u}$ and $D \bar{u} \cdot D \bar{u}$ the results were established up to the critical exponent, namely the existence of a solution was proved for every $s > s_{c}$, see \cite{g1}. This was possible  because of the following fact : the Fourier transforms of solutions for Schr\"odinger equations concentrate near the paraboloid $\tau=\xi^{2}$ and the most difficult estimates are the ones when the interacting elements are localized in frequency near the paraboloid and the result falls back near the paraboloid. The effect of the complex conjugation is that it replaces the paraboloid $\tau=\xi^{2}$ by the symmetric $\tau=-\xi^{2}$. The interaction via convolution of two paraboloids of type $\tau=-\xi^{2}$ is localized in the region $\tau \leq 0$ so it does not intersect the paraboloid $\tau=\xi^{2}$ (except the point $(0,0)$). This way the most difficult interactions mentioned above do not occur and the problem becomes easier.   

Let us return to the general form of the quadratic (D-NLS). The analysis of the problem brings the conclusion that the ``worst'' interactions are the orthogonal ones, i.e. those between waves which travel in orthogonal directions. In one-dimension this is not possible and this is why the problem becomes interesting once $n \geq 2$. As always, one tries to understand what happens when $n=2$ and then attempt to replicate the argument in higher dimensions. This is why in this work we will specialize to the case of two-dimension quadratic (D-NLS).   
 
Our goal is to obtain local well-posedness for initial data $u_{0} \in H^{s}$, for any $s > s_{c}$. With the techniques we involve it is not possible to prove such a result unless something else comes into play. Often one considers first what happens for spherically symmetric initial data; we choose to use a bit of spherical symmetry. We also know that decay is needed when we deal with quadratic (D-NLS). 

We define the differential operator: 
\beq \label{r}
\mathcal{R}f=(x_{1} \partial_{x_{2}} - x_{2} \partial_{x_{1}})f 
\eeq

\noindent
and the pseudodifferential operators in the left calculus:

\beq \label{de}
\mathcal{D}f= (1+\frac{\langle x \rangle^{2}}{\mu+\langle D \rangle^{2}})^{\frac{1}{4}+\frac{\e}{2}}f \ \ \mbox{with symbol} \ \ (1+\frac{\langle x \rangle^{2}}{\mu+\langle(\xi,\tau) \rangle^{2}})^{\frac{1}{4}+\frac{\e}{2}}
\eeq

\noindent 
for some $0 < \e < \q$. For a generic space of functions $\mathcal{X}$ we define:

\beq \label{S}
\mathcal{D}\mathcal{R}\mathcal{X}= \{f \in \mathcal{X}: \mathcal{D}f\in \mathcal{X} \ \mbox{and} \  \mathcal{D}\mathcal{R}f \in \mathcal{X}\}
\eeq

We denote  by $\chi_{[0,T]}$ a smooth approximation of the characteristic function of $[0,T]$ such that $\chi_{[0,T]}(t)=1, \ \forall t \in [0,T]$. We will always consider $\chi_{[0,T]}$ as a function of time, in other words by $\chi_{[0,T]}$ we mean  $\chi_{[0,T]}(t)$.

The space for our initial data is $\mathcal{D}\mathcal{R}H^{s}$. We dedicate section \ref{ds} to the definition of the spaces  $\mathcal{D}\mathcal{R}Z^{s,5}$ (for the solutions) and $\mathcal{D}\mathcal{R}W^{s}$ (for the inhomogeneity). These spaces satisfy the linear estimate:

\begin{t1} \label{le1}

If $g \in \mathcal{D}\mathcal{R}H^{s}$ and $f \in \mathcal{D}\mathcal{R}W^{s}$, then the solution of:

\begin{equation} \label{E}
\begin{cases}
\begin{aligned} 
&iu_{t}-\Delta u = f \\
&u(x,0) =g(x)
\end{aligned}
\end{cases}
\end{equation}

satisfies  $\chi_{[0,1]} u \in \mathcal{D}\mathcal{R}Z^{s,5} \cap C_{t}\mathcal{D}\mathcal{R}H_{x}^{s}$.

\end{t1}

To each quadratic nonlinearity we associate is the standard way the bilinear form $B(u,v)$. The bilinear estimate is the next key result:

\begin{t1} \label{bg}
If $s > s_{c}$, we have the global bilinear estimate:

\beq \label{bg1}
||B(u,v)||_{\mathcal{R}\mathcal{D}W^{s}} \leq C_{\e,s} ||u||_{\mathcal{R}\mathcal{D}Z^{s}}||v||_{\mathcal{R}\mathcal{D}Z^{s}} 
\eeq

\end{t1}

Once we have the above two results, a standard fixed point argument gives us the main result:

\begin{t1} \label{T}

Assume n=2. Given any $s > s_{c}$ and $T > 0$, there exists $\delta > 0$ such that for every $u_{0} \in \mathcal{D}\mathcal{R}H^{s}$ with $\delta_{0}=||u_{0}||_{\mathcal{D}\mathcal{R}H^{s}}  < \delta$ , the quadratic (D-NLS) has a unique solution $u$ in $C([0,T]:\mathcal{D}\mathcal{R}H^{s}) \cap \mathcal{R}\mathcal{D} Z^{s,5}$ with Lipschtiz dependence on the initial data.

\end{t1}

A sketch of the proof goes as follows. Let $B(u,v)$ be the bilinear form:

\beq \label{bf}
B(u,v)=\sum_{i,j \in \{1,2\}} c_{ij} u_{x_{i}} v_{x_{j}}
\eeq

\noindent
where $c_{ij}$ are constant complex numbers. We intend to obtain bilinear estimates for $B(u,v)$ and $B(u,\bar{v})$ since this way we cover the theory for all quadratic polynomials of type $P(\nabla{u},\nabla{\bar{u}})$, except for those of type $P(\nabla{\bar{u}})$. For the last ones the theory had been developed previously, as we remarked before.  

 We start with the Bourgain space $X^{s,\q,1}$ as the candidate for $Z^{s}$ and $X^{s,-\q,1}$ as a candidate for $W^{s}$. We split the Fourier space in pieces according to the size of $(\xi,\tau)$ and its distance to the paraboloid ($\tau=\xi^{2}$). Taking into account that the product becomes convolution under the Fourier transform, we see how the pieces interact and try to recover the $X^{s,-\q,1}$ structure for $B(u,v)$. This goes fine as long as we recover information which is at some distance from the paraboloid and it breaks down very close to the paraboloid - we catch a logarithm of the high frequency which cannot be controlled. To remedy this we come up with a more delicate decomposition of the part of the Fourier space which is at distance less than $1$ from the paraboloid. More exactly we introduce a wave packet decomposition and we measure the packets in $L^{\infty}_{t}L^{2}_{x}$. Then the target space $W^{s}$ is also modified at distance less than $1$ from paraboloid, i.e. we also have a wave packet decomposition and the packets are measured in $L^{1}_{t}L^{2}_{x}$. We have to recover a $L^{1}_{t}$ structure on the packets for $B(u,v)$ and this is why we need to involve the extra decay and spherical symmetry.

All along the argument we do involve some spherical symmetry and decay in the bilinear estimates and this is why our spaces will be of type $\mathcal{R}\mathcal{D}Z^{s}$ and $\mathcal{R}\mathcal{D}W^{s}$. See the next section for the definitions.

Once the bilinear estimates are fixed, then a standard fixed point argument gives us the result of Theorem \ref{T}. 

One can easily adapt our argument for the bilinear forms of type:

\beq \label{bf1}
B(u,v)=\sum_{j=1}^{2} c_{j} u v_{x_{j}}
\eeq

 This is because the basic estimates are derived for the bilinear form $\tilde{B}(u,v)=u \cdot v$ and then we "over-estimate" the size of $\nabla$, see the beginning of section \ref{best} for more details. Thus we are entitled to claim the result for the quadratic polynomials of type $P(u,\nabla{u})$,  $P(\bar{u},\nabla{u})$ and $P(u,\nabla{\bar{u}})$.

We think that further analysis should reveal that without assuming some spherical symmetry there is no way one can get existence all the way down to the critical exponent. Without spherical symmetry, but involving decay, we expect a positive result for $s > s_{c}+1$ and a negative one for $s < s_{c}+1$. This work is in progress and it will be the main core of the second part of this paper. 
 
The spaces we use in this paper are in some way the counterpart of the ones involved in dealing with the wave maps equation, see \cite{ta1} and \cite{tao1}. Our spaces are a bit more difficult since they involve phase-space localization, rather than phase localization which is the case for the wave-maps. 

We conclude the introduction with few open problems. The most obvious thing to ask is what happens at the scaling exponent $s_{c}$. Our techniques lose logarithms of the low frequency in the bilinear estimates and we can eliminate them only by imposing $s > s_{c}$. 
Then the question of the optimality of the decay is another one to ask. Is the factor $\q+\e$ the optimal one? The estimates in Bourgain spaces seem to indicate that we should manage with a factor of $\e$.

In the end the generalization to higher dimensions should be of interest too. We know that the scaling exponent is $\frac{n}{2}$ for the case when both terms come with derivatives and we think it should be possible to get similar results under similar conditions in all dimensions.

Quadratic D-NLS is in itself an interesting problem, but one of the most important reasons to study it comes from the Schr\"odinger maps. They are the natural Schr\"odinger equation when the target space is a complex manifold. It is well-known in the literature that understanding the problem (\ref{E:s}) with nonlinearity $|\nabla{u}|^{2}$ is essential for the study of the Schroedinger-maps equation (at least for the case when the manifold is $\mathbb{S}^{n-1}$). So far the problem has been solved for initial data in $H^{\frac{3}{2}+\e}$ in $\R^{2}$ and when the target manifold is $\mathbb{S}^{2}$ or $\mathbb{H}^{2}$, see \cite{n1}. It is known that the scaling exponent for these problems is $\frac{n}{2}$, so a natural question to ask in dimension $2$ is what happens when the initial data is in $H^{s}$ for $1 \leq s \leq \frac{3}{2}$. The result in higher dimensions is also open. The results in this paper might be a good start in understanding and approaching these problems. 

\vspace{.1in}

The results of this paper were mostly obtained when the author was a PhD student at University of California, Berkeley. The author is greatly indebted to Daniel Tataru, his thesis adviser, for many fruitful conversations and for the constant
encouragement along the way.

\vspace{.2in}

\section{Proof of Main Theorem}

Assuming the results of Theorems \ref{le1} and \ref{bg} we can prove the result of the Theorem \ref{T}. We did not define yet the spaces $Z^{s,5}$ and $W^{s,5}$, but at this time it is enough to take for granted that they are Banach spaces. 

We define the operator $\mathcal{T}_{1}$ by $w=\mathcal{T}_{1}f$ to be the solution of the inhomogeneous  Schr\"odinger equation with zero initial data: 

\begin{equation} \label{e2}
\begin{cases}
\begin{aligned} 
&iw_{t}-\Delta w = f \\
&w(x,0) =0
\end{aligned}
\end{cases}
\end{equation}

We fix $T=1$ and prove that if the initial data is small enough then our problem has a solution. We define the set $K$ to be 

$$K=\{\chi_{[0,1]}w \in \mathcal{R}\mathcal{D}Z^{s,5}: ||\chi_{[0,1]}w||_{\mathcal{R}\mathcal{D}Z^{s,5}} \leq 2||\chi_{[0,1]} e^{it\Delta} u_{0}||_{\mathcal{R}\mathcal{D}Z^{s,5}} \}$$

and the operator $\mathcal{T}: \mathcal{R}\mathcal{D}Z^{s,5} \rightarrow \mathcal{R}\mathcal{D}Z^{s,5}$ by 

$$\mathcal{T} (v) = e^{it\Delta} u_{0} + \mathcal{T}_{1} ( \chi_{[0,1]}^{2} B(v,v))$$

In the hypothesis that $||u_{0}||_{\mathcal{D}\mathcal{R}H^{s}}$ is small enough, we prove that  $\mathcal{T}:K \rightarrow K$ and that $\mathcal{T}$ is a contraction on $K$. This give us the existence of a fixed point for $\mathcal{T}$ which is the solution of our problem in the interval $[0,1]$. This is because we chose $\chi_{[0,1]}$ to be equal to $1$ on $[0,1]$.  

To prove the invariance of $K$ under the action of $\mathcal{T}$ we use (\ref{e3}) and (\ref{e4}):

$$||\chi_{[0,1]}\mathcal{T}u||_{\mathcal{D}\mathcal{R}Z^{s,5}} \leq ||\chi_{[0,1]} e^{it\Delta}u_{0}||_{\mathcal{D}\mathcal{R}Z^{s,5}}+||\chi_{[0,1]} \mathcal{T}_{1}(\chi_{[0,1]}^{2} B(u,u))||_{\mathcal{D}\mathcal{R}Z^{s,5}} \leq $$

$$||\chi_{[0,1]} e^{it\Delta}u_{0}||_{\mathcal{D}\mathcal{R}Z^{s,5}} + C_{5}||\chi_{[0,1]} B(\chi_{[0,1]}u,\chi_{[0,1]}u)||_{\mathcal{D}\mathcal{R}W^{s,5}}$$

Using the bilinear estimate in (\ref{bg1}) we continue with:

$$||\chi_{[0,1]}\mathcal{T}u||_{\mathcal{D}\mathcal{R}Z^{s,5}} \leq ||\chi_{[0,1]} e^{it\Delta}u_{0}||_{\mathcal{D}\mathcal{R}Z^{s,5}}+ C_{\e,s} ||\chi_{[0,1]} u||_{\mathcal{R}\mathcal{D}Z^{s,5}}^{2} \leq $$

$$||\chi_{[0,1]} e^{it\Delta}u_{0}||_{\mathcal{D}\mathcal{R}Z^{s,5}}+ C_{\e,s}||\chi_{[0,1]} e^{it\Delta}u_{0}||^{2}_{\mathcal{R}\mathcal{D}Z^{s,5}}  $$

It is enough to choose $||u_{0}||_{\mathcal{D}\mathcal{R}H^{s}}$ small enough in order to obtain the bound $C_{\e,s}||e^{it\Delta}u_{0}||_{\mathcal{D}\mathcal{R}Z^{s}} \leq C_{1} C_{\e,s}||u_{0}||_{\mathcal{D}\mathcal{R}H^{s}} \leq 1$ which gives us the desired inequality. 

To prove that $\mathcal{T}$ is a contraction we proceed as follows:

$$\mathcal{T}u_{1}-\mathcal{T}u_{2} = \mathcal{T}_{1}(\chi_{[0,1]}^{2} B(u_{1},u_{1}))-\mathcal{T}_{1}(\chi_{[0,1]}^{2} B(u_{2},u_{2}))=$$ 

$$\mathcal{T}_{1}(\chi_{[0,1]}^{2} B(u_{1},u_{1}-u_{2}))+\mathcal{T}_{1}(\chi_{[0,1]}^{2} B(u_{2},u_{1}-u_{2}))$$

followed by the estimates

$$||\chi_{[0,1]}(\mathcal{T}u_{1}-\mathcal{T}u_{2})||_{\mathcal{D}\mathcal{R}Z^{s,5}} \leq $$

$$ C_{\e,s}(||\chi_{[0,1]} u_{1}||_{\mathcal{D}\mathcal{R}Z^{s,5}}+||\chi_{[0,1]} u_{2}||_{\mathcal{D}\mathcal{R}Z^{s,5}})||\chi_{[0,1]}^{2}(u_{1}-u_{2})||_{\mathcal{D}\mathcal{R}Z^{s,5}} \leq$$

$$C_{\e,s}||u_{0}||_{\mathcal{D}\mathcal{R}H^{s}} ||\chi_{[0,1]}(u_{1}-u_{2})||_{\mathcal{D}\mathcal{R}Z^{s,5}} < \q ||\chi_{[0,1]}(u_{1}-u_{2})||_{\mathcal{D}\mathcal{R}Z^{s,5}}$$

\noindent
where again we have to choose $||u_{0}||_{\mathcal{D}\mathcal{R}H^{s}}$ small enough so that we have $C_{\e,s}||u_{0}||_{\mathcal{D}\mathcal{R}H^{s}} < \q$. 

We conclude that $\mathcal{T}$ has a unique fixed point in $K$ which is a solution to our problem.  By rescaling we can obtain the result of the Theorem for any $T > 0$.

\vspace{.2in}

\section{Definition of the spaces} \label{ds}

For each $u$ we denote by $\mathcal{F}u=\hat{u}$ the Fourier transform of $u$. This is always taken with respect to all the variables, unless otherwise specified.

Throughout this paper all the inequalities of type $\leq$ should be understood in the sense $\lesssim$: $A \leq B$ is meant to be $A \lesssim B \Leftrightarrow A \leq CB$ for some constant $C$ which is independent of any possible variable in our problem.   

We say $A \approx B$ if $A \leq CB \leq C^{2}A$ for the same constant $C$. We say that we localize at frequency $2^{i}$ to mean that  in the support of the localized function $|(\xi,\tau)| \in [2^{i-1},2^{i+1}]$.

The paraboloid $P=\{(\xi,\tau): \tau=\xi^{2}\}$ plays a very important role in the geometry of the problem.

In the Schr\"odinger equation time and space scale in a different way, and this suggests to define the norm for $(\xi, \tau)$ by $|(\xi, \tau)|=(|\tau|+ \xi^{2})^{\q}$. In dealing with the quadratic nonlinearity without derivatives the Bourgain space $X^{s,b}$ proved to be a very useful space to work with for appropriate choice of $b$, see \cite{c1}. They are defined in the following way:

$$X^{s,b} =  \{f\in S'; \langle (\xi,\tau) \rangle^{s} \langle \tau-\xi^{2} \rangle^{b} \hat{f} \in L^{2} \}$$

Here and thereafter $\langle x \rangle=(1+|x|^{2})^{\q}$ where $|x|$ is the norm of $x$. The integral defining $X^{s,\q}$ has two weights in it, an elliptic one, $ \langle (\tau,\xi) \rangle^{s}$, and one adapted to the paraboloid, $ \langle \tau-\xi^{2} \rangle^{\q}$. We will employ frequency localized versions of $X^{s,\q}$ which are constructed according to these weights.

Consider $\varphi_{0} : [0,\infty) \rightarrow \R$ to be a nonnegative smooth function such that $\varphi_{0}(x)=1$ on $[0,1]$ and $\varphi_{0}(x)=0$ if $x \geq 2$. Then for each $i \geq 1$ we define $\varphi_{i}: [0,\infty) \rightarrow \R$ by $\varphi_{i}(x)= \varphi_{0}(2^{-i}x)- \varphi_{0}(2^{-i+1}x)$. With the help of  $(\varphi_{i})_{i \geq 0}$ we define the operators $S_{i}$, $S^{\xi}_{i}$ and $S_{i}^{\tau}$ by:

$$\mathcal{F}(S_{i}f)= \hat{f}_{i} = \varphi_{i}(|(\xi, \tau)|) \cdot \hat{f}(\xi, \tau)$$
 
$$\mathcal{F}(S^{\xi}_{i}f)= \hat{f}_{i}^{\xi} = \varphi_{i}(|\xi|) \cdot \hat{f}(\xi, \tau)$$

$$\mathcal{F}(S_{i}^{\tau}f)= \hat{f}_{i}^{\tau} = \varphi_{i}( \sqrt{|\tau|}) \cdot \hat{f}(\xi, \tau)$$

Since $|\tau-\xi^{2}| \approx |(\tau,\xi)| d((\xi,\tau),P)$ (away from zero), then we can chose to localize with respect to $d((\xi,\tau),P)$ instead of $|\tau-\xi^{2}|$. If $|(\xi,\tau)| \approx 2^{i}$, then $|\tau-\xi^{2}|$ ranges in the interval $[0,2^{2i+2}]$, hence $d((\tau,\xi),P)$ ranges in the interval $[0,2^{i+2}]$. The appropriate localization for $|\tau-\xi^{2}|$ is on a dyadic scale of type $2^{k}$ with $k \in \{0,1,..,2i+2\}$. Therefore the appropriate localization for $d((\tau,\xi),P)$ is on a dyadic scale of type $2^{k}$ with $k \in \{-i,-i+1,..,-1,0,1,..,i+2\}$. For $d \in I_{i}=\{2^{-i}, 2^{-i+1}, .., 2^{i+1}\}$, we build a system of functions $\varphi_{i,d} : \R^{2} \rightarrow \R$  having the following property that the support of $\varphi_{i,d}$ is approximately the set 

$$\{(\xi,\tau): |(\xi,\tau)| \approx 2^{i}, d((\tau,\xi),P) \approx d\} \approx \{(\xi,\tau): |(\xi,\tau)| \approx 2^{i}, |\tau -\xi^{2}| \approx d2^{i}\}$$

\noindent
and

$$ \sum_{d \in I_{i}} \varphi_{i,d} (\xi,\tau) = \varphi_{i}(|(\xi,\tau)|)  , \ \forall (\xi,\tau) \in \R^{2} \times \R $$

 The support of $\varphi_{i,2^{-i}}$ should contain not only the points which are at distance $\approx 2^{-i}$ from $P$, but also those at distance less than $2^{-i}$ from $P$.  

We define the operators $S_{i,d}$ by $S_{i,d}f=f_{i,d}=\check{\varphi}_{i,d} * S_{i}f$ and we have $f_{i} = \sum_{d \in I_{i}} f_{i,d}$. In the support of $\hat{f}_{i,d}$ we have $1+|\tau - \xi^{2}| \approx 2^{i}d$.

Sometimes it is useful to localize in a linear way rather than a dyadic way. In these cases we localize with respect to the value of $|\tau - \xi^{2}|$ instead; we will make this clear when we need it.

For each dyadic value $d \in I_{i}$ we define $\varphi_{i,\leq d}=\sum_{d' \in I_{i}: d' \leq d} \varphi_{i,d'}$ and $\varphi_{i,\geq d}=\sum_{d' \in I_{i}: d' \geq d} \varphi_{i,d'}$. The give rise to the operators which localize at distance less and greater than $d$ from $P$: 

$$S_{i, \leq d}f = f_{i, \leq d}= f * \check{\varphi}_{i, \leq d} \ \  \mbox{and} \ \ S_{i, \geq d}f=f_{i, \geq d}= f * \check{\varphi}_{i, \geq d}$$

The part of $\hat{f}$ which is at distance less than $1$ from $P$ plays an important role and this is why we define the global operators:

$$S_{\cdot,\leq 1}f=f_{\cdot, \leq 1}=\sum_{i=0}^{\infty} f_{i, \leq 1} \ \  \mbox{and} \ \ S_{\cdot,\geq 1}f=f_{\cdot,\geq 1}=\sum_{i=0}^{\infty} f_{i, \geq 1}$$

We denote by $A_{i}$ the support in $\R^{2} \times \R$ of $\varphi_{i}(|(\xi,\tau)|)$ and by $A_{i,d}$ the support of $\varphi_{i,d}$. In a similar way we can define $A_{i, \leq d}$ and $A_{i, \geq d}$ to be the support of the operators $S_{i, \leq d}$, respectively $S_{i, \geq d}$.

We work with $X^{s,\q, 1}$ which is defined as follows:

$$ X^{s,\q, 1}_{i} = \{f: \hat{f} \ \mbox{supported in} \ A_{i} \ \mbox{and} \ ||f||_{X^{s,\q, 1}_{i}}=\sum_{d \in I_{i}} ||f_{i,d}||_{X^{s,\q}} \leq \infty \}$$

$$ X^{s,\q, 1}= \{f: f_{i} \in X_{i}^{s,\q,1} \ \mbox{and} \ ||f||_{X^{s,\q, 1}}^{2}= \sum_{i} ||f_{i}||^{2}_{X^{s,\q,1}_{i}} < \infty \}$$

For technical purposes we need:

$$ X^{s,\q}_{i,d} = \{f \in X^{s,\q}: \hat{f} \ \mbox{supported in} \ A_{i,d} \}$$

\noindent 
and, similarly, $ X^{s,\q,1}_{i,\leq d}$ and $ X^{s,\q,1}_{i,\geq d}$.

When we work out the estimates it turns out that $X^{s,\q,1}$ is the right space to measure only the part of the solution whose support in the Fourier space is at distance greater than $1$ from $P$, i.e. the $S_{\cdot,\geq 1}$ part of our solutions. 
 
This is why we introduce also:

$$X^{s,\q,\infty}_{i}=\{f:  \hat{f} \ \mbox{supported in} \ A_{i} \ \mbox{and} \ ||f||_{X^{s,\q, \infty}_{i}}=||f_{i,d}||_{l^{\infty}_{d}(X^{s,\q})} < \infty, d \in I_{i} \} $$

$$X^{s,\q,\infty}= \{f: f_{i} \in X^{s,\q,\infty}_{i}  \mbox{and} \ ||f||^{2}_{X^{s,\q, \infty}}= \sum_{i} ||f_{i}||^{2}_{X^{s,\q,\infty}_{i}} < \infty \}$$

We will measure the $S_{\cdot, \leq 1}$ part of our solutions in $X^{s,\q,\infty}$; in addition to that we will measure it in a space whose construction goes as follows. 

We define the following lattice in the plane $\tau=0$:

$$\Xi = \{\xi=(r,\t): r=n,\  \t=\frac{\pi}{2} \frac{k}{n}, \ n,k \ \mbox{positive integers} \}$$

$\Xi$ is like a lattice in polar coordinates. It has the properties that the distance between any two points is at least $1$ and that for every $\eta \in \R^{2}$ there is a $\xi \in \Xi$ such that $|\xi-\eta| \leq 1$. For each $\xi \in \Xi$ we build a non-negative function $\phi_{\xi}$ to be a smooth approximation of the characteristic function of the cube of size $1$ in $\R^{2}$ centered at $\xi$ and satisfying the natural partition property:

\beq \label{d100}
\sum_{\xi \in \Xi} \phi_{\xi}=1
\eeq

We can easily impose uniforms bounds on the derivatives of the system $(\phi_{\xi})_{\xi \in \Xi}$. For each $\xi \in \Xi$ we define:

$$f_{\xi} = \check{\phi}_{\xi} * f \ \ \ \ \mbox{and} \ \ \ f_{\xi, \leq 1} = \check{\phi}_{\xi} * f_{\cdot,\leq 1}$$

The convolution above is performed only with respect to the $x$ variable. The support of $\hat{f}_{\xi, \leq 1}$ is a almost a parallelepiped having the center $(\xi,\xi^{2}) \in P$ and sizes: $\approx |\xi|$ in the $\tau$ direction and $1$ in the other two directions (normal to $P$ and the completing third one). 

We can also build a system of non-negative functions $\tilde{\phi}_{\xi}$ satisfying: 

- $\tilde{\phi}_{\xi} \cdot \phi_{\xi}=\phi_{\xi}$ and $\sum_{\xi} \tilde{\phi}_{\xi} \leq C$

- the support of $\tilde{\phi}_{\xi}$ is contained in the set $\{\eta: |\eta-\xi| \leq 2\}$

- the system $(\tilde{\phi}_{\xi})_{\xi \in \Xi}$ has uniform bounds on the derivatives. 

For technical purposes we need the following construction. For $n$ positive integer define $\Xi_{n} = \{ \xi \in \Xi: |\xi|=n \}$ and then :

$$f_{n , \leq 1} = \sum_{\xi \in \Xi_{n}} f_{\xi, \leq 1}$$

The next concern is how to measure $f_{\xi, \leq 1}$. We denote by $(Q^{m})_{m \in \Z^{2}}$ the standard partition of $\R^{2}$ in cubes of size $1$; i.e. $Q^{m}$ is centered at $m=(m_{1},m_{2}) \in \Z^{2}$, has its sides parallel to the standard coordinate axis and has size $1$. For each $\xi \in \Xi$, $m \in \Z^{2}$ and $l \in \Z$ we define the tubes:

$$T_{\xi}^{m,l}=\cup_{t \in [l,l+1]} (Q^{m}-2t\xi) \times \{ t \}=$$

$$\{(x-2t\xi_{1}, y-2t\xi_{2},t): (x,y) \in Q^{m} \ \mbox{and} \ t \in [l,l+1] \}$$

Then, for each $\xi \in \Xi$, we define the space $Y_{\xi}$ by the following norm:

$$||f||^{2}_{Y_{\xi}}= \sum_{(m,l) \in \Z^{3}} ||f||^{2}_{L^{\infty}_{t}L^{2}_{x}(T_{\xi}^{m,l})}$$

We have $f=\sum_{\xi \in \Xi} f_{\xi}$ and then we define the space $Y^{s}$ by the norm:

$$||f||^{2}_{Y^{s}}= \sum_{\xi \in \Xi} \langle \xi \rangle^{2s} ||f_{\xi}||^{2}_{Y_{\xi}}$$

For technical reasons we need also:

$$Y_{i}= \{ f \in Y^{0}; \hat{f} \ \mbox{supported in} \ A_{i} \}$$

$$Y_{i, \leq d}= \{ f \in Y^{0}; \hat{f} \ \mbox{supported in} \ A_{i, \leq d} \}$$

\noindent
the last one being defined for any $d \in I_{i}$ with $d \leq 1$.

We localize our solutions in time. If we come with a frequency localization on the top of this we are left with decay in time of our solutions. For this we define $Y_{\xi}^{N}$ and $Y^{s,N}$ by the norms:

$$||f||_{Y_{\xi}^{N}}=||\langle t \rangle^{N}f||_{Y_{\xi}} \ \ \ \mbox{and} \ \ \ ||f||^{2}_{Y^{s,N}}= \sum_{\xi \in \Xi} \langle \xi \rangle^{2s} ||f_{\xi}||^{2}_{Y_{\xi}^{N}}$$

To bring everything together, define $Z^{s,N}$ to be 

$$Z^{s,N}=\{ f \in S': ||f_{\cdot,\geq 1}||_{X^{s,\q,1}}+||f_{\cdot, \leq 1}||_{Y^{s,N}}+||f_{\cdot, \leq 1}||_{X^{s,\q,\infty}} < \infty \}$$

\noindent
with the obvious norm. Our spaces are going to be equipped with some additional structure, namely a bit of spherical symmetry and some decay.  Recalling the definitions in (\ref{r}), (\ref{de}) and (\ref{S}), we are going to measure our solution in $\mathcal{D}\mathcal{R}Z^{s,5}$.

So far we have built the spaces suitable for the solution of (\ref{E:s}). We need also a space for the right hand side of the equation, see Theorem \ref{le1}.

We can easily define $X^{s,-\q,1}$ by simply replacing $\q$ with $-\q$ in the definition of $X^{s,\q,1}$. Then we define $\mathcal{Y}^{s}$ by:

$$||f||^{2}_{\mathcal{Y}^{s,N}}= \sum_{\xi \in \Xi} \langle \xi \rangle ^{2s} ||f_{\xi}||^{2}_{\mathcal{Y}_{\xi}^{N}}$$

\noindent
where $\mathcal{Y}_{\xi}^{N}$ is defined as follows:

$$||f||^{2}_{\mathcal{Y}_{\xi}^{N}}= \sum_{(m,l) \in Z^{3}} || \langle t \rangle^{N} f||^{2}_{L^{1}_{t}L^{2}_{x}(T_{\xi}^{m,l})}$$

Notice that $(\mathcal{Y}_{\xi})^{*}=Y_{\xi}$ since we will use this later for duality purposes.

Back to the right hand side of ($\ref{E}$), we need $f_{\cdot, \geq 1}$ to be in $X^{s,-\q,1}$, but we do not need  $f_{\cdot, \leq 1}$ to be both in $\mathcal{Y}^{s,5}$ and $X^{s,-\q,1}$ in order to recover the $Z^{s,5}$ structure for $\chi_{[0,1]}u$. This is why we introduce $W^{s,N}$ defined by the norm:

$$||f||_{\mathcal{W}^{s,N}}= \inf \{||f_{1}||_{\mathcal{Y}^{s,N}} + ||f_{2}||_{X^{s,-\q,1}}; f=f_{1}+f_{2} \}$$

We measure the right hand side of ($\ref{E}$) in: 

$$W^{s,5} =\{ f \in S': ||f_{\cdot,\leq 1}||^{2}_{\mathcal{W}^{s,5}} + ||f_{\cdot,\geq 1}||^{2}_{X^{s,-\q,1}} < \infty \}$$

Besides $X^{s,b}$ we need the conjugate $\bar{X}^{s,b}$ which is defined as follows:

$$\bar{X}^{s,b} =  \{f\in S'; \langle (\xi,\tau) \rangle^{s} \langle \tau+\xi^{2} \rangle^{b} \hat{f} \in L^{2}  \}$$

We can define all the other elements the same way as above by simply placing a bar on each space and operator, while replacing everywhere $|\tau-\xi^{2}|$ with $|\tau+\xi^{2}|$ and $P$ with $\bar{P}=\{(\xi,\tau): \tau+\xi^{2}=0 \}$.  

We have the following important fact:

$$ f \in X^{s,b} \Longleftrightarrow \bar{f} \in \bar{X}^{s,b}$$

\noindent
and the obvious correspondents for the variants of $X^{s,\q}$ we work with. In addition the dual of $X^{s,b}$ is:

$$(X^{s,b})^{*}=\bar{X}^{-s,-b}$$

\section{The linear estimates}

This section is dedicated to proving the result in Theorem \ref{le1}. As expected from the statement we have to prove the following two Propositions.

\begin{p2} \label{k4}

 The solution of the homogeneous equation satisfies:

\beq \label{e3}
||\chi_{[0,1]} e^{it\D}g||_{\mathcal{D}\mathcal{R}Z^{s,N} \cap C_{t} \mathcal{D}\mathcal{R} H^{s}_{x}} \leq C_{N} ||g||_{\mathcal{D}\mathcal{R}H^{s}}
\eeq

\end{p2}

\begin{p2} \label{k3}

The solution of (\ref{e2}) satisfies:

\beq \label{e4}
||\chi_{[0,1]} w||_{\mathcal{D}\mathcal{R}Z^{s,N} \cap C_{t}\mathcal{D}\mathcal{R} H^{s}_{x}} \leq C_{N} ||f||_{\mathcal{D}\mathcal{R}W^{s}}
\eeq

\end{p2}

 It is well-known that the Schr\"odinger equation is invariant under rotations, therefore  it is enough to prove the results in the two Proposition without rotations. As about conserving the decay, we leave this problem for the end of the section. Therefore we first prove Propositions $\ref{k4}$ and $\ref{k3}$ without involving the $\mathcal{R} \mathcal{D}$ structure. 

The following result brings some important informations about the structure of the spaces we work with and will be very useful for the rest of the section.

\begin{l1} \label{ees} The spaces we involve have the following properties:

a) $X^{s,\q,1} \subset C_{t}H^{s}_{x}$ 

b) $f \in X^{s,\q,1} \Ra \chi_{[0,1]}f \in X^{s,\q,1}$

c) $\mathcal{Y}^{0} \subset X^{0,-\q,\infty}$ and $X^{0,\q,1} \subset Y^{0}$

d) $f \in Y^{s} \Ra \chi_{[0,1]}f \in Y^{s,N} \ \forall N$

e) $f \in Y^{s,N} \Ra f_{\cdot, \leq 1} \in Y^{s,N}$ and, more generally, for any $d \in I_{i}$ we have $f_{i} \in Y^{s,N} \Ra f_{i, \leq d} \in Y^{s,N}$.  

\end{l1}

\begin{proof}

a) and b) are standard result, see for instance \cite{ta2}.

c) Let us fix $i$, $d \in I_{i}$ and $\xi \in \Xi, |\xi| \approx 2^{i}$. Let us denote by $\tilde{\phi}_{\xi,d}=\tilde{\phi}_{\xi} \cdot \varphi_{i,d}$. We have $(f_{\xi})_{i,d}=f_{\xi} * \check{\tilde{\phi}}_{\xi,d}$. One can show (see part e)) that $\check{\tilde{\phi}}_{\xi,d}$ is highly concentrated in $T_{\xi}^{0,0}$ in the following sense:

$$||\chi_{T^{m,l}_{\xi}}\check{\tilde{\phi}}_{\xi,d}||_{L^{2}} \leq C_{N} \langle (m,l) \rangle^{-N} ||\check{\tilde{\phi}}_{\xi,d}||_{L^{2}} \approx C_{N} \langle (m,l) \rangle^{-N} (2^{i}d)^{\q}$$ 

This allows us to conclude that:

$$||(f_{\xi})_{i,d}||^{2}_{L^{2}} \leq \sum_{m,l} ||\chi_{T^{m,l}_{\xi}}f_{\xi} * \check{\tilde{\phi}}_{\xi,d}||^{2}_{L^{2}} \leq \sum_{m,l} ||\chi_{T^{m,l}_{\xi}}f_{\xi} * \check{\tilde{\phi}}_{\xi,d}||^{2}_{L^{2}_{t}L^{\infty}_{x}} \leq$$

$$\sum_{m,l} ||\chi_{T^{m,l}_{\xi}}f_{\xi}||^{2}_{L^{1}_{t}L^{2}_{x}} ||\check{\tilde{\phi}}_{\xi,d}||^{2}_{L^{2}} \approx (2^{i}d)^{\q} ||f_{\xi}||^{2}_{Y_{\xi}}$$

Summing up with respect to the $\xi \in \Xi$, $|\xi| \approx 2^{i}$ gives us:

$$||f_{i,d}||_{L^{2}} \leq (2^{i}d)^{\q}||f_{i}||_{Y^{0}}$$

\noindent
which is enough to conclude $\mathcal{Y}^{0} \subset X^{0,-\q,\infty}$. By duality we obtain $X^{0,\q,1} \subset Y^{0}$.

d) We observe that $(\chi_{[0,1]}f)_{\xi}=\chi_{[0,1]}f_{\xi}$. Since $f_{\xi} \in Y_{\xi}$, then it follows immediately that $t^{N} \chi_{[0,1]}f_{\xi} \in Y_{\xi}$ for any $N$. Summing this up with respect to $\xi \in \Xi$ gives us the claim.

e) We fix $i$, $d \in I_{i}$ and $\xi \in \Xi, |\xi| \approx 2^{i}$ and denote by $\tilde{\phi}_{\xi,\leq d}=\tilde{\phi}_{\xi} \cdot \varphi_{i,\leq d}$. We observe that $(f_{\xi})_{i, \leq d}=f_{\xi} * \check{\tilde{\phi}}_{\xi,\leq d}$. $\phi_{\xi,\leq d}$ is a smooth approximation of the characteristic function of the set $\{(\eta,\tau): |\eta-\xi| \leq \q, |\tau-\xi^{2}| \leq 2^{i}d \}$ which, geometrically, is approximately a parallelepiped. Then $\check{\tilde{\phi}}_{\xi,\leq d}$ is highly localized in a subset of $T^{0,0}_{\xi}$, namely $\cup_{t \in [0,(2^{i}d)^{-1}]} (Q^{0}-2t\xi) \times \{ t \}$ the dual parallelepiped. In the particular setup that $d=2^{-i}$, we can quantify this high localization as follows:

$$||\check{\tilde{\phi}}_{\xi,\leq 2^{-i}}||_{L^{2}(T^{m,l}_{\xi})} \leq C_{N} \langle (m,l) \rangle^{-N} ||\check{\tilde{\phi}}_{\xi,\leq 2^{-i}}||_{L^{2}} = C_{N} \langle (m,l) \rangle^{-N} ||\tilde{\phi}_{\xi,\leq 2^{-i}}||_{L^{2}}$$

In the general case we should construct a new family of tubes by splitting each $T^{m,l}_{\xi}$ in $2^{i}d$ subtubes, by splitting the time interval in equal intervals, and then have a similar estimate on these tubes. We skip this formalization, since it does not bring anything illuminating and requires complicated notations. One key observation, besides the fact that $\check{\tilde{\phi}}_{\xi,\leq d}$ is highly localized in $\cup_{t \in [0,(2^{i}d)^{-1}]} (Q^{0}-2t\xi) \times \{ t \}$ is that:

$$||\check{\tilde{\phi}}_{\xi,\leq d}||_{L^{1}(\cup_{t \in [0,(2^{i}d)^{-1}]} (Q^{0}-2t\xi) \times \{ t \})} \leq (2^{i}d)^{-\q} ||\check{\tilde{\phi}}_{\xi,\leq d}||_{L^{2}}= (2^{i}d)^{-\q} ||\tilde{\phi}_{\xi,\leq d}||_{L^{2}}\approx 1$$

On behalf of all these facts, a straightforward argument gives us:

$$||\check{\tilde{\phi}}_{\xi,\leq d}||_{L^{1}(T^{m,l}_{\xi})}  \leq C_{N} \langle (m,l) \rangle^{-N}$$

This allows us to estimate:

$$||f_{\xi} * \check{\tilde{\phi}}_{\xi,\leq d}||^{2}_{Y_{\xi}} \leq \sum_{m,l} \left( \sum_{m',l'} ||\chi_{T^{m,l}_{\xi}} f_{\xi} * \chi_{T^{m',l'}_{\xi}} \check{\tilde{\phi}}_{\xi,\leq d}||_{L^{\infty}_{t}L^{2}_{x}} \right)^{2} \leq $$

$$\sum_{m,l} ||\chi_{T^{m,l}_{\xi}} f_{\xi}||^{2}_{L^{\infty}_{t}L^{2}_{x}}=||f_{\xi}||^{2}_{Y_{\xi}}$$

Summing up with respect to $\xi \in \Xi$, $|\xi| \approx 2^{i}$ gives us the general result $f_{i,\leq d} \in Y^{s}$. If for each $i$ we take $d=1$ we can sum up with respect to all $\xi \in \Xi$ and obtain $f_{\cdot, \leq 1} \in Y^{s}$. One can easily observe that the decay in time structure is conserved in all the above computations, therefore we also get $f_{i,\leq d} \in Y^{s,N}$ and $f_{\cdot, \leq 1} \in Y^{s,N}$.

\end{proof}

\subsection{Proof of Proposition \ref{k4}} 

We denote by $v=e^{it\D}g$. It is well known that $v \in C_{t}H^{s}_{x}$ and that $\chi_{[0,1]}v \in X^{s,b}$ for any $b$, see for example \cite{k5}. As a consequence $\chi_{[0,1]}v \in X^{s,\q,1}$.

The delicate part of the proof is to show that $(\chi_{[0,1]}v)_{\cdot \leq 1} \in Y^{s}$. For this purpose we decompose the initial data in the following way (see (\ref{d100})):

$$g= \sum_{\xi \in \Xi} \sum_{m} \check{\tilde{\phi}}_{\xi} * (\chi_{Q^{m}} g_{\xi}) =  \sum_{\xi \in \Xi} \sum_{m} g_{\xi}^{m}$$

Then we can write:

\beq \label{s20}
v =\sum_{\xi \in \Xi} v_{\xi} =\sum_{\xi \in \Xi} \sum_{m} v_{\xi}^{m}
\eeq

\noindent
where $v_{\xi}^{m}=e^{it\Delta}g_{\xi}^{m}$. We fix $\xi=\xi^{0} \in \Xi$ and recall the well-known energy conservation for the homogeneous Schr\"odinger equation.

\begin{l1}
For each $\xi^{0}$ and $m$ we have:

\beq \label{s5}
||v_{\xi^{0}}^{m}||_{L^{\infty}_{t}L^{2}_{x}} =  ||g_{\xi^{0}}^{m}||_{L^{2}}
\eeq

\end{l1}

$\chi_{[0,1]} v_{\xi^{0}}^{m}$ has a phase-space concentration: it is highly concentrated in $T_{\xi_{0}}^{m}$ in space and in a neighborhood of size $1$ around $(\xi_{0},\xi_{0}^{2})$ in frequency. The next Lemma bellow
prepares step by step the estimates necessary to prove this claim. Before that we need to define $P_{j,m}(x,t,D)$, $j=1,2$ and $m \in \Z$, to be the differential operators:

$$P_{j,m} = x_{j}-m_{j} - i 2t D_{x_{j}} \ \ \ \mbox{with symbols} \ \ \ p_{j,m} = x_{j}-m_{j} + 2t\xi_{j}$$

\begin{l1} $i\partial_{t} - \D$ commutes with both $P_{j,m}(x,t,D)$, j=1,2. 

For each $m \in \Z$ and $j=1,2$ we have the estimates:

 \beq \label{s6}
||P^{n}_{j,m}(x,t,D)v_{\xi^{0}}^{m}||_{L^{\infty}_{t}L^{2}_{x}} =  ||P^{n}_{j,m}(x,0,D)g_{\xi^{0}}^{m}||_{L^{2}}
\eeq

\beq \label{s7}
||P_{j,m}^{n}(x,0,D)g_{\xi^{0}}^{m}||_{L^{2}} \leq C_{n} ||\chi_{Q^{m}} g_{\xi^{0}}||_{L^{2}}
\eeq

\beq \label{s10}
||(x_{j}-m_{j} +2t \xi^{0}_{j})^{n} \chi_{[0,1]} v_{\xi^{0}}^{m}||_{L^{\infty}_{t}L^{2}_{x}} \leq  C_{n} ||\chi_{Q^{m}} g_{\xi^{0}}||_{L^{2}}
\eeq

\beq \label{s13}
||\chi_{[0,1]} v_{\xi^{0}}^{m}||_{L^{\infty}_{t}L^{2}_{x}(T^{m',0}_{\xi^{0}})} \leq C_{n} \langle m-m' \rangle^{-n}  C_{n} ||\chi_{Q^{m}} g_{\xi^{0}}||_{L^{2}}
\eeq

\end{l1}

\begin{proof}

The fact that $i\partial_{t} - \D$ commutes with both $P_{j,m}(x,t,D)$, j=1,2 can be verified by direct computation. As a consequence $P_{j,m}^{n}v_{\xi_{0}}^{m}$ is also a solution of the homogeneous equation and then the energy conservation gives us (\ref{s6}).

In order to prove (\ref{s7}) we start with:

$$P_{j,m}^{n}(x,0,D)g_{\xi^{0}}^{m} = (x_{j}-m_{j})^{n} \tilde{\phi}_{\xi^{0}}(D) (\chi_{Q^{m}} g_{\xi^{0}}) $$

Standard calculus gives us:

$$(x_{j}-m_{j})^{n} \tilde{\phi}_{\xi^{0}}(D)=\sum_{k=1}^{n} \frac{\partial^{k} \tilde{\phi}_{\xi^{0}}}{\partial \xi_{j}^{k}}(D) (x_{j}-m_{j})^{k} $$

A simple computation shows that $||\frac{\partial^{k} \tilde{\phi}_{\xi^{0}}}{\partial \xi_{j}^{k}}||_{L^{\infty}_{\xi}} \leq C_{k}$. In addition we have that $|x_{1}-m_{1}| \leq 2$ in the support of $\chi_{Q^{m}}$. Therefore we can conclude that for all $k$'s in the above sum we have: 

$$||\frac{\partial^{k} \tilde{\phi}_{\xi^{0}}}{\partial \xi_{j}^{k}}(D) (x_{j}-m_{j})^{k}\chi_{Q^{m}} g_{\xi^{0}})||_{L^{2}} \leq  C_{k} ||\chi_{Q^{m}} g_{\xi^{0}}||_{L^{2}}$$

This is enough to justify the claim of (\ref{s7}). Next we expand:

$$(x_{j}-m_{j}+2t\xi_{j}^{0})^{n}=(x_{j}-m_{j}- i2tD_{x_{j}} +(i2tD_{x_{j}}+2t\xi_{j}^{0}))^{n}=$$

$$\sum_{\a+\b +\g \leq n} C_{\a\b\g}P_{j,m}^{\a}(x,t,D) (i2tD_{x_{j}}+2t\xi_{j}^{0})^{\beta}t^{\g}$$

We took into account that the commutator $[P_{j,m}^{\a}(x,t,D),i2tD_{x_{j}}+2t\xi_{j}^{0})]$ has symbol $C \cdot t$, for some constant C which can be explicitly computed. We make two observations: in the support of $\hat{v}_{\xi^{0}}^{m}$ we have $|\xi^{0}_{1}-\xi_{1}| \leq 1$ and $|2t \chi_{[0,1]}| \leq 2$. Then using (\ref{s6}) and (\ref{s7}) gives us (\ref{s10}). 

The proof of (\ref{s13}) is a direct consequence of (\ref{s10}).

\end{proof}

We continue now with the proof of Proposition 1.  Recalling (\ref{s20}) we estimate:

\begin{align*}
 &||\chi_{[0,1]} v_{\xi^{0}}||^{2}_{Y} = \sum_{m'} ||\chi_{[0,1]} v_{\xi^{0}}||^{2}_{L^{\infty}_{t}L^{2}_{x}(T^{m',0}_{\xi^{0}})} \\ 
 &\leq \sum_{m'} \left( \sum_{m} ||\chi_{[0,1]} v^{m}_{\xi^{0}}||_{L^{\infty}_{t}L^{2}_{x}(T^{m',0}_{\xi^{0}})}  \right)^{2} \\ 
 &\leq \sum_{m'} \left( \sum_{m} C_{n} \langle m-m' \rangle^{-n}   ||\chi_{Q^{m}} g_{\xi^{0}}||_{L^{2}}  \right)^{2} \\
 &\leq C^{2}_{n}\sum_{m'} \sum_{m}  \langle m-m' \rangle^{-n}  ||\chi_{Q^{m}} g_{\xi^{0}}||^{2}_{L^{2}}  \\ 
 &\leq C^{2}_{n} \sum_{m} ||\chi_{Q^{m}} g_{\xi^{0}}||^{2}_{L^{2}} \leq C^{2}_{n}  || g_{\xi^{0}}||^{2}_{L^{2}}
\end{align*}

In the above computations we used twice the inequality:

$$\sum_{m} \langle m-m' \rangle^{-n}\leq C$$

\noindent
which is true as long as $n \geq 3$. At the level of $Y^{s}$ the above estimate becomes $||\chi_{[0,1]} v_{\xi^{0}}||_{Y^{s}} \leq ||g_{\xi^{0}}||_{H^{s}}$. Summing up with respect to $\xi^{0} \in \Xi$ gives us $||\chi_{[0,1]} v||_{Y^{s}} \leq ||g||_{H^{s}}$. Notice also that we get for free the estimate $||\chi_{[0,1]} v||_{Y^{s,N}} \leq ||g||_{H^{s}}$, see part d, Lemma \ref{ees}.

From part e) of Lemma \ref{ees} we can conclude $||(\chi_{[0,1]} v)_{\cdot, \leq 1}||_{Y^{s,N}} \leq ||g||_{H^{s}}$

\vspace{.1in}

\subsection{Proof of Proposition \ref{k3}}

\noindent
\vspace{.1in}

The claim is that if $f \in W^{s}$ then $\chi_{[0,1]}w \in Z^{s,N} \cap C_{t}H^{s}_{x}$. We rewrite (\ref{e2}):

$$iw_{t}-\D w=f_{\cdot,\geq 1}+f_{\cdot,\leq 1}=f_{\cdot,\geq 1}+f^{1}+f^{3}$$

\noindent
where $f^{1} \in X^{s,-\q,1}$ and $f^{3} \in \mathcal{Y}^{s}$. We decompose more:

$$iw_{t}-\D w=f_{\cdot,\geq 1}+\sum_{j}f^{1}_{j, \geq 2^{-j}}+ \sum_{j}f^{1}_{j, \leq 2^{-j}} +f^{3}$$

The solution of (\ref{e2}) can be written:

$$w=w^{1}+w^{2}+w^{3}-e^{it\D}v$$

\noindent
where $w^{1},w^{2},w^{3}$ are given by:

\beq \label{ee1}
(\tau-\xi^{2})\hat{w}^{1}=\hat{f}_{\cdot,\geq 1}+\sum_{j}\hat{f}^{1}_{j, \geq 2^{-j}}
\eeq

\begin{equation} \label{ee2} 
\begin{cases}
\begin{aligned} 
&(i\partial_{t} - \Delta) w^{2}=\sum_{j}f^{1}_{j, \leq 2^{-j}}\\
&w^{2}(x,0)=0
\end{aligned}
\end{cases}
\end{equation}

\begin{equation} \label{ee3} 
\begin{cases}
\begin{aligned} 
&(i\partial_{t} - \Delta) w^{3}=f^{3}\\
&w^{3}(x,0)=0
\end{aligned}
\end{cases}
\end{equation}

The correction factor $e^{it\D}v$ is present since the solution of (\ref{ee1}) does not necessarily come with zero initial data. But we prove  that it corresponds to an initial data in $v \in H^{s}$ and this justifies our correction.     

We will focus on showing the following properties:

\beq \label{ee6}
w^{1}, \chi_{[0,1]}w^{2} \in X^{s,\q,1} \ \ \mbox{and} \ \ \chi_{[0,1]}w^{3} \in Z^{s,5} \cap C([-1,2]:H^{s})
\eeq

If we assume for the moment the properties (\ref{ee6}), we can prove the result of the Proposition \ref{k3}.

\begin{proof}[Proof of Proposition \ref{k3}]

From part a) and b) of Lemma \ref{ees} we can conclude that $\chi_{[0,1]}w^{1},\chi_{[0,1]}w^{2} \in C([-1,2]:H^{s})$, therefore $\chi_{[0,1]}(w^{1}+w^{2}+w^{3}) \in C([-1,2]:H^{s})$. Thus $(w^{1}+w^{2}+w^{3})(0)=v \in H^{s}$ which implies $w=w^{1}+w^{2}+w^{3}-e^{it\D}v$ is the solution of (\ref{e2}). Since $v \in H^{s}$, we can use the result in Proposition \ref{k4} to obtain $\chi_{[0,1]}e^{it\D}v \in Z^{s,5} \cap C([-1,2]:H^{s})$. 

Part b) of Lemma \ref{ees} gives us that $\chi_{[0,1]}w^{1} \in X^{s,\q,1}$ and part e) implies  $(\chi_{[0,1]}w^{1}+\chi_{[0,1]}w^{2})_{\cdot, \leq 1} \in X^{s,\q,1} \subset X^{s,\q,\infty}$ and  $(\chi_{[0,1]}w^{1}+\chi_{[0,1]}w^{2})_{\cdot, \geq 1} \in X^{s,\q,1}$.

From  Lemma \ref{ees}, part c),  $\chi_{[0,1]}w^{1}+\chi_{[0,1]}w^{2} \in  X^{s,\q,1} \subset Y^{s}$ and using part d) we obtain that $\chi_{[0,1]}w^{1}+\chi_{[0,1]}w^{2} \in  X^{s,\q,1} \subset Y^{s,5}$. In the end we invoke part e) of the same Lemma to conclude with  $(\chi_{[0,1]}w^{1}+\chi_{[0,1]}w^{2})_{\cdot, \leq 1} \in Z^{s,5}$.

\end{proof}

We continue with the proof of the claims in (\ref{ee6}). 

\vspace{.1in}

\begin{bfseries} 

{\mathversion{bold}$w^{1} \in X^{s,\q,1} $}.

\end{bfseries}

The equation (\ref{ee1}) can be written in the form:

$$\hat{w}_{1}=\frac{1}{\tau-\xi^{2}}(\hat{f}_{\cdot,\geq 1}+\sum_{j}\hat{f}^{1}_{j, \geq 2^{-j}})$$

\noindent
and the right hand side is localized in a region where $|\tau-\xi^{2}| \geq 1$, therefore:

$$||w^{1}||_{X^{s,\q,1}} \leq ||f_{\cdot,\geq 1}+\sum_{j} f^{1}_{j, \geq 2^{-j}} ||_{X^{s,-\q,1}}$$

\begin{bfseries} 

{\mathversion{bold}$\chi_{[0,1]}w^{2} \in X^{s,\q,1} $}

\end{bfseries}

We write $w^{2}=\sum_{j} w^{2}_{j}$, where $w^{2}_{j}$ solves the equation:

\begin{equation} \label{ee11} 
\begin{cases}
\begin{aligned} 
&(i\partial_{t} - \Delta) w^{2}_{j}=f^{1}_{j, \leq 2^{-j}}\\
&w^{2}(x,0)=0
\end{aligned}
\end{cases}
\end{equation}

\noindent
for each $j$. The standard energy estimate gives us:

$$||\chi_{[0,1]} w^{2}_{j}||_{L^{2}} \leq ||f^{1}_{j,\leq 2^{-j}}||_{L_{t}^{1}L^{2}_{x} ([-1,2] \times \R^{2})} \leq ||f^{1}_{j,\leq 2^{-j}}||_{L^{2}}$$

 $\hat{f}^{1}_{j,\leq 2^{-j}}$ is localized in $A_{j, \leq 2^{-j}}$, hence $\hat{w}^{2}_{j}$ is localized in $A_{j, \leq 2^{-j}}$. The cut in time spreads the support, but $\widehat{\chi_{[0,1]}w^{2}_{j}}$ is highly localized in $A_{j, \leq 2^{-j}}$:

\beq \label{hh1} 
|| \chi_{[\xi^{2}+k, \xi^{2}+k+1]} \widehat{\chi_{[0,1]}w^{2}_{j}} (\xi,\cdot)||_{L^{2}_{\tau}} \leq C_{N} \langle k \rangle^{-N} \sum_{\a=0}^{N} ||\frac{\partial^{N-\a} \chi_{[0,1]}}{\partial t^{N-\a}} w^{2,\a}_{j} ||_{L^{2}_{t}}
\eeq

\noindent
where $\chi_{[\xi^{2}+k, \xi^{2}+k+1]}$ is a smooth approximation of the characteristic function of the interval $[\xi^{2}+k, \xi^{2}+k+1]$ (in the $\tau$ variable) and $\hat{w}^{2,\a}_{j}=(\tau-\xi^{2})^{\a}\hat{w}^{2}_{j} (\xi,\tau)$. This estimate can be proved by using the commutator identity:

$$(i \partial_{t} - \xi^{2})^{N} \chi_{[0,1]}w^{2}_{j}=\sum_{\a} {N \choose \a} \frac{\partial^{N-\a} \chi_{[0,1]}}{\partial t^{N-\a}} (i \partial_{t} - \xi^{2})^{\a} w^{2}_{j} $$  

We denote by $\hat{f}^{1,\a}_{j}=(\tau-\xi^{2})^{\a}\hat{f}^{1}_{j} (\xi,\tau)$. Since $|\tau-\xi^{2}| \leq 1$ in the support of $\hat{f}^{1}_{j,\leq 2^{-j}}$ it follows that

$$||\frac{\partial^{N-\a} \chi_{[0,1]}}{\partial t^{N-\a}} w^{2,\a}_{j}||_{L^{2}} \leq ||f^{1,\a}_{j,\leq 2^{-j}}||_{L_{t}^{1}L^{2}_{x} ([-1,2] \times \R^{2})} \leq ||f^{1}_{j,\leq 2^{-j}}||_{L^{2}}$$

 We have all ingredients to claim that $\chi_{[0,1]}w^{2} \in X^{s,\q,1} $.

\vspace{.1in}

\begin{bfseries} 

{\mathversion{bold}$\chi_{[0,1]}w^{3} \in Z^{s,5} \cap C([-1,2]:H^{s})$}.

\end{bfseries}

 We decompose the inhomogeneous term in the following way (see (\ref{d100})):

$$f^{3}= \sum_{\xi \in \Xi} \sum_{(m,l) \in \Z^{3}} \check{\tilde{\phi}}_{\xi} * (\chi_{T_{\xi}^{m,l}} f_{\xi}^{3}) =  \sum_{\xi \in \Xi} \sum_{(m,l) \in \Z^{3}} f_{\xi}^{m,l}$$

The solution of $(\ref{ee2})$ can be written as: 

\beq \label{i3}
w^{3}= \sum_{\xi \in \Xi} \sum_{(m,l) \in \Z^{3}} w_{\xi}^{m,l}
\eeq

\noindent
where for each $\xi^{0} \in \Xi$, $w_{\xi^{0}}^{m,l}$ satisfies the equation

\begin{equation} \label{i555} 
\begin{cases}
\begin{aligned} 
&(i\partial_{t} - \Delta) w_{\xi^{0}}^{m,l}=f_{\xi^{0}}^{m,l}\\
&w_{\xi^{0}}^{m,l}(x,0)=0
\end{aligned}
\end{cases}
\end{equation}

We want estimates for $\chi_{[0,1]}w^{3}$, hence we want estimates for $\chi_{[0,1]}w_{\xi^{0}}^{m,l}$. Since $f_{\xi^{0}}^{m,l}$ is essentially supported in the time interval $[l,l+1]$ it follows that $\chi_{[0,1]}w_{\xi^{0}}^{m,l}=0$ for $l \ne 0,\pm 1$. Therefore, all the estimates bellow are for $l = 0,\pm 1$

For the solution of (\ref{i555}) we have the standard energy estimate:

\beq \label{i1}
||w^{m,l}_{\xi^{0}}||_{L^{\infty}_{t}L^{2}_{x}} \leq  ||f_{\xi^{0}}^{m,l}||_{L^{1}_{t}L^{2}_{x}}
\eeq

From this point on the argument follows the same steps as the one for the homogeneous equation. The only difference is that we  work with $P_{j,m}(x,t,D)$ instead of $P_{j,m}(x,0,D)$. 

 We list a sequence of results whose proofs can be derived in a similar way as their correspondent in the homogeneous case.

\begin{l1} \label{i33}

For each $m \in \Z$, $l = 0,\pm 1$ and $j=1,2$ we have the estimates:

\beq \label{i110}
||P_{j,m}^{n}(x,t,D)w_{\xi^{0}}^{m,l}||_{L^{\infty}_{t}L^{2}_{x}} \leq C_{N} ||P_{j,m}^{n}(x,t,D)f_{\xi^{0}}^{m,l}||_{L^{1}_{t}L^{2}_{x}}
\eeq
 
\beq \label{i4}
||P_{j,m}^{n}(x,t,D)f_{\xi^{0}}^{m,l}||_{L^{1}_{t}L^{2}_{x}} \leq C_{n} || \chi_{T^{m,l}_{\xi^{0}}} f_{\xi^{0}}^{3}||_{L^{1}_{t}L^{2}_{x}}
\eeq

\beq \label{i9}
||(x_{j}-m_{j} +2t \xi^{0}_{j})^{n} \chi_{[0,1]} w_{\xi^{0}}^{m,l}||_{L^{\infty}_{t}L^{2}_{x}} \leq   C_{n} || \chi_{T^{m,l}_{\xi^{0}}} f_{\xi^{0}}^{3}||_{L^{1}_{t}L^{2}_{x}}
\eeq

\beq \label{i12}
||\chi_{[0,1]} w_{\xi^{0}}^{m,l}||_{L^{\infty}_{t}L^{2}_{x}(T^{m',0}_{\xi^{0}})} \leq C_{n} \langle m-m' \rangle^{-n} || \chi_{T^{m,l}_{\xi^{0}}} f_{\xi^{0}}^{3}||_{L^{1}_{t}L^{2}_{x}}
\eeq

\beq \label{i34}
||\chi_{[0,1]} w_{\xi^{0}}^{3}||_{Y^{n}} \leq C_{n} ||f_{\xi^{0}}^{3}||_{\mathcal{Y}^{n}}
\eeq

\end{l1}

Using (\ref{i34}) and performing the summation with respect to $\xi \in \Xi$ implies that $||\chi_{[0,1]} w^{3}||_{Y^{s,n}} \leq C_{n} ||f^{3}||_{\mathcal{Y}^{s,n}}$. Part e) of Lemma \ref{ee6} gives us that $||(\chi_{[0,1]} w^{3})_{\cdot, \leq 1}||_{Y^{s,n}} \leq ||\chi_{[0,1]} w^{3}||_{Y^{s,n}} \leq C_{n} ||f^{3}||_{\mathcal{Y}^{s,n}}$. 

The solution of (\ref{i555}) has a better property, namely:

\beq \label{i11}
||w^{m,l}_{\xi^{0}}||_{C_{t}L^{2}_{x}} \leq  ||f_{\xi^{0}}^{m,l}||_{L^{1}_{t}L^{2}_{x}}
\eeq

Going through the same machinery as in Lemma \ref{i33} gives us $\chi_{[0,1]} w^{3} \in C_{t}H^{s}_{x}$.

We can also write $w^{3}=\sum_{j} w^{3}_{j}$ where 

\begin{equation} \label{i666} 
\begin{cases}
\begin{aligned} 
&(i\partial_{t} - \Delta) w_{j}^{3}=f_{j}^{3}\\
&w_{j}^{3}(x,0)=0
\end{aligned}
\end{cases}
\end{equation}

For fixed $d \in I_{j}$, $d > 2^{-j}$ we have $\hat{w}^{3}_{j,d}=(\tau-\xi^{2})^{-1}\hat{f}^{3}_{j,d}$, hence:

$$  ||w^{3}_{j,d}||_{X^{s,\q}} \leq ||f^{3}_{j,d}||_{X^{s,-\q}} \leq ||f^{3}_{j}||_{X^{s,-\q, \infty}} \leq ||f^{3}_{j}||_{\mathcal{Y}^{s}}$$

In the last estimate we have used part c) of Lemma \ref{ees} and the fact that both $f^{3}$ and $w^{3}$ are localized at distance less than $1$ from $P$. This was important since $X^{s,\q,\infty}$ has an improved weight, i.e. $\langle (\xi,\tau) \rangle^{s}$, over $\mathcal{Y}^{s}$, i.e. $\langle \xi \rangle^{s}$, and this makes a difference away from $P$. Hence we can claim:

$$  ||w^{3}_{j,\geq 2^{-j}}||_{X^{s,\q,\infty}} \leq ||f^{3}_{j}||_{\mathcal{Y}^{s}}$$

For $d=2^{-j}$ we have:

$$||\chi_{[0,1]}w^{3}_{j,\leq 2^{-j}}||_{L^{2}} \leq ||\chi_{[0,1]}w^{3}_{j}||_{L^{2}}+ ||\chi_{[0,1]}w^{3}_{j, \geq 2^{-j}}||_{L^{2}} \leq $$

$$||\chi_{[0,1]}w^{3}_{j}||_{Y^{0}}+||w^{3}_{j, \geq 2^{-j}}||_{L^{2}} \leq ||f^{3}_{j}||_{\mathcal{Y}^{0}}+||w^{3}_{j, \geq 2^{-j}}||_{X^{0,\q,\infty}} \leq ||f^{3}_{j}||_{\mathcal{Y}^{0}}$$

Since $\hat{w}^{3}_{j,\leq 2^{-j}}$ is concentrated in a region where $|\tau-\xi^{2}| \leq 1$ and $\hat{\chi}_{[0,1]}$ is highly concentrated in $[-1,1]$ a similar argument to the one in {\mathversion{bold}$\chi_{[0,1]}w^{2} \in X^{s,\q,1} $} gives us that actually:

$$||\chi_{[0,1]}w^{3}_{j,\leq 2^{-j}}||_{X^{s,\q,1}} \leq ||f^{3}_{j}||_{\mathcal{Y}^{s}}$$

We are left with proving that $(\chi_{[0,1]} w^{3}_{j, \geq 2^{-j}})_{\cdot, \geq 1} \in X^{s,\q,1}$. We know that $w^{3}_{j, \geq 2^{-j}} \in X^{s,\q,\infty}$ and $\hat{w}^{3}_{j, \geq 2^{-j}}$ is supported at distance less than $1$ from $P$. The key idea is that $\hat{\chi}_{[0,1]}$ is highly localized in $[-1,1]$. To formalize a bit, if we denote by $h_{k}(\tau)=h(\tau) \cdot \chi_{[k-\q,k+\q]}(\tau)$ (here $\chi_{[k-\q,k+\q]}$ is a smooth approximation of the characteristic function of the interval $[k-\q,k+\q]$) we have:

$$||(\hat{\chi}_{[0,1]} * h_{k})_{k'}||_{L^{2}} \leq C_{N} \langle k-k'\rangle^{-N} ||h_{k}||_{L^{2}}$$

Applying this to $h_{\xi}=\hat{w}^{3}_{j, \geq 2^{-j}}(\xi, \cdot)$, for each $\xi$ with $|\xi| \approx 2^{j}$, it is a matter of algebraic computations to obtain:

$$||(\chi_{[0,1]} w^{3}_{j, \geq 2^{-j}})_{\cdot, \geq 1}||_{X^{s,\q,1}} \leq ||w^{3}_{j, \geq 2^{-j}}||_{X^{s,\q,\infty}} \leq ||f^{3}_{j}||_{\mathcal{Y}^{s}}$$

\vspace{.1in}

\subsection{Conservation of decay} \label{sub}
\noindent
\vspace{.1in}

  In this section we complete the proof of (\ref{e3}) and (\ref{e4}) in the sense that we add the decay structure. The decay we use is scaled properly for the Schr\"odinger equation and this is why, in principle, it should be easily conserved.

Before we start we need to introduce some new localization operators. For each $k \in \N$ we define the lattice:

\beq \label{Xi}
\Xi^{k}=\{\xi=(n2^{-k},\t): \t=\frac{\pi}{2} \frac{l}{n}; n \in \N, l \in \Z \}
\eeq

\noindent
and for each $\xi \in \Xi^{k}$ we build  the corresponding $\phi^{k}_{\xi}$ to be a smooth approximation of the characteristic function of the cube centered at $\xi$ and with sizes $2^{-k}$. We also assume that the system $(\phi_{\xi}^{k})_{\xi \in \Xi^{k}}$ forms a partition of unity in $\R^{2}$. For every $\xi \in \Xi^{k}$ we define:

$$\mathcal{F}(S_{\xi}^{k}f)= \phi^{k}_{\xi} \cdot \hat{f}$$

For each $k$ we define $\Xi^{k,*} \subset \Xi^{k}$ to be the subset of those $\xi$ with $|\xi| \approx 2^{k}$ and such that we have:

$$\sum_{k} \sum_{\xi \in \Xi^{k,*}} \phi^{k}_{\xi} =1$$

For each $l \in \Z$ we can easily construct a function $\chi_{[l-\q,l+\q]}$ to be a smooth approximation of the characteristic function of the interval $[l-\q,l+\q]$ and such that the system $(\chi_{[l-\q,l+\q]})_{l \in \Z}$ form a partition of unity in $\R$. For any $\xi \in \Xi^{k}$ with $|\xi| \leq 2^{k+1}$ we consider those $l \in \Z$ with the property $|(\xi,l)| \approx 2^{k}$ and define the operators: 

$$\mathcal{F}(T_{\xi,l}f) (\xi,\tau)=\hat{f}_{\xi,l}= \phi_{[l-\q,l+\q]}(\tau) \phi^{k}_{\xi}(\xi) \hat{f}(\xi,\tau) $$  

This is one example when we localize in a linear way with respect to the size of $\tau-\xi^{2}$ rather than a dyadic way, as commented in the section 2. 

The main result of this section which will help us to prove the conservation of decay is the following:

\begin{p2} \label{dec}
 a) We have the approximation:

\beq \label{d1111}
||\mathcal{D}f||^{2}_{X^{s,\q} \cap Y^{s,N}} \approx   \sum_{i} \sum_{\xi \in \Xi^{i,*}}  \sum_{l \in \Z: |(\xi,l)| \approx 2^{i}} \langle (\xi,l)\rangle^{2s} \langle l-\xi^{2} \rangle ||\langle \frac{x}{2^{i}} \rangle^{\q+\e} f_{\xi,l} ||^{2}_{L^{2}} + 
\eeq

$$\sum_{\xi \in \Xi} 2^{2|\xi|s} ||\langle \frac{x}{\langle \xi \rangle} \rangle^{\q+\e} f_{\xi, \leq 1}||_{Y^{0,N}}^{2}$$

b) The result in a) holds true if we replace $X^{s,\q} \cap Y^{s,N}$ with $X^{s,\q,1}_{\cdot, \geq 1} \cap X^{s,\q, \infty}_{\cdot, \leq 1} \cap Y^{s,N}_{\cdot, \leq 1}$.

\end{p2}

We did not write down the result claimed in part b) because of its complicated formulation. The constant in the $\approx$ relation depends on $\mu$, but we choose to ignore it. $\mu$ is present only for technical reasons and it does not affect the computations with more than a $C_{\mu}$. This is why we keep track of it only when we do symbolic calculus and discard it later on.

The result of Proposition \ref{dec} are of commutator type. The norm of $\mathcal{D}f$ in $X^{s,\q} \cap Y^{s,N}$ is:

$$\sum_{i} \sum_{\xi \in \Xi^{i, *}}  \sum_{l \in \Z: |(\xi,l)| \approx 2^{i}} \langle (\xi,l)\rangle^{2s} \langle l-\xi^{2} \rangle ||(\mathcal{D} f)_{\xi,l} ||^{2}_{L^{2}} + \sum_{\xi \in \Xi} 2^{2|\xi|s} ||(\mathcal{D}f)_{\xi, \leq 1}||_{Y^{0,N}}^{2}$$

The above result claims that we can commute $D$ with the localization $T_{\xi,l}$ and keep the same weights even if $\mathcal{F}(\mathcal{D} f_{\xi,l})$ does not have the same support as $\hat{f}_{\xi,l}$; in fact $\mathcal{F}(\mathcal{D} f_{\xi,l})$ is not even compactly supported. The reason this works is that $\mathcal{F}(\mathcal{D} f_{\xi,l})$ is mainly localized in the support of $T_{\xi,l}$ and decreases rapidly outside this support. 

Another fact which is implicit in the statement of Proposition \ref{dec} is that if $f$ is localized at frequency $2^{i}$ then $\mathcal{D}f$ can be seen as $\langle \frac{x}{2^{i}} \rangle^{\q+\e} f$; in other words we can freeze the frequency part of the symbol of $\mathcal{D}$ and see the symbol as a multiplier.

We now turn to the proof of Proposition \ref{dec}. We deal with a commutator problem. The symbol of $\mathcal{D}$ is $d(x, \xi, \tau)=(1+\frac{x^{2}}{\mu +|\tau|+\xi^{2}})^{\frac{1}{4}+\e}$ commutes with frequency localizations in $\tau$. Hence we can ignore the time component of the problem for a while and consider $\tau$ as a parameter in the expression of the symbols bellow. 

We have to recall some basics of the theory of the hypoelliptic operators as developed, for instance, in \cite{ho}. An symbol $p(x,\xi)$ is said to be hypoelliptic if it satisfies the following condition:

\beq \label{h1}
|\partial_{\xi}^{\a} \partial_{x}^{\b} p(x,\xi)| \leq C_{\a \b} |p(x,\xi)| \langle \xi \rangle^{-|\a|}, \ \ |\xi| \geq C
\eeq

\noindent
for some $C$. The operator pseudo-differential operator $P(x,D)$ with symbol $p(x,\xi)$ is invertible in the sense that there is $P^{-1}(x,D)$ with the properties:

$$PP^{-1}=I+R_{1}  \ \ \ \ \ \ P^{-1}P=I+R_{2}$$

\noindent
where $R_{1}, R_{2}$ are of order $-\infty$. In addition the symbol $q(x,\xi)$ of $P^{-1}(x,D)$ satisfies:

\beq \label{h2}
|\partial_{\xi}^{\a} \partial_{x}^{\b} q(x,\xi)| \leq C_{\a \b} |p(x,\xi)|^{-1} \langle \xi \rangle^{-|\a|}, \ \ |\xi| \geq C
\eeq

$d(x,\xi,\tau)$ is hypoelliptic. Moreover we added the constant $\mu >> 1$ with a sole purpose: to be able to take $C=0$ in (\ref{h1}). If we denote by $e(x,\xi,\tau)$ the symbol of $\mathcal{D}^{-1}$ (defined up to an operator of order $- \infty$), $e(x,\xi,\tau)$ has bounds of type (\ref{h2}) where $p(x,\xi)$ is replaced by $d(x,\xi,\tau)$. 

We can prove our first result:

\begin{l1} We have the estimate:
\beq \label{dec22}
||\mathcal{D}f_{i}||_{L^{2}} \approx ||\langle \frac{x^{2}}{\mu+2^{2i}}\rangle^{\frac{1}{4}+\frac{\e}{2}} f_{i}||_{L^{2}} \approx C_{\mu} ||\langle \frac{x}{2^{i}}\rangle^{\q+\e} f_{i}||_{L^{2}}
\eeq
\end{l1}

\begin{proof}
The second estimate is obvious, hence we have to deal only with the first one. We first prove the inequality:

$$||\mathcal{D}f_{i}||_{L^{2}} \leq ||\langle \frac{x^{2}}{\mu+2^{2i}}\rangle^{\frac{1}{4}+\frac{\e}{2}} f_{i}||_{L^{2}}$$

If we denote by $h=\langle \frac{x^{2}}{\mu+2^{2i}}\rangle^{\frac{1}{4}+\frac{\e}{2}} f_{i}$ we observe that $f_{i}=\phi_{i}(D) \langle \frac{x^{2}}{\mu+2^{2i}} \rangle^{-\frac{1}{4}-\frac{\e}{2}} h$, so we have to show that:

$$||\mathcal{D}\phi_{i}(D) \langle \frac{x^{2}}{\mu+2^{2i}} \rangle^{-\frac{1}{4}-\frac{\e}{2}} h||_{L^{2}} \leq ||h||_{L^{2}}$$

It is enough to show that the composition $\mathcal{D}\phi_{i}(D) \langle \frac{x^{2}}{\mu+2^{2i}} \rangle^{-\frac{1}{4}-\frac{\e}{2}}$ is in the standard class $\Psi^{0}_{1,0}$. The symbol of $\mathcal{D}\phi_{i}(D)$ is $d(x,\xi,\tau) \cdot \phi_{i}(\xi,\tau)$ and satisfies:

$$|\partial_{\xi}^{\a} \partial_{x}^{\b} (d(x,\xi,\tau) \cdot \phi_{i}(\xi,\tau))| \leq C_{\a \b} \langle \frac{x^{2}}{\mu+2^{2i}} \rangle^{\frac{1}{4}+\frac{\e}{2}} \langle \xi \rangle^{-|\a|}$$

$\langle \frac{x^{2}}{\mu+2^{2i}} \rangle^{-\frac{1}{4}-\frac{\e}{2}}$ is obviously a hipoelliptic operator and this is enough to invoke the standard theory for composition of pseudo-differential operators to obtain that $\mathcal{D}\phi_{i}(D) \langle \frac{x^{2}}{\mu+2^{2i}} \rangle^{-\frac{1}{4}-\frac{\e}{2}}$ is in the standard class $\Psi^{0}_{1,0}$.

We are left with the second inequality:

$$||\mathcal{D}f_{i}||_{L^{2}} \geq ||\langle \frac{x^{2}}{\mu+2^{2i}}\rangle^{\frac{1}{4}+\frac{\e}{2}} f_{i}||_{L^{2}}$$

If we denote by $h=\mathcal{D} f_{i}$ we observe that $f_{i}=\phi_{i}(D) \mathcal{D}^{-1} h$, so we have to show that:

\beq \label{dd11}
||\langle \frac{x^{2}}{\mu+2^{2i}} \rangle^{\frac{1}{4}+\frac{\e}{2}} \phi_{i}(D) \mathcal{D}^{-1} h||_{L^{2}} \leq ||h||_{L^{2}}
\eeq

The symbol $q(x,\xi,\tau)$ of the operator $\phi_{i}(D) \mathcal{D}^{-1}$ has an asymptotical expansion of type:

$$q(x,\xi,\tau)=\sum_{\a} \frac{1}{\a!} (i\partial_{\xi})^{\a} \phi_{i}(\xi,\tau) \partial_{x} e(x,\xi,\tau)$$

 The bounds on its derivatives which are of type:

$$|\partial_{\xi}^{\a} \partial_{x}^{\b} q(x,\xi,\tau)| \leq C_{\a \b} \langle \frac{x^{2}}{\mu+2^{2i}} \rangle^{-\frac{1}{4}-\frac{\e}{2}} \langle \xi \rangle^{-|\a|}$$

Here we have used the bounds we have for $e(x,\xi,\tau)$ from the standard theory. Again this is enough to conclude that $\langle \frac{x^{2}}{\mu+2^{2i}} \rangle^{\frac{1}{4}+\frac{\e}{2}} \phi_{i}(D) \mathcal{D}^{-1} \in \Psi^{0}_{1,0}$ and then conclude that (\ref{dd11}) is true. 

\end{proof}

\begin{l1} \label{decay}

If $f$ is localized at frequency $2^{i}$ and $\xi_{\a} \in \Xi^{i}$, $|\xi_{\a}| \leq 2^{i}$, we have the following estimates:

\beq \label{d1000}
||S_{\xi_{\b}} \mathcal{D}S_{\xi_{\a}}f||_{L^{2}} \leq C_{N} (2^{i}|\xi_{\b}-\xi_{\a}|)^{-N} ||\mathcal{D}S_{\xi_{\a}}f||_{L^{2}}
\eeq

\beq \label{d1077}
||\mathcal{D}S_{\xi_{\b}} \mathcal{D}^{-1}S_{\xi_{\a}}f||_{L^{2}} \leq C_{N} (2^{i}|\xi_{\b}-\xi_{\a}|)^{-N} ||S_{\xi_{\a}}f||_{L^{2}}
\eeq

\end{l1}

\begin{proof}

We start with proving (\ref{d1000}) for $N=1$. We have:

\beq \label{d565}
(D_{x_{1}}-\xi^{1}_{\a}) \mathcal{D}S_{\xi_{\a}}f = \mathcal{D}(D_{x_{1}}-\xi^{1}_{\a}) S_{\xi_{\a}} f + [D_{x_{1}}-\xi^{1}_{\a}, \mathcal{D}] S_{\xi_{\a}}f
\eeq

We deal separately with each term on the right hand side of (\ref{d565}). The standard calculus gives us that the symbol of $[D_{x_{1}}-\xi^{1}_{\a}, \mathcal{D}]$ is:

$$ r(x,\xi,\tau)= - i \frac{\partial}{\partial x_{\a}} d(x,\xi,\tau) =  -i  \frac{2 x_{1}}{\mu +|\tau|+\xi^{2}} (1+\frac{x^{2}}{\mu +|\tau|+\xi^{2}})^{-1} d(x,\xi,\tau)$$

$r(x,\xi,\tau)$ is hypoelliptic of one class smoother than $d(x,\xi,\tau)$ in the following sense:

$$|\partial_{\xi}^{\a} \partial_{x}^{\b} r(x,\xi,\tau)| \leq C_{\a \b} |d(x,\xi)| \langle \sqrt{\mu+|\tau|+ \xi^{2}} \rangle^{-1-|\a|}$$

Then emulating a similar argument as in the proof of (\ref{dec22}) we obtain:

\beq \label{d566}
||[D_{x_{1}}-\xi^{1}_{\a}, \mathcal{D}] S_{\xi_{\a}}f||_{L^{2}} \leq 2^{-i} ||\langle \frac{x^{2}}{\mu+2^{2i}} \rangle^{\frac{1}{4}+\frac{\e}{2}} S_{\xi_{\a}}f|| \approx 2^{-i} ||\mathcal{D} S_{\xi_{\a}}f||
\eeq

For the second term on the right hand side of (\ref{d565}) we want to prove a similar estimate:

\beq \label{d567}
||\mathcal{D}(D_{x_{1}}-\xi^{1}_{\a}) S_{\xi_{\a}} f||_{L^{2}} \leq 2^{-i} ||\mathcal{D} S_{\xi_{\a}} f||_{L^{2}}
\eeq

The underlying idea is that in the support of $S_{\xi_{\a}}$ we have that $|\xi^{1}-\xi^{1}_{\a}| \leq 2^{-i}$ ($\xi^{1}-\xi^{1}_{\a}$ is the symbol of $D_{x_{1}}-\xi^{1}_{\a}$). We need an argument just a bit more subtle. Let $(Q_{i}^{m})_{m \in \Z^{2}}$ be a system of of cubes centered at $(2^{i}+\mu)m$ and of sizes $2^{i}+\mu$ which form a partition of $\R^{2}$. Let also $\chi_{Q_{i}^{m}}$ to be a smooth approximation of the characteristic function of $Q_{i}^{m}$ and such that $(\chi_{Q_{i}^{m}})_{m \in \Z^{2}}$ forms a partition of unity in $R^{2}$. We decompose:

$$S_{\xi_{\a}} f = \sum_{m \in \Z^{2}} \chi_{Q_{i}^{m}} S_{\xi_{\a}} f$$

If $\tilde{S}_{\xi_{\a}}$ is a localizing operator similar to $S_{\xi_{\a}}$ in the sense that $\tilde{S}_{\xi_{\a}} S_{\xi_{\a}}=S_{\xi_{\a}}$ and the support of $\tilde{S}_{\xi_{\a}}$ is a cube of size $2 \times 2^{-i}$ centered at $\xi_{\a}$, then:

$$(D_{x_{1}}-\xi^{1}_{\a}) S_{\xi_{\a}} f = \sum_{m' \in \Z^{2}} \sum_{m \in \Z^{2}} \chi_{Q^{m'}_{i}} \tilde{S}_{\xi_{\a}} (D_{x_{1}}-\xi^{1}_{\a})  \chi_{Q_{i}^{m}} S_{\xi_{\a}} f$$

 $\tilde{S}_{\xi_{\a}} (\xi^{1}-\xi^{1}_{\a})$ is localized on a scale $2^{-i} \times 2^{-i}$ which is greater or equal than the dual scale produced by localizations corresponding to $\chi_{Q_{i}^{m}}$. Hence a standard argument gives us:
 
 $$||\chi_{Q^{m'}_{i}} \tilde{S}_{\xi_{\a}} (D_{x_{1}}-\xi^{1}_{\a})  \chi_{Q_{i}^{m}} S_{\xi_{\a}} f||_{L^{2}} \leq $$
 
 $$C_{N} \langle m-m' \rangle^{-N} ||\tilde{S}_{\xi_{\a}} (D_{x_{1}}-\xi^{1}_{\a}) \chi_{Q_{i}^{m}} S_{\xi_{\a}} f||_{L^{2}} \leq C_{N} 2^{-i} \langle m-m' \rangle^{-N} || \chi_{Q_{i}^{m}} S_{\xi_{\a}} f||_{L^{2}}$$
 
 In the last line we have use the above observation that $||\tilde{S}_{\xi_{\a}} (D_{x_{1}}-\xi^{1}_{\a}) h||_{L^{2}} \leq 2^{-i} ||h||_{L^{2}}$. We can continue with:
 
\begin{align*}
& ||\mathcal{D}(D_{x_{1}}-\xi^{1}_{\a}) S_{\xi_{\a}} f||^{2}_{L^{2}} \approx ||\langle \frac{x}{\mu+2^{i}} \rangle^{\q+\e} (D_{x_{1}}-\xi^{1}_{\a}) S_{\xi_{\a}} f||^{2}_{L^{2}} \\ 
&\approx \sum_{m' \in \Z^{2}} \langle m' \rangle^{1+2\e} ||\chi_{Q^{m'}_{i}} \tilde{S}_{\xi_{\a}} (D_{x_{1}}-\xi^{1}_{\a}) \sum_{m} \chi_{Q_{i}^{m}}  S_{\xi_{\a}} f||^{2}_{L^{2}} \\
&\leq \sum_{m' \in \Z^{2}} \langle m' \rangle^{1+2\e} \left( \sum_{m \in \Z^{2}} C_{N} 2^{-i} \langle m-m' \rangle^{-N} || \chi_{Q_{i}^{m}} S_{\xi_{\a}} f||_{L^{2}}   \right)^{2} \\
&\leq C^{2}_{N}2^{-2i} \sum_{m' \in \Z^{2}} \sum_{m \in \Z^{2}} \langle m' \rangle^{1+2\e} \langle m-m' \rangle^{-2N+4} || \chi_{Q_{i}^{m}} S_{\xi_{\a}} f||^{2}_{L^{2}}  \\
&\leq C^{2}_{N}2^{-2i} \sum_{m \in \Z^{2}} \sum_{m' \in \Z^{2}} \langle m \rangle^{1+2\e}  \langle m-m' \rangle^{-2N+6} || \chi_{Q_{i}^{m}} S_{\xi_{\a}} f||^{2}_{L^{2}} \\
&\leq C^{2}_{N} 2^{-2i} \sum_{m \in \Z^{2}} \langle m \rangle^{1+2\e} ||\chi_{Q_{i}^{m}} S_{\xi_{\a}} f||^{2}_{L^{2}} \approx C^{2}_{N} 2^{-2i} ||\mathcal{D} S_{\xi_{\a}} f||^{2}_{L^{2}}\\
\end{align*}

This finishes the proof of (\ref{d567}). From (\ref{d565}),(\ref{d566}) and (\ref{d567}) we obtain:

$$||(D_{x_{1}}-\xi^{1}_{\a}) \mathcal{D} S_{\xi_{\a}} f||_{L^{2}} \leq 2^{-i} ||\mathcal{D} S_{\xi_{\a}} f||_{L^{2}}$$

In a similar manner we obtain a similar estimate when $D_{x_{1}}-\xi^{1}_{\a}$ is replaced by $D_{x_{2}}-\xi^{2}_{\a}$ and these two estimates together are enough to justify (\ref{d1000}) in the case $N=1$. The argument for general $N$ follows in a similar manner.

Now we turn our attention to proving (\ref{d1077}). The proof is similar in spirit to the argument for (\ref{d1000}), just that it involves more computations. This is why we outline only the main steps. Let us assume that $|\xi_{\b}| \approx |\xi_{\a}| \approx 2^{i}$. We decompose:

$$\mathcal{D}S_{\xi_{\b}}\mathcal{D}^{-1} S_{\xi_{\a}}=\sum_{m \in \Z^{2}} \sum_{\bar{m} \in \Z^{2}} \mathcal{D}S_{\xi_{\b}} \chi_{Q^{i}_{\bar{m}}} \mathcal{D}^{-1} S_{i} \chi_{Q^{m}_{i}} S_{\xi_{\a}}$$

The hard part of the argument is to prove that:

$$||\mathcal{D}S_{\xi_{\b}} \chi_{Q_{i}^{\bar{m}}} \mathcal{D}^{-1} S_{i} \chi_{Q_{i}^{m}} S_{\xi_{\a}}f||_{L^{2}} \leq C_{N} (2^{i}|\xi_{\b}-\xi_{\a}|)^{-N} \langle m-\bar{m} \rangle^{-N} ||\chi_{Q_{i}^{m}} S_{\xi_{\a}}f||_{L^{2}}$$

This can be achieved by using the results (and the arguments used in their proofs) in (\ref{dec22}) and (\ref{d1000}). Then we can put the above estimates together and sum them up to obtain (\ref{d1077}). A similar approach would give the estimate in the case $|\xi_{\b}| \approx |\xi_{\a}|$.

\end{proof}

We prove in detail an easier variant of (\ref{d1111}) in the sense that we evaluate all terms in $X^{s,\q}$: 

\begin{l1}
We have:

\beq \label{d3001}
||\mathcal{D}f||^{2}_{X^{s,\q}} \approx   \sum_{i} \sum_{\xi \in \Xi^{i,*}}  \sum_{l \in \Z: |(\xi,l)| \approx 2^{i}} \langle (\xi,l)\rangle^{2s} \langle l-\xi^{2} \rangle ||\langle \frac{x}{2^{i}} \rangle^{\q+\e} f_{\xi,l} ||^{2}_{L^{2}} 
\eeq

\end{l1}

\begin{proof}

We have the decomposition:

$$\mathcal{D}f=\sum_{i} \sum_{\xi \in \Xi^{i,*}}  \sum_{l \in \Z: |(\xi,l)| \approx 2^{i}} \mathcal{D}  f_{\xi,l}$$

A direct consequence of (\ref{d1000}) is that:

$$||(\mathcal{D}  f_{\xi,l})_{\eta,k}||_{L^{2}} \leq C_{N} \langle |\xi||(\xi,l)-(\eta,k)|\rangle^{-N} ||\mathcal{D}  f_{\xi,l}||_{L^{2}}$$

\noindent 
and this is enough to justify:

$$||\mathcal{D}f||^{2}_{X^{s,\q}} \leq   \sum_{i} \sum_{\xi \in \Xi^{i,*}}  \sum_{l \in \Z: |(\xi,l)| \approx 2^{i}} \langle (\xi,l)\rangle^{2s} \langle l-\xi^{2} \rangle ||\langle \frac{x}{2^{i}} \rangle^{\q+\e} f_{\xi,l} ||^{2}_{L^{2}} 
$$

If we write $\mathcal{D}f=g$, the reverse inequality is equivalent to:

$$  \sum_{i} \sum_{\xi \in \Xi^{i,*}}  \sum_{l \in \Z: |(\xi,l)| \approx 2^{i}} \langle (\xi,l)\rangle^{2s} \langle l-\xi^{2} \rangle ||\mathcal{D} (\mathcal{D}^{-1} g)_{\xi,l} ||^{2}_{L^{2}} \leq ||g||_{X^{s,\q}}$$

This can be easily deduced from:

$$||\mathcal{D} (\mathcal{D}^{-1} g_{\eta,k})_{\xi,l}||_{L^{2}} \leq C_{N} \langle |\xi||(\xi,l)-(\eta,k)|\rangle^{-N} ||g_{\eta,k}||_{L^{2}}$$

\noindent
which is a direct consequence of (\ref{d1077})

\end{proof}

\begin{proof}[Proof of Proposition \ref{dec}]

It can be easily seen that in the argument for (\ref{d3001}) we can replace $X^{s,\q}$ with $X^{s,\q,1}_{\geq 1} \cap X^{s,\q,\infty}_{\leq 1}$. 

In addition to what has been done so far, we need only to embed the $Y$ structure in the computations. We need results similar to (\ref{d1000}) and (\ref{d1077}). For instance, the equivalent of (\ref{d1000}) is: 

\beq \label{d1001}
||(\mathcal{D}f_{\xi, \leq 1})_{\eta, \leq 1}||_{Y} \leq C_{N} \langle \xi-\eta \rangle^{-N} ||\mathcal{D}f_{\xi, \leq 1}||_{Y}
\eeq

We can copy verbatim the argument we provided for the proof of (\ref{d1000}) with replacing $L^{2}$ with $L^{\infty}_{t}L^{2}_{x}$ and work on localized tubes in the physical space. Just that in this way we end up measuring $ (\mathcal{D} f_{\xi, \leq 1})_{\eta, \leq 1}$ in $L^{\infty}_{t}L^{2}_{x}(T_{\eta}^{m,l})$, instead of measuring it in $L^{\infty}_{t}L^{2}_{x}(T_{\xi}^{m,l})$. 

One the other hand, we  see a difference between the system $(T_{\eta}^{m,l})_{(m,l) \in \Z^{3}}$ and its associated norm and the system $(T_{\xi}^{m,l})_{(m,l) \in \Z^{3}}$ and its associated norm only if $|\xi -\eta| \geq 1$. From the projected gain of a factor $\langle \xi -\eta \rangle^{-N}$ we can easily spare a factor of $ \langle \xi-\eta \rangle^{-1}$ to make this transition. 

Notice that the factor $\langle \xi-\eta \rangle^{-N}$ was $ \langle |\xi|(\xi-\eta) \rangle^{-N}$ in the context of (\ref{d1000}). The reason is that in the proof of (\ref{d1000}) we deal with $L^{2}$ theory while now we deal with $L^{\infty}L^{2}$ functions. The spared factors of $|\xi|^{-\q}$ is enough to make the transition between $L^{\infty}_{t}L^{2}_{x}$ and $L^{2}_{x,t}$. We leave out the rest of the details of the proof. 

\end{proof}

\begin{proof}[Proof of (\ref{e3}) and (\ref{e4}) with decay]
If we wanted to provide a complete proof of how decay is preserved in (\ref{e3}) and (\ref{e4}) we should rewrite many of the computations we made for the proofs without decay. We choose instead to outline only the main steps and leave the computations aside.

\vspace{.1in}

\it{Step 1.} \rm If $\xi_{0} \in \Xi^{*}$ and $l_{0} \in \Z$ such that $|l_{0}-\xi_{0}^{2}| \geq 1$, then:

$$\hat{w}_{\xi_{0},l_{0}}(\xi,\tau)=\frac{\hat{f}_{\xi_{0},l_{0}}(\xi,\tau)}{\tau-\xi^{2}}$$

Therefore $w_{\xi_{0},l_{0}}=f_{\xi_{0},l_{0}} * \mathcal{F}^{-1}(\frac{\tilde{\phi}_{\xi_{0},l}}{\tau-\xi^{2}})$ where $\tilde{\phi}_{\xi_{0},l}$ is a smooth approximation of the characteristic function of the support of the operator $T_{\xi_{0},l_{0}}$ and $\hat{f}_{\xi_{0},l_{0}} \cdot \tilde{\phi}_{\xi_{0},l}= \hat{f}_{\xi_{0},l_{0}}$ for any $f$. A simple computation shows that $\mathcal{F}^{-1}(\frac{\tilde{\phi}_{\xi_{0},l}}{\tau-\xi^{2}})$ is concentrated in a cube of sizes $\langle \xi \rangle \times \langle \xi \rangle \times 1$ around the origin and that $||\mathcal{F}^{-1}(\frac{\tilde{\phi}_{\xi_{0},l}}{\tau-\xi^{2}})||_{L^{1}} \leq |l_{0}-\xi_{0}^{2}|^{-1}$. This is enough to justify that:

 $$||\mathcal{D}w_{\xi_{0},l_{0}}||_{L^{2}} \leq |l_{0}-\xi_{0}|^{-1}||\mathcal{D}f_{\xi_{0},l_{0}}||_{L^{2}}$$
 
 Using the results from Proposition \ref{dec} we can get the conservation of decay at the $X^{s,\q}$ level.
 
 \vspace{.1in}
 
 \it{Step 2.} \rm In dealing with the $Y^{s,N}$ structure we notice that the weights coming from the decay we impose are constant on the tubes $T_{\xi}^{m,l}$. Then it is a routine exercise to show that they can be easily absorbed in the computations in the proofs of Propositions (\ref{e3}) and (\ref{e4}). 

\vspace{.1in}
 
 \it{Step 3.} \rm Last thing to prove is that the homogeneous equation preserves the decay condition in $X^{s,\q,1}$. This is not that straightforward since when we argued about the $X^{s,\q,1}$ structure without decay, we did not involve any space decomposition. This time this should be done in the same spirit with the decomposition used for proving the $Y^{s,N}$ structure for the homogeneous equation. The localizations in the physical space do not have to be at the level of $T_{\xi}^{m,l}$, but rather at the level of $Q^{m} \times [l,l+1]$ on which the decay is like a constant. Then we go on and argue in a similar manner as before.

\end{proof}

The last result we provide in this section is a Lemma which will be useful in proving the conservation of decay in the bilinear estimates. For any $j \in \Z$, let us denote by $d_{j}(x)=(1+\frac{|x|^{2}}{\mu+2^{2j}})^{\frac{1}{4}+\frac{\e}{2}}$ to be the weight corresponding to $\mathcal{D}$ localized at frequency $2^{j}$. We have the result:

\begin{l1} For any function $f \in L^{2}$ we have:

\beq \label{dec50}
||d_{k}S_{k}f||_{L_{2}} \leq C_{N} 2^{(k-j)(\q+\e)} \sum_{i \in \N} 2^{-|i-k|N} ||S_{i}d_{j}f||_{L^{2}}
\eeq

\begin{proof}

We write:

$$d_{k}S_{k}f=d_{k}S_{k}d_{j}^{-1}d_{j}f=\sum_{i}d_{k}S_{k}d_{j}^{-1}S_{i}d_{j}f$$

We denote by $h=d_{j}f$. In the same spirit of argument for (\ref{d1077}) we can prove the following:

$$||d_{k}S_{k}d_{j}^{-1}S_{i}h||_{L^{2}} \leq C_{N} 2^{(k-j)(\q+\e)} 2^{-|i-k|N} ||h||_{L^{2}}$$

The factor $2^{(k-j)(\q+\e)}$ appears here because we freeze the weights $d_{k}$ and $d_{j}^{-1}$ unlike in (\ref{d1077}). If $i=k$, then it is almost like estimating $||d_{k}d_{j}^{-1}h||_{L^{2}}$ and here it is obvious where the factor comes from, since $\frac{d_{k}}{d_{j}} \approx 2^{(k-j)(\q+\e)}$.

\end{proof}

\end{l1}

\section{Bilinear estimates in $\mathcal{R}X^{s,\q,1}$} \label{best}

The objective of this section is to obtain the bilinear estimates for $B(u,v)$ and $B(u,\bar{v})$ in $\mathcal{R}X^{s,\q,1}$, where $B$ is of type (\ref{bf}). We introduce the additional bilinear form $\tilde{B}(u,v)=u \cdot v$. If $\hat{u}$ is localized in $A_{i}$ we use the estimate $||\nabla{u}||_{L^{2}} \leq 2^{i} ||u||_{L^{2}}$. $X^{s,\pm \q,1}$ are $L^{2}$ like on dyadic pieces, hence if $\hat{u}$ is localized in $A_{i}$ and $\hat{v}$ is localized in $A_{j}$ we use the estimates:

\beq \label{bb}
 ||\tilde{B}(u,v)||_{X} \leq C ||u||_{X'} ||v||_{X''}   \Ra  ||B(u,v)||_{X} \leq 2^{i+j} C ||u||_{X'} ||v||_{X''}
\eeq

\beq \label{bbb}
 ||\tilde{B}(u,\bar{v})||_{X} \leq C ||u||_{X'} ||v||_{X''}   \Ra  ||B(u,\bar{v})||_{X} \leq 2^{i+j} C ||u||_{X'} ||v||_{X''}
\eeq

Here $X,X',X''$ are of type $X^{s,\pm \q, a}$ ($a \in \{1,2,\infty\}$). The constant $C$ may depend on $u,v$, more exactly of their localizations. The key thing is once we have estimates for $\tilde{B}$, we obtain estimates for $B$ by simply bringing in the correction factor of $2^{i+j}$.  

If $B$ were of type (\ref{bf1}) the correction factor would be only $2^{j}$ and this justifies why we can claim the estimates for bilinear estimates of type (\ref{bf1}).

The main results we claim are listed in the following theorem.  

\begin{t1} \label{tb1}

a)  Assume that $i \leq j$. We have the bilinear estimate:

\beq \label{b1}
||B(u,v)||_{\mathcal{R}X^{s,-\q,1}_{k}} \leq j^{\frac{3}{2}} 2^{(1-s)i} 2^{(k-j)s} ||u||_{\mathcal{R}X^{s,\q,1}_{i}} ||v||_{\mathcal{R}X^{s,\q,1}_{j}}
\eeq

b) Assume that $5i \leq j$. We have the following bilinear estimates:

\beq \label{b4}
||B(u,v)||_{\mathcal{R}X^{s,-\q,1}_{j, \geq 2^{-i}}} \leq i^{\frac{3}{2}} 2^{(1-s)i} ||u||_{\mathcal{R}X^{s,\q,1}_{i}} ||v||_{\mathcal{R}X^{s,\q,1}_{j, \geq 2^{-i}}}
\eeq

Both estimates (\ref{b1}) and (\ref{b4}) remain valid if $B(u,v)$ is replaced by $B(\bar{u},v)$ or $B(u,\bar{v})$. 

\end{t1}

The theorem tells us that at the dyadic level the spaces $\mathcal{R}X^{s,\q,1}$ with $1 < s$ are suitable for the bilinear estimates as long as $j$ can be controlled by a power of $2^{i}$ . For instance the estimates are good as long as $j \leq Ci$. 

If this does not happen, the second part of the Theorem tells us that $\mathcal{R}X^{s,\q,1}$ is still good enough to measure the high frequency at distance greater than $2^{-i}$ from the parabola. This means we need additional information close to $P$. 

\vspace{.2in}

\subsection{Basic Estimates}
\noindent

\vspace{.1in}

We start with a simple result stating how two parabolas interact under convolution. We need few technical definitions. 

Throughout this section functions are defined on Fourier space (they should be thought as Fourier transforms). This is why we use the standard coordinates $(\xi,\tau)$. Also the operator $S^{\xi}_{k}$ should be understood as $S^{\xi}_{k}f=\varphi_{k}(|\xi|) \cdot f(\xi,\tau)$.

For each $c \in \R$ denote by $P_{c}= \{(\xi,\tau): \tau-\xi^{2}=c\}$ and  by $\bar{P}_{c}=\{(\xi,\tau): \tau+\xi^{2}=c \}$. For simplicity $P=P_{0}$ and $\bar{P}=\bar{P}_{0}$. 

Denote by $\delta_{P_{c}}=\delta_{\tau-\xi^{2}=c}$ the standard surface measure associated to the parabola $P_{c}$. With respect to this measure, the restriction of $f:\R^{3} \rightarrow \C$ to $P_{c}$ has norm:

$$||f||_{L^{2}(P_{c})}=\left( \int |f|^{2}(\xi,\xi^{2}+c) \sqrt{1+4|\xi|^{2}} d\xi  \right)^{\q}$$

Throughout this section we make the following convention:

$$i \cong j \ \ \mbox{if} \ \  |i-j| \leq 3 \ \ \mbox{and} \ \ i \ncong j \ \ \mbox{if} \ \ |i-j| \geq 4$$

The estimates for the case $i=j$ are generic for the case $i \cong j$. Hence, we will list or prove only the case $i=j$ since the other ones are similar and can be proved the same way.

We want estimates for $f \delta_{P} * g \delta_{P}$ restricted to $P_{c}$. We assume $f$ is localized at $2^{i}$ and $g$ at $2^{j}$. We want to measure the part of the restriction of the convolution which is localized at $2^{\max{(i,j)}}$. We obtain a good result for the those $c$ with $|c| \leq 2^{i+j-2}$, i.e. at distance at most $2^{\min{(i,j)}-2}$ from $P$, while for the rest of the interaction we provide only a global $L^{2}$ estimate.

The next Proposition states the main ingredient for the bilinear estimates.

\begin{p2} \label{bep1}

a) Let $f,g \in \mathcal{R}L^{2}(P)\  \mbox{or} \ \mathcal{R}L^{2}(\bar{P})$ such that $f$ is localized at $|\xi| \approx 2^{i}$ and $g$ at $|\xi| \approx 2^{j}$. If $\min{(i,j)} \geq 1$ and $|c| \leq 2^{i+j-2}$ we have

\beq \label{gg} 
||S_{\max{(i,j)}}^{\xi}(f \delta_{P^{1}} * g \delta_{P^{2}}) ||_{L^{2}(P_{c})} \leq 2^{\frac{\min{(i,j)}}{2}} ||f||_{\mathcal{R}L^{2}(P^{1})}||g||_{L^{2}(P^{2})} 
\eeq

\noindent
where $P^{1}, P^{2} \in \{P,\bar{P}\}$, except for the case $P_{1}=P_{2}=\bar{P}$.

b) Also we have the global estimate

\beq \label{ge}
||f \delta_{P^{1}} * g \delta_{P^{2}} ||_{L^{2}} \leq 2^{\min{(i,j)}} ||f||_{L^{2}(P^{1})}||g||_{L^{2}(P^{2})}
\eeq
 
\noindent
where $P^{1}, P^{2} \in \{P,\bar{P}\}$, except for the case $P_{1}=P_{2}=\bar{P}$.

\end{p2}

Remark. In the above Proposition we do not provide the ingredients for high-high interactions with outcome at low frequencies (see the $S^{\xi}_{\max{(i,j)}}$). These type of estimates are going to be proved later using duality.  

\begin{proof}
$a)$ {\mathversion{bold} $f \delta_{P} * g \delta_{P}$}.

We make the choice $i \leq j$. We have to prove two results according to the case when we use rotations on the high or the low frequency.

We simplify also the arguments of $f$ and $g$. For $f \in L^{2}(P)$, instead of using the full argument $f(\xi, \xi^{2})$ we reduce it to $f(\xi)$, where $\xi \in \R^{2}$. We can also use polar coordinates $\xi=(\r,\t) \in (0,\infty) \times [0,2 \pi)$  and denote by $f(\r,\t)=f(\xi)$. With these reductions we have:

$$||f||^{2}_{L^{2}(P)}=\int |f(\xi)|^{2} \sqrt{1+4\xi^{2}} d\xi= \int |f(\r,\t)|^{2} \sqrt{1+4\r^{2}} \r d\r d\t$$
 
Because $f$ is localized at $|\xi| \approx 2^{i}$, the domain of integration for the $\r$ variable is $(2^{i-1},2^{i+1})$. We have a similar formula for $g$ with the corresponding observation about the domain of integration for the $\r$ variable. 

We want to estimate the convolution of two measures, hence we need to derive  the formula which gives us the value of $f \delta_{\tau=\xi^{2}} * g \delta_{\tau=\xi^{2}}$ at $(\zeta_{1}, \zeta_{2},\zeta_{3})$. If $h:\R^{3} \rightarrow \C$ is a smooth function which decays rapidly at $\infty$, then from the definition of convolution we have:

$$(f \delta_{\tau=\xi^{2}} * g \delta_{\tau=\xi^{2}})h=$$

$$\int f(\xi) g(\eta) h(\xi+\eta, \xi^{2}+\eta^{2}) \sqrt{1+4\xi^{2}}  \sqrt{1+4\eta^{2}} d\xi d\eta =$$

$$ \int f(\r_{1}, \t_{1}) g(\r_{2}, \t_{2}) h(\r_{1} \cos{\t_{1}}+\r_{2} \cos{\t_{2}},\r_{1} \sin{\t_{1}}+\r_{2} \sin{\t_{2}}, \r_{1}^{2}+\r_{2}^{2}) $$

$$\sqrt{1+4\r_{1}^{2}} \sqrt{1+4\r_{2}^{2}} \r_{1} \r_{2} d{\r_{1}} d \r_{2} d\t_{1} d\t_{2} $$

If we write $(f \delta_{\tau=\xi^{2}} * g \delta_{\tau=\xi^{2}})h=\int l(\z) h(\z) d \z$ this gives us $f \delta_{\tau=\xi^{2}} * g \delta_{\tau=\xi^{2}}=l$. Motivated by this, we introduce the change of variables $\r_{2}(\r_{1}), \t_{1}(\r_{1}),\t_{2}(\r_{1}) \rightarrow (\z_{1},\z_{2},\z_{3})$:

\begin{eqnarray} \label{sys}
  \left\{
        \begin{array}{rl}
		 \r_{1}\cos{\t_{1}}+\r_{2}\cos{\t_{2}}=\z_{1}  \\
		 \r_{1}\sin{\t_{1}}+\r_{2}\sin{\t_{2}}=\z_{2}  \\
		 \r_{1}^{2}+\r_{2}^{2}=\z_{3}
	\end{array}\right.
\end{eqnarray}

\noindent
where $\r_{1}$ is seen as a parameter. Computing the Jacobian of the transformation, we obtain $d \z_{1} d \z_{2} d \z_{3}= 2\r_{1} \r_{2}^{2} |\sin{(\t_{1}-\t_{2})}| d \r_{2} d \t_{1} d \t_{2}$, therefore:

$$(f \delta_{\tau=\xi^{2}} * g \delta_{\tau=\xi^{2}}) (\zeta_{1}, \zeta_{2},\zeta_{3}) =$$

$$ \int f(\r_{1}, \t_{1}) g(\r_{2}, \t_{2}) \sqrt{1+4\r_{1}^{2}} \sqrt{1+4\r_{2}^{2}} \frac{d \r_{1}}{2\r_{2} |\sin{(\t_{1}-\t_{2})}|}$$

\noindent
where $\r_{2}(\r_{1}), \t_{1}(\r_{1}),\t_{2}(\r_{1})$ solve (\ref{sys}) with parameter $\r_{1}$.

Since we evaluate the result on $P_{c}$, we are interested in the points satisfying $\z_{3}=\z^{2}_{1}+\z^{2}_{2}+c$. In polar coordinates, this condition becomes  $2\r_{1}\r_{2}\cos{(\t_{1}-\t_{2})}=-c$. Then we have:

$$\frac{\sqrt{1+4\r_{1}^{2}} \sqrt{1+4\r_{2}^{2}}}{\r_{1} \r_{2} |\sin{(\t_{1}-\t_{2})}|} = \frac{\sqrt{1+4\r_{1}^{2}} \sqrt{1+4\r_{2}^{2}}}{\sqrt{\r^{2}_{1} \r_{2}^{2}-\frac{c^{2}}{4}}} \approx 4 $$

\noindent
since $|c| \leq 2^{i+j-2}$ and $2^{i-1} \leq \r_{1}$, $2^{j-1} \leq \r_{2}$ . We obtain:

$$|f \delta_{\tau=\xi^{2}} * g \delta_{\tau=\xi^{2}} (\zeta_{1}, \zeta_{2},\zeta^{2}_{1}+\zeta^{2}_{2}+c)| \leq \int |f(\r_{1}, \t_{1})| |g(\r_{2}, \t_{2})|  \r_{1} d \r_{1} $$

Next we estimate $f \delta_{\tau=\xi^{2}} * g \delta_{\tau=\xi^{2}}$ in $L^{2}(P_{c})$. From (\ref{sys}) we can easily derive that $d\z_{1}d\z_{2} = \r_{2}d\r_{2}d\t_{2}$. Notice that the third equation in the system gave us $2\r_{1}\r_{2}\cos{(\t_{1}-\t_{2})}=-c$ from where we can express $\t_{1}$ in terms of $\r_{1},\r_{2},\t_{2}$. Therefore we have:

$$||f \delta_{\tau=\xi^{2}} * g \delta_{\tau=\xi^{2}}||^{2}_{L^{2}(P_{c})} $$

$$\leq \int \left( \int |f(\r_{1}, \t_{1})| \cdot |g(\r_{2}, \t_{2})|  \r_{1} d \r_{1} \right)^{2} \sqrt{1+4(\z^{2}_{1}+\z^{2}_{2})} d \z_{1} d \z_{2}$$

$$\leq 2^{j}  \left( \int \sup_{\t_{1}} |f(\r_{1},\t_{1})|^{2} \r_{1}d\r_{1} \right) \left( \int |g(\r_{2},\t_{2})|^{2}  \r_{1} d \z_{1} d \z_{2} d \r_{1} \right) $$

$$\approx \left( \int \sup_{\t_{1}} |f(\r_{1},\t_{1})|^{2} \r_{1}d\r_{1}  \right)  \left( \int |g(\r_{2},\t_{2})|^{2} \sqrt{1+4\r_{2}^{2}} \r_{2}\r_{1} d \r_{2}  d \r_{1} \right) $$

$$\leq 2^{-i}  \left( \int \sup_{\t_{1}} |f(\r_{1},\t_{1})|^{2} \sqrt{1+4\r_{1}^{2}} \r_{1}d\r_{1} \right) \left( \int ||g||^{2}_{L^{2}(P)} \r_{1} d \r_{1} \right) $$

$$ \leq 2^{i} ||f||^{2}_{\mathcal{R}L^{2}} ||g||^{2}_{L^{2}(P)}$$

In the last line we have taken advantage of the rotations via the estimate:

$$ \int \sup_{\t_{1}} |f(\r_{1},\t_{1})|^{2} \sqrt{1+4\r_{1}^{2}} \r_{1}d\r_{1} \leq ||f||^{2}_{\mathcal{R}L^{2}}$$

If we want to use the rotations on $g$, then we use the change of variables $(\z_{1}, \z_{2}) \rightarrow (\r_{2}, \t_{1})$ satisfying $d\z_{1}d\z_{2} = \r_{2}d\r_{2}d\t_{1}$. We obtain:

$$||f \delta_{\tau=\xi^{2}} * g \delta_{\tau=\xi^{2}}||^{2}_{L^{2}(P_{c})} \leq$$

$$\int \left( \int |f(\r_{1}, \t_{1})| \cdot |g(\r_{2}, \t_{2})|  \r_{1} d \r_{1} \right)^{2} \sqrt{1+4(\z^{2}_{1}+\z^{2}_{2})} d \z_{1} d \z_{2} \leq$$

$$2^{j} \left( \int |f(\r_{1},\t_{1})|^{2} \r_{1}d\r_{1} d\t_{1} \right) \left( \int \sup_{\t_{2}} |g(\r_{2},\t_{2})|^{2}  \r_{1} \r_{2} d\r_{2} d \r_{1} \right) \leq$$

$$ 2^{-i} ||f||^{2}_{L^{2}(P)} \left( \int \sup_{\t_{2}} |g(\r_{2},\t_{2})|^{2} \sqrt{1+4\r_{2}^{2}} \r_{1} \r_{2} d \r_{2}  d \r_{1} \right) \leq$$

$$2^{-i} ||f||^{2}_{L^{2}(P)} \left( \int ||g||^{2}_{\mathcal{R}L^{2}(P)} \r_{1} d \r_{1} \right) \leq 2^{i} ||f||^{2}_{L^{2}} ||g||^{2}_{\mathcal{R}L^{2}(P)} $$

\vspace{.1in}

{\mathversion{bold} $f \delta_{\bar{P}} * g \delta_{P}$}

\vspace{.1in}

We pursue the same argument we used for $f \delta_{P} * g \delta_{P}$.

The value of $f \delta_{\tau=-\xi^{2}} * g \delta_{\tau=\xi^{2}}$ at $(\zeta_{1}, \zeta_{2},\zeta_{3})$ is given by: 

$$(f \delta_{\tau=-\xi^{2}} * g \delta_{\tau=\xi^{2}}) (\zeta_{1}, \zeta_{2},\zeta_{3}) =$$

$$ \int f(\r_{1}, \t_{1}) g(\r_{2}, \t_{2}) \sqrt{1+4\r_{1}^{2}} \sqrt{1+4\r_{2}^{2}} \frac{d \r_{1}}{\r_{2} |\sin{(\t_{1}-\t_{2})}|}$$

\noindent
where $\r_{2}(\r_{1}), \t_{1}(\r_{1}),\t_{2}(\r_{1})$ solve the system with parameter $\r_{1}$:

\begin{eqnarray} \label{sys2}
  \left\{
        \begin{array}{rl}
		 \r_{1}\cos{\t_{1}}+\r_{2}\cos_{\t_{2}}=\z_{1}  \\
		 \r_{1}\sin{\t_{1}}+\r_{2}\sin_{\t_{2}}=\z_{2}  \\
		 -\r_{1}^{2}+\r_{2}^{2}=\z_{3}
	\end{array}\right.
\end{eqnarray}

We want to evaluate the result on $P_{c}$ therefore we are interested in the points satisfying $\z_{3}=\z^{2}_{1}+\z^{2}_{2}+c$. In polar coordinates, this condition becomes  $2\r_{1}\r_{2}\cos{(\t_{1}-\t_{2})}=-c-2\r^{2}_{1}$. Since we localize the interaction at $2^{j}$, we have $\r_{2}^{2}-\r_{1}^{2} \geq 2^{2j-2}$. Then we can estimate:

$$\frac{\sqrt{1+4\r_{1}^{2}} \sqrt{1+4\r_{2}^{2}}}{\r_{1} \r_{2} |\sin{(\t_{1}-\t_{2})}|} = \frac{\sqrt{1+4\r_{1}^{2}} \sqrt{1+4\r_{2}^{2}}}{\sqrt{\r^{2}_{1} \r_{2}^{2}-\frac{(c+2 \r_{1}^{2})^{2}}{4}}} \leq 8 $$

From this point on we can continue like in the case of $f \delta_{\tau=\xi^{2}} * g \delta_{\tau=\xi^{2}}$.

b) We assume $i \leq j$ and we have to consider two cases.

\begin{bfseries} Case 1: {\mathversion{bold}  $i \ncong j$} \end{bfseries} 

{\mathversion{bold} $f \delta_{P} * g \delta_{P}$}

The strategy here is to prove an estimate of type:

\begin{equation} \label{eq1}
||f \delta_{\tau=\xi^{2}} * g \delta_{\tau=\xi^{2}} * h||_{L^{\infty}} \leq 2^{i} ||f||_{L^{2}(P)}||g||_{L^{2}(P)} ||h||_{L^{2}} 
\end{equation}

\noindent
for any $h \in L^{2}(\R^{3})$. We have:

$$(f \delta_{\tau=\xi^{2}} * g \delta_{\tau=\xi^{2}} * h)(z_{1},z_{2},z_{3})=$$

$$ \int h(z_{1}-\xi_{1}-\eta_{1},z_{2}-\xi_{2}-\eta_{2}, z_{3} - \xi^{2}-\eta^{2}) f(\xi) g(\eta) \sqrt{1+4\xi^{2}} \sqrt{1+4\eta^{2}} d\xi d\eta $$

\noindent
where $\xi=(\xi_{1},\xi_{2}), \eta=(\eta_{1},\eta_{2})$. A direct use of Schwartz inequality gives us:

$$|(f \delta_{\tau=\xi^{2}} * g \delta_{\tau=\xi^{2}} * h)(z_{1},z_{2},z_{3})| \leq ||f||_{L^{2}(P)}||g||_{L^{2}(P)} \cdot$$
 $$\left( \int |h|^{2}(z_{1}-\xi_{1}-\eta_{1},z_{2}-\xi_{2}-\eta_{2}, z_{3} - \xi^{2}-\eta^{2}) \sqrt{1+4\xi^{2}} \sqrt{1+4\eta^{2}} d\xi d\eta\right)^{\q} \leq $$

$$ 2^{\frac{i+j}{2}} ||f||_{L^{2}(P)}||g||_{L^{2}(P)}  \left( \int |h|^{2}(z_{1}-\xi_{1}-\eta_{1},z_{2}-\xi_{2}-\eta_{2}, z_{3} - \xi^{2}-\eta^{2}) d\xi d\eta\right)^{\q}$$

We use the change of variables $(\xi_{1},\eta_{1},\eta_{2}) \rightarrow (\z_{1},\z_{2},\z_{3})$ given by the system:

\begin{eqnarray} \label{sys3}
  \left\{
        \begin{array}{rl}
		 z_{1}-\xi_{1}-\eta_{1}=\z_{1}  \\
		 z_{2}-\xi_{2}-\eta_{2}=\z_{2}  \\
		 z_{3}-\xi^{2}-\eta^{2}=\z_{3}
	\end{array}\right.
\end{eqnarray}

This Jacobian of this transformation is $\q(\eta_{1}-\xi_{1})^{-1}$. If we were to integrate over a region where $|\eta_{1}| \geq |\eta_{2}|$, then we would get $|\eta_{1} -\xi_{1}| \geq 2^{j-2}$ (here it is important that $i \leq j-3$) which would lead us to:

$$|(f \delta_{\tau=\xi^{2}} * g \delta_{\tau=\xi^{2}} * h)(z_{1},z_{2},z_{3})| \leq 2^{\frac{i}{2}} ||f||_{L^{2}(P)}||g||_{L^{2}(P)}  \left( \int ||h||_{L^{2}}^{2} d\xi_{1} \right)^{\q} \leq$$

$$2^{i} ||f||_{L^{2}(P)}||g||_{L^{2}(P)} ||h||_{L^{2}}$$

\noindent
the last inequality being justified by the fact that we integrate over a region where $|\xi| \approx 2^{i}$. 

The way to fix the proof is to split $g = g_{1}+g_{2}$ where $g_{1}$ is localized in a region where $|\eta_{1}| \geq |\eta_{2}|$ and $g_{2}$ is localized in a region where $|\eta_{2}| \geq |\eta_{1}|$. For $g_{1}$ we apply the above argument, while for $g_{2}$ we use the change of variables $(\xi_{2},\eta_{1},\eta_{2}) \rightarrow (\z_{1},\z_{2},\z_{3})$ given by the same system (\ref{sys3}). The Jacobian is given in this case by $\q(\eta_{2}-\xi_{2})^{-1}$ and we can argue in the same way as before.

By adding the results we obtain for $g_{1}$ and $g_{2}$, we get (\ref{eq1}). 

The estimates for $f \delta_{\bar{P}} * g \delta_{P}$ and $f \delta_{\bar{P}} * g \delta_{\bar{P}}$ can be obtained in a similar way. 
 
\begin{bfseries} Case 2: {\mathversion{bold}  $i \cong j$} \end{bfseries}

In this case we make use of the the Strichartz estimate:

$$|| \int a(\xi)  e^{i(x \cdot \xi + t \cdot \xi^{2})} d\xi  ||_{L^{4}} \leq C  ||a||_{L^{2}_{\xi}} $$

In our case, $f$ and $g$ are localized at $\approx 2^{j}$, therefore the above inequality gives us $||\mathcal{F}^{-1}(u \delta_{P})||_{L^{4}} \leq 2^{\frac{j}{2}} ||u||_{L^{2}(P)}$ and the similar one for $v$. 

$$||f \delta_{P}* g \delta_{P}||_{L^{2}} = ||\mathcal{F}^{-1}(f \delta_{P}) \cdot \mathcal{F}^{-1}(g \delta_{P})||_{L^{2}} \leq$$

$$||\mathcal{F}^{-1}(f \delta_{P})||_{L^{4}}||\mathcal{F}^{-1}(g \delta_{P})||_{L^{4}} \leq 2^{j} ||f||_{L^{2}(P)} ||g||_{L^{2}(P)}$$

The Strichartz estimate is valid also for $\bar{P}$, i.e. :

$$|| \int a(\xi,-\xi^{2})  e^{i(x \cdot \xi - t \cdot \xi^{2})} d\xi  ||_{L^{4}} \leq C  ||a||_{L^{2}_{\xi}(\bar{P})} $$

\noindent
therefore we get also the estimates involving $\delta_{\bar{P}}$. 

\end{proof}

We are interested in a general result for $f \delta_{P^{1}} * g \delta_{P^{2}}$ restricted to a $P_{c}$, where $P^{1},P^{2} \in \{P_{c}, \bar{P}_{c}: c \in \R \}$ . The result of Proposition 1 is true with a simple modification on the condition imposed on $c$ : $|c \pm c_{1} \pm c_{2}| \leq 2^{i+j-2}$. If we work out the details of the proof, we can see the way the $\pm$ in this condition are related to the $\pm$ in $P^{1}= \{\tau \pm \xi^{2}=c_{1} \}$ and $P^{2}= \{ \tau \pm \xi^{2}=c_{2} \}$. On the other hand, when we will apply this result we will have the condition $|c \pm c_{1} \pm c_{2}| \leq 2^{i+j-2}$ fulfilled for all choices of signs. Therefore we will not be concerned about the sign connection. 

Another observation is that we will not apply this result for the extreme values for $c$'s. More exactly, if we deal with $f \delta_{P_{c_{1}}}$ and $f$ is localized at $2^{i}$ we always take $|c_{1}| < 2^{2i-2}$. This guarantees that in the support of  $f \delta_{P_{c_{1}}}$ we have $|\xi| \approx 2^{i}$ and not only that $|(\xi,\tau)| \approx 2^{i}$. When we dealt with $\delta_{P}$ we have $\tau=\xi^{2}$ therefore we had for free that $|(\xi,\tau)| \approx \sqrt{|\tau|} \approx |\xi|$.

\begin{p2} \label{pp2}

Let $f \in \mathcal{R}L^{2}(P_{c_{1}}) \ \ \mbox{or} \ \ \mathcal{R}L^{2}(\bar{P}_{c_{1}})$ and $g \in \mathcal{R}L^{2}(P_{c_{2}})$ or $\mathcal{R}L^{2}(\bar{P}_{c_{2}})$ such that $f$ is localized at $2^{i}$, $g$ at $2^{j}$ and $\min{(i,j)} \geq 1$, $|c_{1}| < 2^{2i-2}$, $|c_{2}| < 2^{2j-2}$. Then all the results listed in Proposition 1 hold true with the obvious adjustments. The adjustment for $c$ is  $|c \pm c_{1} \pm c_{2}| \leq 2^{i+j-2}$. 

\end{p2}

\vspace{.2in}

\subsection{Bilinear estimates on dyadic regions}
\noindent
\vspace{.1in}

We first derive the estimates in $X^{0,\q}$. The advantages are that $X^{0,\q}$ and $\bar{X}^{0,-\q}$ are dual to each other and we do not carry $s$ in all the computations. For a bilinear estimate we use the notation:

$$ \mathcal{X} \cdot \mathcal{Y} \rightarrow \mathcal{Z}$$

\noindent
which means that we seek for an estimate $||B(u,v)||_{\mathcal{Z}} \leq C ||u||_{\mathcal{X}} \cdot ||v||_{\mathcal{Y}}$. Here the constant $C$ may depend on some variables, like the frequency where the functions are localized.

Our function spaces involve rotations, therefore we need estimates of type:

$$\mathcal{R} \mathcal{X} \cdot \mathcal{R} \mathcal{Y} \rightarrow \mathcal{R} \mathcal{Z}$$

$\mathcal{R}$ is a first-order differential operator, i.e. it satisfies

$$\mathcal{R}(u \cdot v) = \mathcal{R}u \cdot v + u \cdot \mathcal{R}v$$

Therefore we have the implication:

\beq \label{t2}
\mathcal{R} \mathcal{X} \cdot \mathcal{Y} \rightarrow  \mathcal{Z} \ \ \mbox{and} \ \  \mathcal{X} \cdot \mathcal{R} \mathcal{Y} \rightarrow  \mathcal{Z} \Rightarrow \mathcal{R} \mathcal{X} \cdot \mathcal{R} \mathcal{Y} \rightarrow \mathcal{R} \mathcal{Z}
\eeq

\noindent

We also have:

\beq \label{t3}
\mathcal{X} \cdot \mathcal{Y} \rightarrow  \mathcal{Z} \Rightarrow \mathcal{R} \mathcal{X} \cdot \mathcal{R} \mathcal{Y} \rightarrow \mathcal{R} \mathcal{Z}
\eeq

\noindent
with the same constant in the estimate. We already saw in the Proposition (1) that there are regimes where we do have estimates without the use of rotations. These type of estimates are good for duality purposes. We do not want to involve in duality factors of type $\mathcal{R} X^{0,\q}$. 

A standard way of writing down each case looks like:

\vspace{.1in}

{\mathversion{bold} $ X_{i, d_{1}}^{0,\q} \cdot X_{j, d_{2}}^{0,\q} \rightarrow  X_{j, d_{3}}^{0,-\q}, \ \ \ \ i \leq j-4$  }

\vspace{.1in}

This means that for $u \in X_{i, d_{1}}^{0,\q}$ and $v \in X_{j, d_{2}}^{0,\q}$ we estimate the part of $B(u,v)$ (or $\tilde{B}(u,v)$) whose Fourier transform is supported in $A_{j, d_{3}}$. Formally we estimate $\mathcal{F}^{-1}(\chi_{A_{j,d_{3}}}\mathcal{F}(B(u,v)))$. This is going to be the only kind of ``abuse'' in notation which we make throughout the paper, i.e. considering $||B(u,v)||_{X_{j,d_{3}}^{s,\q}}$ even if $\mathcal{F}(B(u,v))$ is not supported in $A_{j,d_{3}}$. We choose to do this so that we do not have to relocalize every time in $A_{j,d_{3}}$. 

 Sometimes we prove estimates via duality: 

$$X \cdot Y \rightarrow Z \Longleftrightarrow X \cdot (Z)^{*} \rightarrow (Y)^{*}$$

Another simple property is the following:

$$X \cdot Y \rightarrow Z \Longleftrightarrow \bar{X} \cdot \bar{Y} \rightarrow \bar{Z}$$

\begin{p2} \label{u1}

 For the operator $B$ we have the following dyadic estimates:

\beq \label{aa1}
||B(u,v)||_{\mathcal{R}X^{0,-\q}_{k,d_{3}}} \leq 2^{\min{(i,j)}} ||u||_{\mathcal{R}X^{0,\q}_{i,d_{1}}} ||v||_{\mathcal{R}X^{0,\q}_{j,d_{2}}}
\eeq

\beq \label{aa2}
||B(\bar{u},v)||_{\mathcal{R}X^{0,-\q}_{k,d_{3}}} \leq 2^{\min{(i,j)}}  ||u||_{\mathcal{R} X^{0,\q}_{i,d_{1}}} ||v||_{\mathcal{R}X^{0,\q}_{j,d_{2}}}
\eeq

\noindent
where the parameters involved are restricted as follows: 

- $\min{(i,j)} \geq 1$

- $i \leq j-5 \Rightarrow d_{1} \leq 2^{i-2}$

- $j \leq i-5 \Rightarrow d_{2} \leq 2^{j-2}$

- $i=j,j \pm1$ and $k \leq j-5 \Rightarrow d_{3} \leq 2^{k-2}$

\end{p2}

\begin{proof}

We should make some commentaries about the statement above. First we make the choice $i \leq j$. If $i \leq j-2$, then the result is localized at frequency $\approx 2^{j}$. There is something to estimate only if $k=j, j \pm1$. 

It is only when $i=j-1,j$ that we have parts of the result at lower frequencies and then we have to provide estimates for all $k \leq j+1$.

We first deal with the case when we measure the outcome at the high frequency and then with the case when we have to measure the outcome at lower frequencies (here $i=j-1,j$). 

If we localize in a region where $|\xi| \approx 2^{k}$, the parabolas $P_{c}$ make an angle of $\approx 2^{-k}$ with the $\tau$ axis, so we have the following relation between measures:

$$d\xi d\tau \approx 2^{-k} dP_{c} dc$$

If $d \leq 2^{i-3}$ then in $A_{i,d}$ we have $|\xi| \approx 2^{i}$. Therefore for $l \leq i-3$:

\beq \label{e100}
||u||^{2}_{X^{0,\pm \q}_{i,2^{l}}} \approx 2^{\pm(l+i)} \int_{b=2^{l}}^{2^{l+1}} ||\hat{u}||^{2}_{L^{2}(P_{b2^{i}})} db 
\eeq

At this time we are ready to start the estimates.

\vspace{.1 in}

{\mathversion{bold} $ \mathcal{R}X_{i, d_{1}}^{0,\q} \cdot  \mathcal{R}X_{j, d_{2}}^{0,\q} \rightarrow  \mathcal{R}X_{j, d_{3}}^{0,-\q}$}

\vspace{.1 in}

\begin{bfseries} Case 1: $d_{1} \leq 2^{i-3}$ \end{bfseries}

\begin{bfseries} Subcase 1.1: $d_{2}, d_{3} \leq 2^{i-3}$ \end{bfseries}

In this case we get $|2^{j}d_{3}-2^{j}d_{2} -2^{i}d_{1}| \leq 2^{i+j-2}$. Therefore, if $c \in [2^{j-1}d_{3},2^{j+1}d_{3}]$ we can apply the result of Proposition \ref{pp2} to evaluate

$$||\hat{u} * \hat{v}||_{L^{2}(P_{c})} \leq \int_{I_{1}} \int_{I_{2}} ||\hat{u} \delta_{P_{b_{1}2^{i}}} * \hat{v} \delta_{P_{b_{2}2^{j}}} ||_{L^{2}(P_{c})} db_{1} db_{2} \leq $$

$$ \int_{I_{1}} \int_{I_{2}} 2^{\frac{i}{2}} ||\hat{u}||_{\mathcal{R}L^{2}(P_{b_{1}2^{i}})} ||\hat{v}||_{L^{2}(P_{b_{2}2^{j}})} db_{1} db_{2} \leq $$

$$2^{\frac{i}{2}} \left( \int_{I_{1}} (1+b_{1}2^{i})^{-1} db_{1} \right)^{\q} ||u||_{\mathcal{R}X^{0,\q}_{i,d_{1}}} \left( \int_{I_{2}} (1+b_{2}2^{j})^{-1} db_{2} \right)^{\q} ||v||_{X^{0,\q}_{j,d_{2}}} \approx $$

$$2^{-\frac{j}{2}} ||u||_{\mathcal{R}X^{0,\q}_{i,d_{1}}} ||v||_{X^{0,\q}_{j,d_{2}}}$$

Here we used the fact that $I_{1} \approx [\frac{d_{1}}{2},2d_{1}]$ which gives us $\int (1+b_{1}2^{i})^{-1} db_{1} \approx 2^{-i} $. Same thing for the integral with respect to $b_{2}$. (\ref{e100}) gives us:

$$||\tilde{B}(u,v)||^{2}_{X^{0,-\q}_{j,d_{3}}} \approx (2^{j}d_{3})^{-1} \int_{c=2^{j-1}d_{3}}^{2^{j+1}d_{3}}  ||\hat{u} * \hat{v}||^{2}_{L^{2}(P_{c})} 2^{-j} dc \leq 2^{-2j} ||u||^{2}_{\mathcal{R}X^{0,\q}_{i,d_{1}}} ||v||^{2}_{X^{0,\q}_{j,d_{2}}}$$ 

Notice that the Proposition 1 allows us to move the rotations on $v$ in all the computations above. Therefore we get also:

$$||\tilde{B}(u,v)||_{X^{0,-\q}_{j,d_{3}}} \leq 2^{-j} ||u||_{X^{0,\q}_{i,d_{1}}} ||v||_{\mathcal{R}X^{0,\q}_{j,d_{2}}}$$ 

 The estimates for $B$ are obtained by using the principle in (\ref{bb}). 

\begin{bfseries} Subcase 1.2: $d_{3} \geq 2^{i-2}$ and $d_{2} \leq 2^{j-2}$ \end{bfseries}

Making use to the global $L^{2}$ estimate for convolutions, see (\ref{ge}), gives us:

$$||\hat{u} * \hat{v}||_{L^{2}} \leq \int_{I_{1}} \int_{I_{2}} ||\hat{u} \delta_{P_{b_{1}2^{i}}} * \hat{v} \delta_{P_{b_{2}2^{j}}} ||_{L^{2}} db_{1} db_{2}$$

$$ \int_{I_{1}} \int_{I_{2}} 2^{i}||\hat{u}||_{L^{2}(P_{b_{1}2^{i}})} ||\hat{v}||_{L^{2}(P_{b_{2}2^{j}})} db_{1} db_{2} \leq 2^{\frac{i-j}{2}}  ||u||_{X^{0,\q}_{i,d_{1}}} ||v||_{X^{0,\q}_{j,d_{2}}}$$

Next

$$||\tilde{B}(u,v)||_{X^{0,-\q}_{j,d_{3}}} \approx  (2^{j} d_{3})^{-\q} ||\hat{u} * \hat{v}||_{L^{2}(A_{j,d_{3}})} \leq 2^{-j} ||u||_{X^{0,\q}_{i,d_{1}}} ||v||_{X^{0,\q}_{j,d_{2}}}$$

\noindent
where we use the fact that $d_{3} \geq 2^{i-2}$. 

\begin{bfseries} Subcase 1.3: $d_{2} \geq 2^{i-2}$ and $d_{3} \leq 2^{j-2}$ \end{bfseries}

This estimate for this case can be deduced by duality from the estimate:

$$ X_{i, d_{1}}^{0,\q} \cdot  \bar{X}_{j, d_{3}}^{0,\q} \rightarrow  \bar{X}_{j, d_{2}}^{0,-\q} \Leftrightarrow  \bar{X}_{i, d_{1}}^{0,\q} \cdot  X_{j, d_{3}}^{0,\q} \rightarrow  X_{j, d_{2}}^{0,-\q} $$

\noindent
and it is important that we get this estimate without rotations since otherwise we would not be able to use duality. The proof of the last estimate is treated in {\mathversion{bold} $ \mathcal{R}\bar{X}_{i, d_{1}}^{0,\q} \cdot  \mathcal{R}X_{j, d_{2}}^{0,\q} \rightarrow  \mathcal{R}X_{j, d_{3}}^{0,-\q}$}, Subcase 1.2.

\begin{bfseries} Subcase 1.4: $d_{2},d_{3} \geq 2^{j-2}$  \end{bfseries}

In this case we use a much simpler argument. For reference we call this the $L^{1} * L^{2} \rightarrow L^{2}$ argument. It goes as follows:

$$||\hat{u}||_{L^{1}} \leq 2^{\frac{3}{2}i} d_{1}^{\q} ||\hat{u}||_{L^{2}} \leq 2^{i} ||u||_{X_{i,d_{1}}^{0,\q}} $$

 Then we continue with:

$$||\tilde{B}(u,v)||_{X^{0,-\q}_{j,d_{3}}} \approx  2^{-j} ||\hat{u} * \hat{v} ||_{L^{2}(A_{j,d_{3}})} \leq $$

$$ 2^{-j} ||\hat{u}||_{L^{1}} \cdot ||\hat{v}||_{L^{2}} \leq 2^{i-2j} ||u||_{X^{0,\q}_{j,d_{1}}} ||v||_{X^{0,\q}_{j,d_{2}}}$$

\begin{bfseries} Case 2: $d_{1} \geq 2^{i-2}$  \end{bfseries}

 We have to deal only with the case $i \geq j-5$ is obtained by duality:

$$ \bar{X}_{j, d_{3}}^{0,\q} \cdot  X_{j, d_{2}}^{0,\q} \rightarrow  \bar{X}_{i, d_{1}}^{0,-\q} \Leftrightarrow  X_{j, d_{3}}^{0,\q} \cdot  \bar{X}_{j, d_{2}}^{0,\q} \rightarrow  X_{i, d_{1}}^{0,-\q}$$

The last estimate is treated in the next group of estimates. It is important that $d_{1} \geq 2^{i-2}$, since it is one of the cases when rotations are not needed.

\vspace{.1 in}

{\mathversion{bold} $ \mathcal{R}\bar{X}_{i, d_{1}}^{0,\q} \cdot  \mathcal{R}X_{j, d_{2}}^{0,\q} \rightarrow  \mathcal{R}X_{j, d_{3}}^{0,-\q}$}

\vspace{.1 in}

\begin{bfseries} Case 1: $d_{1} \leq 2^{i-3}$ \end{bfseries}

\begin{bfseries} Subcase 1.1: $d_{2}, d_{3} \leq 2^{i-3}$ \end{bfseries}

This case is totally similar to Subcase 1.1 in the first estimate because we have all the necessary ingredients. 

\begin{bfseries} Subcase 1.2: $d_{3} \geq 2^{i-2}$ and $d_{2} \leq 2^{j-2}$ \end{bfseries}

Same situation, this case is similar to Subcase 1.2 in the first estimate. Notice that the result is obtained without the use of rotations. 

\begin{bfseries} Subcase 1.3: $d_{2} \geq 2^{i-2}$ and $d_{3} \leq 2^{j-2}$ \end{bfseries}

This estimate for this case can be deduced by duality from the estimate:

$$ \bar{X}_{i, d_{1}}^{0,\q} \cdot  \bar{X}_{j, d_{3}}^{0,\q} \rightarrow  \bar{X}_{j, d_{2}}^{0,-\q} \Leftrightarrow  X_{i, d_{1}}^{0,\q} \cdot  X_{j, d_{3}}^{0,\q} \rightarrow  X_{j, d_{2}}^{0,-\q}$$

The last estimate was proved in {\mathversion{bold} $ \mathcal{R}X_{i, d_{1}}^{0,\q} \cdot  \mathcal{R}X_{j, d_{2}}^{0,\q} \rightarrow  \mathcal{R}X_{j, d_{3}}^{0,-\q}$}, Subcase 1.2.

\begin{bfseries} Subcase 1.4: $d_{2},d_{3} \geq 2^{j-2}$  \end{bfseries}

Use the $L^{1} * L^{2} \rightarrow L^{2}$ argument.

\begin{bfseries} Case 2: $d_{1} \geq 2^{i-2}$ \end{bfseries}

Notice that we work in the hypothesis $i \geq j-5$. Then use duality to get the estimate from $ \bar{X}_{j, d_{3}}^{0,\q} \cdot  \bar{X}_{j, d_{2}}^{0,\q} \rightarrow  X_{i, d_{1}}^{0,-\q}$. This estimate can be easily treated as if it were $\nabla \bar{X}_{j, d_{3}}^{0,\q} \cdot \nabla \bar{X}_{j, d_{2}}^{0,\q} \rightarrow  \bar{X}_{i, d_{1}}^{0,-\q}$, since $d_{1} \geq 2^{i-2}$. The conjugate of this estimate has been treated before.  

\vspace{.1 in}

{\mathversion{bold} $ \mathcal{R}X_{i, d_{1}}^{0,\q} \cdot  \mathcal{R}\bar{X}_{j, d_{2}}^{0,\q} \rightarrow  \mathcal{R} X_{j, d_{3}}^{0,-\q}$}

\vspace{.1 in}

This estimate can be proved going through the same steps as for the estimate $ \mathcal{R}X_{i, d_{1}}^{0,\q} \cdot  \mathcal{R} X_{j, d_{2}}^{0,\q} \rightarrow  \mathcal{R} X_{j, d_{3}}^{0,-\q}$. The underlying idea is that Proposition \ref{bep1} provides the estimates for $f \delta_{P} * g \delta_{\bar{P}}$ too and this is what we need here.

\vspace{.1 in}

\begin{bfseries} High - High interactions with output at low frequencies  \end{bfseries} 

\vspace{.1 in}

{\mathversion{bold} $ \mathcal{R} X_{j, d_{1}}^{0,\q} \cdot  \mathcal{R}X_{j, d_{2}}^{0,\q} \rightarrow  \mathcal{R} X_{i, d_{3}}^{0,-\q}$}

\vspace{.1 in}

Conjugation and duality give us:

$$ X_{i, d_{3}}^{0,\q} \cdot  \mathcal{R}\bar{X}_{j, d_{2}}^{0,\q} \rightarrow  X_{j, d_{1}}^{0,-\q} \Rightarrow \bar{X}_{i, d_{3}}^{0,\q} \cdot  \mathcal{R}X_{j, d_{2}}^{0,\q} \rightarrow  \bar{X}_{j, d_{1}}^{0,-\q}  \Rightarrow  X_{j, d_{1}}^{0,\q} \cdot  \mathcal{R}X_{j, d_{2}}^{0,\q} \rightarrow  X_{i, d_{3}}^{0,-\q}$$

$$ \mathcal{R} \bar{X}_{j, d_{1}}^{0,\q} \cdot  X_{i, d_{3}}^{0,\q} \rightarrow  X_{j, d_{2}}^{0,-\q} \Rightarrow \mathcal{R} X_{j, d_{1}}^{0,\q} \cdot  \bar{X}_{i, d_{3}}^{0,\q} \rightarrow  \bar{X}_{j, d_{2}}^{0,-\q} \Rightarrow \mathcal{R}X_{j, d_{1}}^{0,\q} \cdot  X_{j, d_{2}}^{0,\q} \rightarrow  X_{i, d_{3}}^{0,-\q}$$

\noindent
and this is enough to justify the estimate. 

With one exception though: $i \geq j-5$ and $d_{3} \geq 2^{i-2}$. This exception is treated in the next two cases.

\begin{bfseries} Case 1: $d_{1}, d_{2} \leq 2^{j-2}$ \end{bfseries}

The argument is similar to Subcase 1.2 in the previous estimates. 

Making use of (\ref{ge}) we obtain:

$$||\tilde{B}(u,v)||_{L^{2}} \leq \int_{I_{1}} \int_{I_{2}} ||\hat{u} \delta_{P_{b_{1}2^{i}}} * \hat{v} \delta_{P_{b_{2}2^{j}}} ||_{L^{2}} db_{1} db_{2} \leq$$

$$ \int_{I_{1}} \int_{I_{2}} 2^{j} ||\hat{u}||_{L^{2}(P_{b_{1}2^{j}})} ||\hat{v}||_{L^{2}(P_{b_{2}2^{j}})} db_{1} db_{2} \leq   ||u||_{X^{0,\q}_{j,d_{1}}} ||v||_{X^{0,\q}_{j,d_{2}}}$$

Next

$$||\tilde{B}(u,v)||_{X^{0,-\q}_{i,d_{3}}} \approx  (2^{2j})^{-\q} ||\tilde{B}(u,v)||_{L^{2}(A_{i,d_{3}})} \leq 2^{-j} ||u||_{X^{0,\q}_{j,d_{1}}} ||v||_{X^{0,\q}_{j,d_{2}}}$$

\noindent
where we use the fact that $d_{3} \geq 2^{i-2} \geq 2^{j-7}$. 

\begin{bfseries} Case 2: $\max{(d_{1}, d_{2})} \geq 2^{j-2}$ \end{bfseries}

This case is similar to Subcase 1.4 in the previous estimates and uses the trivial $L^{1} * L^{2} \rightarrow L^{2}$ argument. We skip the rest of the details. 

Notice that in both cases we did not make use of rotations. 

\vspace{.1 in}

{\mathversion{bold} $ \mathcal{R} X_{j, d_{1}}^{0,\q} \cdot  \mathcal{R}\bar{X}_{j, d_{2}}^{0,\q} \rightarrow  \mathcal{R} X_{i, d_{3}}^{0,-\q}$}

\vspace{.1 in}

In the same way as above, duality gives us the estimates as claimed in the Theorem, except for the case: $i \geq j-5$ and $d_{3} \geq 2^{i-2}$ ! The exception is treated as in the Case 1 and Case 2 above.

\end{proof}

Proposition \ref{u1} has some restrictions on the parameters involved. Therefore we have to deal with the cases left out.

\begin{p2} \label{ss}

a) Assume $1 \leq i \leq j-5$ and $d_{1} \geq 2^{i-2}$. Then we have the following estimates on dyadic pieces:

\beq \label{s113}
||B(u,v)||_{\mathcal{R}X^{0,-\q}_{j,d_{3}}} \leq 2^{i} j^{\q}||u||_{\mathcal{R}X^{0,\q}_{i,d_{1}}} ||v||_{\mathcal{R}X^{0,\q}_{j,d_{2}}}
\eeq

The same estimate holds true if $B(u,v)$ is replaced by $B(\bar{u},v)$ or $B(u,\bar{v})$. 

Assume $i=j,j \pm 1$ and $d_{3} \geq 2^{k-2}$. Then we have the estimates: 

\beq \label{s1131}
||B(u,v)||_{\mathcal{R}X^{0,-\q}_{k,d_{3}}} \leq  2^{j} j^{\q} ||u||_{\mathcal{R}X^{0,\q}_{i,d_{1}}} ||v||_{\mathcal{R}X^{0,\q}_{j,d_{2}}}
\eeq 

Both estimates hold true if $B(u,v)$ is replaced by $B(\bar{u},v)$ or $B(u,\bar{v})$.

b) For any $j, d_{2}, d_{3}$ we have:

 \beq \label{s313}
||B(u,v)||_{\mathcal{R}X^{0,-\q}_{j,d_{3}}} \leq j^{\q}||u||_{\mathcal{R}X^{0,\q}_{0}} ||v||_{\mathcal{R}X^{0,\q}_{j,d_{2}}}
\eeq

The same estimate holds true if $B(u,v)$ is replaced by $B(\bar{u},v)$ or $B(u,\bar{v})$. 

\end{p2}

We need to prepare the geometrical setup to approach Proposition 4.

\vspace{.1in}

\begin{bfseries} Preparation \end{bfseries}

In what follows we work for a while with functions of two variables, the idea being that we work on sections with $\tau=constant$.

We fix $i \leq j-5$ and $i \leq k \leq j$. Recall the definition of $\Xi^{k}$ from (\ref{Xi}) and the related entities. For each positive integer $n$ we define:

\beq \label{s120}
\Xi_{n}^{k}= \{\xi=(n 2^{-k},\t): \t=\frac{\pi}{2} \cdot \frac{l}{n}; l \in \Z \}
\eeq

For every $\xi \in \Xi^{k}$ we define:

$$f^{k}_{\xi}= \phi^{k}_{\xi} \cdot f$$

\noindent
and for fixed $n$ we define:

$$f^{k}_{n}= \sum_{\xi \in \Xi^{k}_{n}} f_{\xi}^{k}$$

Because $k$ is fixed we drop the upper index $k$ from $f^{k}_{\xi}$ and $f^{k}_{n}$ and write only $f_{\xi}$ and $f_{n}$.

We fix $n_{2} \approx 2^{k+j}$ and $n_{3} \approx 2^{k+j}$. We want to estimate the term:

$$\phi^{k}_{n_{3}} (\sum_{n_{1} \leq 2^{2i+k-2}} g_{n_{1}} * f_{n_{2}})$$

The result we expect to obtain is the following:

\begin{p2} \label{k100}
 Assume that $f,g \in \mathcal{R}L^{2}(\R^{2})$. Then we have the estimate:

\beq \label{s111}
||\phi^{k}_{n_{3}} (\sum_{n_{1} \leq 2^{i+k-2}} g_{n_{1}} * f_{n_{2}})||_{L^{2}} \leq 2^{-k} (i+k)^{\q}||\sum_{n_{1}} g_{n_{1}}||_{\mathcal{R}L^{2}} ||f_{n_{2}}||_{L^{2}}
\eeq

\noindent
and also the estimate with $\mathcal{R}$ moved on the second factor on the right.

\end{p2}

\vspace{.1in}

\begin{bfseries} Geometry of interactions \end{bfseries}

In order to approach the above stated problem we decompose:

 $$\phi^{k}_{n_{3}} (\sum_{n_{1} \leq 2^{i+k-2}} g_{n_{1}} * f_{n_{2}}) = \sum_{n_{1} \leq 2^{i+k-2}} \sum_{\xi \in \Xi_{n_{1}}^{k}} \sum_{\eta \in \Xi_{n_{2}}^{k}} \phi_{n_{3}}^{k} ( g_{\xi} * f_{\eta})$$
   
The first thing we have to understand is under what conditions the quantity $\phi^{k}_{n_{3}} ( g_{\xi} * f_{\eta})$ is nontrivial. In other words, under what condition between $\xi \in \Xi_{n_{1}}^{k}$ and $\eta \in \Xi_{n_{2}}^{k}$ we have that the support of $\phi^{k}_{\xi} * \phi^{k}_{\eta}$ intersects the support of $\phi_{n_{3}}^{k}$. 

\begin{l1}
If $\a$ is the angle between $\xi \in \Xi^{k}_{n_{1}}$ and $\eta \in \Xi^{k}_{n_{2}}$ then $\phi_{n_{3}}^{k} ( g_{\xi} * f_{\eta})$ is nontrivial iff 

\beq \label{s1}
2|\xi||\eta| \cos{\a} + \xi^{2} \in I(n_{2},n_{3},k)
\eeq

\noindent
where $I(n_{2},n_{3},k)$ is an interval of length $ 2^{-2k+2} n_{3}\approx 2^{j-k}$.

\end{l1}

\begin{proof}
$\phi^{k}_{\xi}$ is supported in a cube centered at $\xi$ and of sizes $2^{-k} \times 2^{-k}$, so each $\xi'$ in its support can be written as $\xi'=\xi+\e_{1}$ with $|\e_{1}| \leq 2^{-k}$. In a similar way   $\phi^{k}_{\eta}$ is supported in a tube centered at $\eta$ and of sizes $2^{-k} \times 2^{-k}$, so each $\eta'$ in its support can be written as $\eta'=\eta+\e_{2}$ with $|\e_{2}| \leq 2^{-k}$. 

Therefore $\phi^{k}_{\xi} * \phi^{k}_{\eta}$ is supported in a cube centered at $\xi+\eta$ and each $\xi'$ in its support can be written as $\xi'=\xi +\e'$ with $|\e'| \leq 2^{-k+1}$. 

We want this cube to be contained in the support of $\phi^{k}_{n_{3}}$.  The condition we have to impose on the center is that $|\xi+\eta| \in 2^{-k}[n_{3}-1,n_{3}+1]$ which is equivalent to $\eta^{2} + 2 \xi \cdot \eta + \xi^{2} \in  2^{-2k}[n_{3}^{2} - 2n_{3}+1,n_{3}^{2} + 2n_{3}+1]$. Our interval $I(n_{2},n_{3},k)$ in (\ref{s1}) is $[2^{-2k}(n_{3}^{2} - 2n_{3}+1)-n^{2}_{2},2^{-2k}(n_{3}^{2} + 2n_{3}+1)-n^{2}_{2}]$.

The condition is also sufficient. 

\end{proof}

We want to solve (\ref{s1}). Basically this is an inclusion equation in $|\xi|=n_{1}$ and $\a$ and we do expect the solution to come as intervals. If the solution interval contains only values of $\a$ away from zero then we have an easy characterization:

\begin{l1}

If the interval solution $I_{n_{1}}$ of (\ref{s1}) contains angles with $|\a| \geq \frac{\pi}{4} > 0$ then $|I_{n_{1}}| \approx|n_{1}|^{-1} $. 

\end{l1}

\begin{proof}

We have $2|\xi| |\eta| \cos{\a} + |\xi|^{2} \in I(n_{2},n_{3},k)$ therefore $\D (2 |\xi| |\eta| \cos{\a} + |\xi|^{2}) \approx 2^{j-k}$. The only variable here is $\a$ hence we get $|\xi|^{-1}2^{-k} \approx \D \cos{\a} \approx \D \a \sin{\a}$ and since $|\sin{\a}| \geq |\sin{\frac{\pi}{4}}| > \q$ we get $\D \a \approx |\xi|^{-1} 2^{-k} =n_{1}$.  

Note. Here $\Delta$ should be understood as a measure of the variation. It is different than the Laplacian; hopefully this will not create any confusion, as it is only in this section that $\Delta$ is seen in this way.

\end{proof}

This piece of information will suffice for these cases.

The problem becomes more complicated when we deal with solutions of (\ref{s1}) giving us values of $\a$ close to $0$. Let us pick $n_{1}^{*}$ such that the solution interval which contains $0$. 

\begin{l1}
The length of $I_{n_{1}^{*}}$ is $\approx \langle n_{1}^{*} \rangle^{-\q}$.

\end{l1}

\begin{proof}
If $\a$ is another angle in the solution interval, then $|\cos{\a} - 1| \leq |\xi|^{-1}2^{-k}$. This implies $2 \sin^{2}{\frac{\a}{2}} \leq |\xi|^{-1} 2^{-k}=(n_{1}^{*})^{-1}$ which gives us $|\a| \leq (n_{1}^{*})^{-\q}$. So the interval has length $\approx (n_{1}^{*})^{-\q}$. The case $n_{1}^{*}$ is trivial.
\end{proof}

The next question one should ask is for what values of $n_{1}$ do we still get that the solution interval contains angles less than $\frac{\pi}{4}$?

\begin{l1}

$I_{n_{1}}$ contains angles less than $\frac{\pi}{4}$ only if

\beq \label{s15}
 n_{1}^{*} -4 \leq n_{1} \leq 4 n_{1}^{*}
\eeq

\end{l1}

\begin{proof}

If $n_{1} < n_{1}^{*}-4$ and $\xi \in \Xi^{k}_{n_{1}}$, $\xi^{*} \in \Xi^{k}_{n^{*}_{1}}$then 

$$2|\xi||\eta|\cos{\a} + |\xi|^{2} \leq 2 |\xi| |\eta| + |\xi|^{2} \leq$$

$$ (|\xi|-|\xi^{*}|)(|\xi|+|\xi^{*}|+ 2|\eta|) + 2|\xi^{*}||\eta| + (\xi^{*})^{2} \leq$$

$$2(n_{1}-n_{1}^{*})2^{-k}|\eta| + 2|\xi^{*}||\eta| + (\xi^{*})^{2}$$

 We have $2(n_{1}-n_{1}^{*}) 2^{-k} |\eta| \leq - 2 \cdot 2^{-k+2} |\eta| \leq -  2^{-2k+2} n_{3} = -|I(n_{2},n_{3},k)|$. Therefore we do not get any solutions in this case. 

If $4n_{1}^{*} < n_{1}$ and $|\a| \leq \frac{\pi}{4}$ then 

$$2n_{1}|\eta|\cos{\a} + \xi^{2} \geq n_{1}|\eta| + n_{1}^{2} \geq$$

$$ (n_{1}-n_{1}^{*})(n_{1} + n_{1}^{*} + |\eta|) + 2n_{1}^{*}|\eta| + (n_{1}^{*})^{2} $$

In the same way as before we can prove that $(n_{1}-n_{1}^{*})(n_{1} + n_{1}^{*} + |\eta|) \geq |I(n_{2},n_{3},k)|$, therefore we do not have solutions under these conditions. 

\end{proof}

We are interested in characterizing all possible solution intervals for different values of $n_{1} \in [n_{1}^{*}-4, 4 n_{1}^{*}]$. Let fix $d$ positive integer such that $1 \leq 2^{d} \leq (n_{1}^{*})^{\q}$. 

\begin{r1} \label{r1}
The possible set of values for $d$ has cardinality  $\approx ln{((n_{1}^{*})^{\q})} \leq i+k $. 
\end{r1}

\begin{l1}

If $I_{n_{1}}$ contains $\a$'s with $|\a| \approx 2^{d} (n_{1}^{*})^{-\q}$ then $|I_{n_{1}}| \approx 2^{-d} (n_{1}^{*})^{-\q}$

\end{l1}

\begin{proof}
 We have  $|\xi|^{-1} 2^{-k} \approx \D \cos{\a} \approx \D \a \sin{\a} \approx 2^{d} (n^{*}_{1})^{-\q} \D \a$ from which we conclude that $\D \a \approx 2^{-d} (n_{1}^{*})^{-\q}$. Here we make use of the fact that $ n_{1}^{*} -4 \leq |\xi| \leq 4 n_{1}^{*}$. 
 
\end{proof}

This suggests to split $I_{d}=[2^{d-1}(n_{1}^{*})^{-\q}, 2^{d+1}(n_{1}^{*})^{-\q}] = \bigcup_{l=1}^{2^{2d}} I^{l}_{d}$ where $|I^{l}_{d}| \approx 2^{-d} (n_{1}^{*})^{-\q}$. This gives us

$$[0,\frac{\pi}{4}] = \bigcup_{d} I_{d} = \bigcup_{d} \bigcup_{l=1}^{2^{2d}} I_{d}^{l}$$

We obtain this way a map $n_{1} \rightarrow  I_{n_{1}}= I_{d}^{l} \rightarrow(d,l)$. We denote this map by $h$. $h$ is ``almost'' injective in the following sense:

\begin{l1}
  For any $(d,l)$ we have $|h^{-1}(d,l)| \leq 4$ 
\end{l1}

\begin{proof}

Suppose that we have $h(n_{1})=h(\bar{n}_{1})$ for $|n_{1}-\bar{n}_{1}| \geq 4$. This implies that there are $\xi \in \Xi_{n_{1}}^{k}$ and $\bar{\xi} \in \Xi_{\bar{n}_{1}}^{k}$ such that $2|\xi||\eta| \cos{\a} + \xi^{2} \in I(n_{2},n_{3},k)$ and $2 |\bar{\xi}| |\eta| \cos{\a} + \bar{\xi}^{2} \in I(n_{2},n_{3},k)$ for the same $\a \in I_{d}^{l}$. On the other hand we have: 

$$|2 |\xi| |\eta| \cos{\a} + \xi^{2} - 2 |\bar{\xi}| |\eta| \cos{\a} - \bar{\xi}^{2}| = |(|\xi|-|\bar{\xi}|)(2|\eta| \cos{\a} + |\xi| + |\bar{\xi}|)| \geq  $$

$$|n_{1}-n_{2}| 2^{-k} |\eta| \geq |I(n_{2},n_{3},k)|$$

This is in contradiction with the fact that both quantities are in $I(n_{2},n_{3},k)$. 

Therefore for every $(d,l)$ there are at most $4$ $n$'s such that $h(n)=(d,l)$. 

\end{proof}

\begin{proof}[Proof of Proposition \ref{k100}]

\noindent

\begin{bfseries} Case 1 \end{bfseries}. We first deal with the $n_{1}$'s for which $I_{n_{1}}$ contains angles less than $\frac{\pi}{4}$ and then with the others.

Let's assume $h(n_{1})=(d,l)$. Denote by $m_{1}=2^{-k} n_{1}$ and $m_{2}=2^{-k} n_{2}$. We prefer these substitutions because we know that $\xi \in \Xi^{k}_{n_{1}}$ implies $|\xi|=m_{1}$ and $\eta \in \Xi^{k}_{n_{2}}$ implies $|\eta|=m_{2}$. 

We split $\Xi^{k}_{n_{1}}$ and $\Xi^{k}_{n_{2}}$ in angular sectors of size $\approx 2^{-d} m^{-\q}_{1} 2^{-\frac{k}{2}}$:

$$\Xi^{k}_{n_{1}} = \cup_{l} {A_{l}} \ \ \mbox{and} \ \ \Xi^{k}_{n_{2}} = \cup_{l} {B_{l}}$$

\noindent
with the following properties:

- any two $\xi$'s in the same $A_{l}$ make an angle of at most $\approx 2^{-d} m^{-\q}_{1} 2^{-\frac{k}{2}}$

- any two $\eta$'s in the same $B_{l}$ make an angle of at most $\approx 2^{-d} m^{-\q}_{1} 2^{-\frac{k}{2}}$

- the angle between a $\xi \in A_{l}$ and an $\eta \in B_{l}$ (for the same $l$!) is $\approx 2^{d} m^{-\q}_{1} 2^{-\frac{k}{2}}$; more exactly it is an angle in $I_{d}^{l}$   

We have:

$$ \phi^{k}_{n_{3}} ( g_{n_{1}} * f_{n_{2}}) = \sum_{\xi \in \Xi_{n_{1}}^{k}} \sum_{\eta \in \Xi_{n_{2}}^{k}} \phi_{n_{3}} ( g_{\xi} * f_{\eta}) \approx \sum_{l} \sum_{\xi \in A_{l}} g_{\xi} * \sum_{\eta \in B_{l}} f_{\eta}$$

We use a simple estimate:

$$||\sum_{\xi \in A_{l}} g_{\xi} * \sum_{\eta \in B_{l}} f_{\eta}||_{L^{2}} \leq ||\sum_{\xi \in A_{l}} g_{\xi}||_{L^{1}}||\sum_{\eta \in B_{l}} f_{\eta}||_{L^{2}}$$

The support of $\sum_{\xi \in A_{l}} g_{\xi}$ has sizes $2^{-k} \times m_{1} \cdot 2^{-d} m^{-\q}_{1} 2^{-\frac{k}{2}}$:

$$||\sum_{\xi \in A_{l}} g_{\xi}||_{L^{1}} \leq ( 2^{-d} m^{\q}_{1} 2^{-\frac{k}{2}})^{\q} 2^{-\frac{k}{2}} ||\sum_{\xi \in A_{l}} g_{\xi}||_{L^{2}}$$

The supports of $\sum_{\xi \in A_{l}} g_{\xi} * \sum_{\eta \in B_{l}} f_{\eta}$ are disjoint with respect to $l$ because it has an angular localization depending on $l$, therefore:

$$||\phi^{k}_{n_{3}} ( g_{n_{1}} * f_{n_{2}})||_{L^{2}}^{2} \approx \sum_{l} ||\sum_{\xi \in A_{l}}  g_{\xi}* \sum_{\eta \in B_{l}} f_{\eta})||^{2}_{L^{2}} \leq $$

$$\sum_{l} 2^{-d} m^{\q}_{1} 2^{-\frac{k}{2}} 2^{-k} ||\sum_{\xi \in A_{l}} g_{\xi}||^{2}_{L^{2}} ||\sum_{\eta \in B_{l}} f_{\eta}||^{2}_{L^{2}} \leq$$

$$2^{-d} m^{\q}_{1} 2^{-\frac{k}{2}} 2^{-k} ||f_{n_{2}}||^{2}_{L^{2}} \sup_{l} ||\sum_{\xi \in A_{l}} g_{\xi}||^{2}_{L^{2}} $$

$$\mbox{or}$$

$$2^{-d} m^{\q}_{1} 2^{-\frac{k}{2}} 2^{-k} ||g_{n_{1}}||^{2}_{L^{2}} \sup_{l} ||\sum_{\eta \in B_{l}} f_{\eta}||^{2}_{L^{2}}$$

For the first option in the estimate above we use the rotations on $g$:

$$\sup_{l} ||\sum_{\xi \in A_{l}} g_{\xi}||^{2}_{L^{2}} \leq 2^{-d} m^{-\q}_{1} 2^{-\frac{k}{2}} ||g_{n_{1}}||^{2}_{\mathcal{R}L^{2}}$$

For the second option we use the rotations on $f$:

$$\sup_{l} ||\sum_{\eta \in B_{l}} f_{\eta}||^{2}_{L^{2}} \leq 2^{-d} m^{-\q}_{1} 2^{-\frac{k}{2}} || f_{n_{2}}||^{2}_{\mathcal{R}L^{2}}$$

The argument continues the same way regardless to whether we choose to use rotations on $f$ or $g$.  We choose to use rotations on $g$:

$$||\phi^{k}_{n_{3}} ( g_{n_{1}} * f_{n_{2}})||_{L^{2}} \leq 2^{-d} 2^{-k} ||g_{n_{1}}||_{\mathcal{R}L^{2}} ||f_{n_{2}}||_{L^{2}}$$

Next, we can add up with respect to those $n_{1}$ such that $h(n_{1})=(d,l)$ with $d$ fixed. We know there are about $2^{2d}$ such $n_{1}$'s, therefore we get:

$$||\phi^{k}_{n_{3}} ( \sum_{n_{1} \in h^{-1}(d, \cdot)} g_{n_{1}} * f_{n_{2}})||_{L^{2}} \leq 2^{-k} ||\sum_{n_{1} \in h^{-1}(d, \cdot)} g_{n_{1}}||_{\mathcal{R}L^{2}} ||f_{n_{2}}||_{L^{2}}$$

In the end we need to sum up with respect to $d$ for which we know there are about $i+k$ values (see Remark 1), so we get:

$$||\phi_{n_{3}} ( \sum_{n_{1}} g_{n_{1}} * f_{n_{2}})||_{L^{2}} \leq 2^{-j} (i+k)^{\q}||\sum_{n_{1}} g_{n_{1}}||_{\mathcal{R}L^{2}} ||f_{n_{2}}||_{L^{2}}$$

\begin{bfseries} Case 2 \end{bfseries}. We deal with the $n_{1}$'s for which $I_{n_{1}}$ contains angles greater than $\frac{\pi}{4}$. We know that the angle localization should be of order $n_{1}^{-1}$. 

$\Xi^{k}_{n_{1}}$ comes already with this localization:

$$\Xi^{k}_{n_{1}} = \cup_{\xi \in \Xi^{k}_{n_{1}}} \{\xi\}$$

We split $\Xi^{k}_{n_{2}}$ in angular sectors of size $\approx n^{-1}_{1}$:

$$\Xi^{k}_{n_{2}} = \cup_{\xi \in \Xi^{k}_{n_{1}}} {B_{\xi}}$$

\noindent
with the following properties:

- any two $\eta$'s in the same $B_{\xi}$ make an angle of at most $\approx n^{-1}_{1}$

- the angle between a $\xi$ and an $\eta \in B_{\xi}$ (for the same $\xi$!) is in the interval which comes out as a solution for (\ref{s1}).  

We have:

$$ \phi^{k}_{n_{3}} ( g_{n_{1}} * f_{n_{2}}) \approx \sum_{\xi} g_{\xi} * \sum_{\eta \in B_{\xi}} f_{\eta}$$

We use a simple estimate:

$$||g_{\xi} * \sum_{\eta \in B_{\xi}} f_{\eta}||_{L^{2}} \leq ||g_{\xi}||_{L^{1}} ||\sum_{\eta \in B_{\xi}} f_{\eta}||_{L^{2}}$$

The support of $g_{\xi}$ has sizes $ \approx 2^{-k} \times 2^{-k}$

$$||g_{\xi}||_{L^{1}} \leq 2^{-k} ||g_{\xi}||_{L^{2}}$$

The supports of $g_{\xi} * \sum_{\eta \in B_{\xi}} f_{\eta}$ are disjoint with respect to $\xi$ due to the angular localization, therefore we get:

$$||\phi_{n_{3}}^{k} ( g_{n_{1}} * f_{n_{2}})||_{L^{2}}^{2} \approx \sum_{\xi} || g_{\xi} * \sum_{\eta \in B_{\xi}} f_{\eta} ||^{2}_{L^{2}} \leq $$

$$\sum_{\xi} 2^{-2k} || g_{\xi}||^{2}_{L^{2}} ||\sum_{\eta \in B_{\xi}}  f_{\eta}||^{2}_{L^{2}} \leq$$

$$ 2^{-2k} || f_{n_{2}}||^{2}_{L^{2}} \sup_{\xi} ||g_{\xi}||^{2}_{L^{2}} \ \ \ \ \ \mbox{or} \ \ \ \  2^{-2k} || g_{n_{1}}||^{2}_{L^{2}} \sup_{\xi} ||\sum_{\eta \in B_{\xi}}  f_{\eta}||^{2}_{L^{2}}$$

As before (see Case 1) we chose to use rotations on $g$; working with rotations on $f$ is completely similar. We continue with: 

$$\sup_{\xi} || g_{\xi}||^{2}_{L^{2}} \leq  n^{-1}_{1} ||g_{n_{1}}||^{2}_{\mathcal{R}L^{2}}$$

Then we end up with the following estimate:

$$||\phi^{k}_{n_{3}} ( g_{n_{1}} * f_{n_{2}})||_{L^{2}} \leq n_{1}^{-\q} 2^{-k} ||g_{n_{1}}||_{\mathcal{R}L^{2}} ||f_{n_{2}}||_{L^{2}}$$

Summing up with respect to $n_{1}$ for which we have at most $2^{i+k}$ values we obtain:

$$||\phi_{n_{3}} ( \sum_{n_{1}} g_{n_{1}} * f_{n_{2}})||_{L^{2}} \leq 2^{-k} (i+k)^{\q} ||\sum_{n_{1}} g_{n_{1}}||_{\mathcal{R}L^{2}} ||f_{n_{2}}||_{L^{2}}$$

Adding up the estimates from Case 1 and Case 2 gives us:

$$||\phi_{n_{3}} ( \sum_{n_{1}} g_{n_{1}} * f_{n_{2}})||_{L^{2}} \leq 2^{-k} (i+k)^{\q} ||\sum_{n_{1}} g_{n_{1}}||_{\mathcal{R}L^{2}} ||f_{n_{2}}||_{L^{2}}$$

\noindent
where now we sum over all $n_{1}$'s.

\end{proof}

An immediate corollary of Proposition $\ref{k100}$ deals with the case when we want to consider the interaction
 $$\sum_{n_{3} \in J} \phi^{k}_{n_{3}} ( \sum_{n_{1}} g_{n_{1}} * \sum_{n_{2} \in I} f_{n_{2}})$$ 

\noindent
where $I$ and $J$ are possible set of values for $n_{2}$ and $n_{3}$  such that we still have the condition $|\eta| \approx 2^{j}$ for every $\eta \in \Xi^{k}_{n_{2}}$ and $|\eta| \approx 2^{j}$ for every $\eta \in \Xi^{k}_{n_{3}}$. We denote by $|I|, |J|$ the cardinal of $I, J$ respectively. 

\begin{c3} In the same conditions as in Proposition 3, we have the estimate:

\beq \label{s112}
||\sum_{n_{3} \in J} \phi^{k}_{n_{3}} ( \sum_{n_{1}} g_{n_{1}} * \sum_{n_{2} \in I} f_{n_{2}})||_{L^{2}} \leq 
\eeq

$$2^{-k} (i+k)^{\q} |I|^{\q} |J|^{\q} ||\sum_{n_{1}} g_{n_{1}}||_{\mathcal{R}L^{2}} ||\sum_{n_{2} \in I} f_{n_{2}}||_{L^{2}}$$

\noindent
and the corresponding estimate when we move $\mathcal{R}$ on the second term.

\end{c3}

\begin{proof}
Fix $n_{3} \in J$. By using Cauchy-Schwartz we get:

$$||\phi^{k}_{n_{3}} ( \sum_{n_{1}} g_{n_{1}} * \sum_{n_{2} \in I} f_{n_{2}})||_{L^{2}} \leq \sum_{n_{2} \in I} ||\phi^{k}_{n_{3}} ( \sum_{n_{1}} g_{n_{1}} *  f_{n_{2}})||_{L^{2}}$$

$$\left( \sum_{n_{2} \in I} 1 \right)^{\q} \left( \sum_{n_{2} \in I} ||\phi^{k}_{n_{3}} ( \sum_{n_{1}} g_{n_{1}} * f_{n_{2}})||^{2}_{L^{2}}\right)^{\q}$$

\noindent
and then use the result of Proposition 3 for each $n_{2} \in I$. Next:

$$||\sum_{n_{3} \in J} \phi^{k}_{n_{3}} ( \sum_{n_{1}} g_{n_{1}} * \sum_{n_{2} \in I} f_{n_{2}})||^{2}_{L^{2}} \approx \sum_{n_{3} \in J} || \phi^{k}_{n_{3}} ( \sum_{n_{1}} g_{n_{1}} * \sum_{n_{2} \in I} f_{n_{2}})||^{2}_{L^{2}} \leq$$

$$\sum_{n_{3} \in J} 2^{-2k} (i+k) |I| ||\sum_{n_{1}} g_{n_{1}}||^{2}_{\mathcal{R}L^{2}} ||\sum_{n_{2} \in I} f_{n_{2}}||^{2}_{L^{2}} = $$

$$ 2^{-2k} (i+k) |I| |J| ||\sum_{n_{1}} g_{n_{1}}||^{2}_{\mathcal{R}L^{2}} ||\sum_{n_{2} \in I} f_{n_{2}}||^{2}_{L^{2}}$$

\end{proof}

\begin{proof}[Proof of Proposition \ref{ss}]

a) In order to prove (\ref{s113}) we apply the result in (\ref{s112}) for $k=j$. 
In terms of the $\tau$ variable  we can derive from (\ref{s112}) a pointwise estimate. $\hat{u}$ is supported in $A_{j,d_{2}}$ whose section with the plane $\tau=\tau_{2}$ is an annulus of thickness $\approx d_{2}$. Similar for $\hat{v}$ and  we have:

$$\hat{u}(\xi,\tau_{2}) = \sum_{n_{2} \in I} \phi^{j}_{n_{2}}(\xi) \hat{u}(\xi,\tau_{2})  $$

$$\hat{v}(\xi,\tau_{3}) = \sum_{n_{3} \in J} \phi^{j}_{n_{3}}(\xi) \hat{u}(\xi,\tau_{3})  $$

\noindent
where $|I| \approx 2^{j}d_{2}$ and $|J| \approx 2^{j}d_{3}$. We have the pointwise estimate:

$$||\chi_{A_{j,d_{3}}} ( \hat{u}(\cdot, \tau_{1}) * \hat{v}(\cdot, \tau_{2}) ||_{L^{2}_{\tau=\tau_{1}+\tau_{2}}} \leq 2^{-j} j^{\q} (2^{2j} d_{2}d_{3})^{\q} ||\hat{u}(\cdot, \tau_{1})||_{\mathcal{R}L^{2}} || \hat{v}(\cdot, \tau_{2})||_{L^{2}}$$

Now we can derive the global estimate:

$$||\chi_{A_{j,d_{3}}} ( \hat{u}(\cdot, \tau_{1}) * \hat{v}(\cdot, \tau_{2}) )||_{L^{2}} \leq  2^{-j} j^{\q} (2^{2j} d_{2}d_{3})^{\q} ||\hat{u}||_{L^{1}_{\tau}\mathcal{R}L^{2}_{\xi}} ||\hat{v}||_{L^{2}} \leq $$

$$2^{-j} j^{\q} (2^{2j} d_{2}d_{3})^{\q} 2^{i} || \hat{u}||_{\mathcal{R}L^{2}} || \hat{v}||_{L^{2}} \approx 2^{-j} j^{\q} (2^{2j} d_{2}d_{3})^{\q} 2^{i} ||\hat{u}||_{\mathcal{R}L^{2}} ||\hat{v}||_{L^{2}}$$

\noindent
the last estimate being justified by the fact that the size of $A_{i,d_{1}}$ in the $\tau$ direction is $\approx 2^{2i}$. At the level of Bourgain spaces, the last estimate becomes:

$$ ||\tilde{B}(u,v)||_{\mathcal{R}X^{0,-\q}_{j,d_{3}}} \approx (2^{j} d_{3})^{-\q} || \hat{u} * \hat{v}||_{L^{2}(A_{j,d_{3}})} \leq$$

$$2^{-j} j^{\q} (2^{j} d_{2})^{\q} 2^{i} || \hat{u}||_{\mathcal{R}L^{2}} || \hat{v}||_{L^{2}} \approx 2^{-j} j^{\q}  ||u||_{\mathcal{R}X^{0,\q}_{i,d_{1}}} ||v||_{X^{0,\q}_{j,d_{2}}}  $$

The estimates for $\tilde{B}(\bar{u},v)$ can be obtained in a similar way. The basic idea is that in this particular setup, when interacting with $v$, $\bar{u}$ behaves like $u$ and this is due to the fact that $i \leq j-5$ and $d_{1} \geq 2^{i-2}$.   

For the estimate for $\tilde{B}(u,\bar{v})$ we need a simple observation: only if $d_{3} \approx 2^{j}$ we have a nontrivial interaction. The rest is trivial.   

The last estimate can be easily derived by duality from the first ones.

b) This is done in a similar way. The only potential difference from the previous argument would be that the condition $d_{1} \geq 2^{i-2}$ is not present here anymore. But this was used to conclude that the weight coming from $(1+|\tau-\xi^{2}|) \approx 2^{2i}$ in the support of $\hat{u}$. On the other hand we deal with the case $i=0$, hence $(1+|\tau-\xi^{2}|) \approx 1$ on the whole support of $\hat{u}$. 

In the end the case $j \leq 5$ is not covered by the previous argument since we do not have anymore the condition $i \leq j-5$ fulfilled. But these cases are essentially reduced to trivial $L^{2}$ estimates which can be easily derived. 

The estimates for $\tilde{B}(\bar{u},v)$ and $\tilde{B}(u,\bar{v})$can be obtained in a similar way.

\end{proof}

We are ready to provide the estimates on dyadic pieces only with respect to the frequency. The next result is the proof of part a) of Theorem \ref{tb1} in the particular case $s=0$.

\begin{c1}

Assume $i \leq j$. We have the following bilinear estimates:

\beq \label{s220}
||B(u,v)||_{\mathcal{R}X^{0,-\q,1}_{k}} \leq j^{\frac{3}{2}} 2^{i} ||u||_{\mathcal{R}X^{0,\q,1}_{i}} ||v||_{\mathcal{R}X^{0,\q,1}_{j}}
\eeq

The same estimate holds true if $B(u,v)$ is replaced by $B(\bar{u},v)$ or $B(u,\bar{v})$.

\end{c1}

\vspace{.1in}

\begin{proof}

We provide the argument in a particular case, namely when $k=j$. 

We fix $d_{3}$ and making use of  (\ref{aa1}) and (\ref{s113})we estimate

$$||B(u,v)||_{\mathcal{R}X^{0,-\q}_{j,d_{3}}} \leq \sum_{d_{1},d_{2}} ||B(u_{\cdot,d_{1}},v_{\cdot,d_{2}})||_{\mathcal{R}X^{0,-\q}_{j,d_{3}}} \leq $$

$$\sum_{d_{1} \leq 2^{i-3}} \sum_{d_{2}} 2^{i} ||u_{\cdot,d_{1}}||_{\mathcal{R}X^{0,\q}_{i,d_{1}}} ||v_{\cdot,d_{2}}||_{\mathcal{R}X^{0,\q}_{j,d_{2}}} + $$ 

$$ \sum_{d_{1} \geq 2^{i-2}}  \sum_{d_{2}} 2^{i} j^{\q} ||u_{\cdot,d_{1}}||_{\mathcal{R}X^{0,\q}_{i,d_{1}}} ||v_{\cdot,d_{2}}||_{\mathcal{R}X^{0,\q}_{j,d_{2}}}\leq $$

$$ 2^{i} ||u||_{\mathcal{R}X^{0,\q,1}_{i}} ||v||_{\mathcal{R}X^{0,\q,1}_{j}} + j^{\q} 2^{i} ||u||_{\mathcal{R}X^{0,\q,1}_{i}} ||v||_{\mathcal{R}X^{0,\q,1}_{j}}$$

Next we sum up with respect to $d_{3}$ and obtain:

$$||B(u,v)||_{\mathcal{R}X^{0,-\q,1}_{j}} \leq  j^{\frac{3}{2}} 2^{i} ||u||_{\mathcal{R}X^{0,\q,1}_{i}} ||v||_{\mathcal{R}X^{0,\q,1}_{j}}$$

The same type of argument gives the rest of the estimates. 

\end{proof}

The estimate (\ref{s220}) and the similar ones for  $B(\bar{u},v)$ or $B(u,\bar{v})$ are good as long as $i \approx j$ since eventually we will be able to control powers of $i$. Otherwise the $j$ factor cannot be controlled in any way, so we cannot close the bilinear estimates in $X^{s,\q,1}$. 

Therefore we are interested in dealing with the case when $5i \leq j$; $5$ was randomly chosen, any constant big enough would be good for our purposes. 

In this case, we still have some ``good'' bilinear estimates in the sense that they do not contain logarithms of the high frequency. If we deal with dyadic pieces at the high frequencies which are at some distance from $P$, then we can obtain an improvement. 

\begin{p2} 

a) Assume $i \leq j-5$, $d_{1} \geq 2^{i-2}$ and $d_{2},d_{3} \geq 2^{-i}$. Then we have the following estimates on dyadic pieces:

\beq \label{s114}
||B(u,v)||_{\mathcal{R}X^{0,-\q}_{j,d_{3}}} \leq 2^{i} i^{\q}||u||_{\mathcal{R}X^{0,\q}_{i,d_{1}}} ||v||_{\mathcal{R}X^{0,\q}_{j,d_{2}}}
\eeq

The same estimate holds true if $B(u,v)$ is replaced by $B(\bar{u},v)$ or $B(u,\bar{v})$.

b) Assume that $d_{2}, d_{3} \geq 1$. Then we have the estimates:

\beq \label{s514}
||B(u,v)||_{\mathcal{R}X^{0,-\q}_{j,d_{3}}} \leq ||u||_{\mathcal{R}X^{0,\q}_{0}} ||v||_{\mathcal{R}X^{0,\q}_{j,d_{2}}}
\eeq

The same estimate holds true if $B(u,v)$ is replaced by $B(\bar{u},v)$ or $B(u,\bar{v})$. 

c) All the estimates in a) and b) hold true without involving rotations with the additional factors:
$2^{i}$ for (\ref{s114}) and none for (\ref{s514}).

\end{p2}

\begin{proof}
 In order to prove (\ref{s114}) we apply the result in (\ref{s112}) for $k=i$. 
In terms of the $\tau$ variable  we can derive from (\ref{s112}) a pointwise estimate. $\hat{u}$ is supported in $A_{j,d_{2}}$ whose section with the plane $\tau=\tau_{2}$ is an annulus of thickness $\approx d_{2}$. Similar thing for $\hat{v}$ and  we have:

$$\hat{u}(\xi,\tau_{2}) = \sum_{n_{2} \in I} \phi^{i}_{n_{2}}(\xi) \hat{u}(\xi,\tau_{2})  $$

$$\hat{v}(\xi,\tau_{3}) = \sum_{n_{3} \in J} \phi^{i}_{n_{3}}(\xi) \hat{u}(\xi,\tau_{3})  $$

\noindent
where $|I| \approx 2^{i}d_{2}$ and $|J| \approx 2^{i}d_{3}$.

From this point on we have the same setup as in the proof of (\ref{s112}), just that $k=i$ instead of $k=j$. With this only correction, the argument there can be copied verbatim now. 

The rest of the estimates are obtained in a similar way.

\end{proof}

At this time we can prove the claim in part b) of the Theorem \ref{tb1} in the particular case $s=0$.

\begin{c1}

 Assume $5i \leq j$. We have the following bilinear estimates:

\beq \label{s130}
||B(u,v)||_{\mathcal{R}X^{0,-\q,1}_{j, \geq 2^{-i}}} \leq i^{\frac{3}{2}} 2^{i} ||u||_{\mathcal{R}X^{0,\q,1}_{i}} ||v||_{\mathcal{R}X^{0,\q,1}_{j,\geq 2^{-i}}}
\eeq

The same estimate holds true if $B(u,v)$ is replaced by $B(\bar{u},v)$ or $B(u,\bar{v})$.

\end{c1}

\begin{proof}

We prove the first estimate, the other ones being treated in a similar way. Without losing the generality we can assume that $u=u_{i}$ and $v=v_{j}$. The main observation is that the small frequency cannot change the distance to $P$ of the high frequency with a factor bigger than $2^{i+1}$. Therefore if we decompose:

$$\sum_{d_{2} \geq 2^{-i}} v_{\cdot, d_{2}} =  \sum_{d_{2}=2^{-i}}^{2^{i+1}}v_{\cdot,d_{2}} + \sum_{d_{2} \geq 2^{i+1}} v_{\cdot,d_{2}}$$

\noindent
we notice that the interaction of $\hat{u}*\sum_{d_{2}=2^{-i}}^{2^{i+1}}\hat{v}_{\cdot,d_{2}}$ is localized at distance less than $2^{i+1}$ from $P$ while $\hat{u}* \hat{v}_{\cdot,d_{2}}$ is localized at distance $\approx d_{2}$ from $P$, for any $d_{2} \geq 2^{i+1}$. Using (\ref{s114}) also we obtain:

$$||B(u,v)||_{X^{0,-\q,1}_{j, \geq 2^{-i}}} \leq \sum_{d_{2}=2^{-i}}^{2^{i+1}} ||B(u, v_{\cdot, d_{2}})||_{X^{0,-\q,1}_{j, \geq 2^{-i}}} + \sum_{d_{2} \geq 2^{i+1}} ||B(u,v_{\cdot, d_{2}})||_{X^{0,-\q}_{j, \geq 2^{-i}}}$$

$$\leq \sum_{d_{3}=2^{-i}}^{2^{i+1}} \sum_{d_{2}=2^{-i}}^{2^{i+1}} ||B(u, v_{\cdot, d_{2}})||_{X^{0,-\q,1}_{j, d_{3}}}+ \sum_{d_{2} \geq 2^{i+1}} ||B(u,v_{\cdot, d_{2}})||_{X^{0,-\q}_{j, \geq 2^{-i}}}$$

$$\leq \sum_{d_{3}=2^{-i}}^{2^{i+1}} \sum_{d_{2}=2^{-i}}^{2^{i+1}} i^{\q} 2^{i} ||u||_{X^{0,\q,1}_{i}} ||v||_{X^{0,\q}_{j,d_{2}}} + \sum_{d_{2} \geq 2^{i+1}} i^{\q} 2^{i} ||u||_{X^{0,\q,1}_{i}} ||v||_{X^{0,\q}_{j,d_{2}}}$$

$$\leq i^{\frac{3}{2}} 2^{i} ||u||_{X^{0,\q,1}_{i}} ||v||_{X^{0,\q,1}_{j}}$$

In the last line we used the fact that there are $\approx 2i$ values for $d_{3}$. 

\end{proof}

\begin{proof}[Proof of Theorem \ref{tb1}]
The statements in the Theorem are the statements in (\ref{s220}) and (\ref{s130}) when we pass from $X^{0,\q,1}$ to $X^{s,\q,1}$ in all the norms involved. 

\end{proof}

The theory of bilinear estimates does not require any decay of type $\mathcal{D}$. This is necessary in the next section when we want to provide bilinear estimates involving the $Y$ spaces. On the other hand when we involve decay we do it globally, therefore it is going to affect the $X^{s,\q}$ spaces too. The question one should ask is whether the bilinear estimates are the same if we involve $\mathcal{D}\mathcal{R}X^{s,\q}$ spaces instead of  $\mathcal{R}X^{s,\q}$ in Theorem \ref{tb1}. The answer is provided in the following Theorem.

\begin{t1} \label{td1}

a)  Assume that $i \leq j$. We have the following estimates:

\beq \label{b30}
||B(u,v)||_{\mathcal{D}\mathcal{R}X^{s,-\q,1}_{k}} \leq j^{\frac{3}{2}} 2^{(1-s)i} 2^{(k-j)(s-\q-\e)} ||u||_{\mathcal{D}\mathcal{R}X^{s,\q,1}_{i}} ||v||_{\mathcal{D}\mathcal{R}X^{s,\q,1}_{j}}
\eeq

b) Assume that $5i \leq j$. We have the following bilinear estimates:

\beq \label{b31}
||B(u,v)||_{\mathcal{D}\mathcal{R}X^{s,-\q,1}_{j, \geq 2^{-i}}} \leq i^{\frac{3}{2}} 2^{(1-s)i} ||u||_{\mathcal{D}\mathcal{R}X^{s,\q,1}_{i}} ||v||_{\mathcal{D}\mathcal{R}X^{s,\q,1}_{j, \geq 2^{-i}}}
\eeq

Both estimates (\ref{b30}) and (\ref{b31}) remain valid if $B(u,v)$ is replaced by $B(\bar{u},v)$ or $B(u,\bar{v})$. 

\end{t1}

\begin{proof}

Another way of stating the result of this theorem is that all results in Theorem \ref{tb1} hold true if the outcome of the interaction is localized at the high frequency and in the case of high - high frequency interactions with outcome at lower frequency we have to replace the factor $2^{(k-j)s}$ by  $2^{(k-j)(s-\q-\e)}$. 

Let's deal with the first case, i.e. when the outcome gets localized at high frequency. In this case $\mathcal{D}$ acts the same way as the multiplication with $d_{j}=(1+\frac{|x|^{2}}{\mu+2^{2j}})^{\frac{1}{4}+\frac{\e}{2}}$ in the physical space (in terms of $L^{2}$ estimates). We have:

$$d_{j} \cdot B(u_{i},v_{j}) = B(u_{i}, d_{j} v_{j}) + B(u_{i},d_{j}) v_{j}$$

We have $\nabla d_{j}=(\frac{2x_{1}}{2^{2j}},\frac{2x_{2}}{2^{2j}})(1+\frac{|x|^{2}}{2^{2j}})^{-1} d_{j}=\overrightarrow{m}_{j}d_{j}$. A straightforward computations gives us that $||\overrightarrow{m}_{j}||_{L^{\infty}} \leq 2^{-j}$. 

For the term $B(u_{i}, d_{j} v_{j})$ we apply the theory of bilinear estimates we have developed so far. There is a potential difficulty since  $d_{j} v_{j}$ is not localized anymore at frequency $2^{j}$. We learned in the argument for Lemma \ref{decay} that $d_{j} v_{j}$ is essentially localized at frequency $2^{j}$. Quantitatively we have shown there that $(d_{j} v_{j})_{k}$ decays rapidly with respect to $|k-j|$. This is enough to perform every necessary summation with respect to $k$.

For the term $B(u_{i},d_{j})v_{j}$ we apply apply the same strategy, just that the situation is far more simple since we do not have the gradient on the high frequency and moreover we have a gain of a $2^{-j}$ from the $L^{\infty}$ norm of $\overrightarrow{m}_{j}$. 

If we have a high-high frequency interaction ($|i-j| \leq 1$) with outcome at a lower frequency then we have to deal with a term of type $d_{k} S_{k} B(u_{j},v_{j})$. The above argument gives us the estimate:

$$||S_{l}d_{j}B(u_{i},v_{j})||_{\mathcal{R}X^{s,-\q,1}} \leq j^{\frac{3}{2}} 2^{(1-s)j} 2^{-|l-j|s} ||u||_{\mathcal{D}\mathcal{R}X^{s,\q,1}_{i}} ||v||_{\mathcal{D}\mathcal{R}X^{s,\q,1}_{j}}$$

Making use of the result in (\ref{dec50}) we obtain:

$$||d_{k}S_{k}B(u_{i},v_{j})||_{\mathcal{R}X^{s,\q,1}}\leq C_{N} 2^{(j-k)(\q+\e)} \sum_{l \in \N}  2^{-|l-k|N}  ||S_{l}d_{j}B(u_{i},v_{j})||_{\mathcal{R}X^{s,-\q,1}} $$

$$ \leq C_{N} 2^{(j-k)(\q+\e)} \sum_{l \in \N} 2^{-|l-k|N} j^{\frac{3}{2}} 2^{(1-s)j} 2^{-|l-j|s} ||u||_{\mathcal{D}\mathcal{R}X^{s,\q,1}_{i}} ||v||_{\mathcal{D}\mathcal{R}X^{s,\q,1}_{j}} $$

$$ \leq j^{\frac{3}{2}} 2^{(1-s)j} 2^{(k-j)(s-\q-\e)} ||u||_{\mathcal{D}\mathcal{R}X^{s,\q,1}_{i}} ||v||_{\mathcal{D}\mathcal{R}X^{s,\q,1}_{j}}$$

We have skipped quite a few steps in this argument. The reason we did so is in order to spare space and avoid redundancy. For instance Lemma \ref{dec50} gives us an estimate in $L^{2}$ and we use it directly at the level of $X^{s,\q,1}$. This could be done rigourously by preparing an analogue of that Lemma for $X^{s,\q,1}$; we did something similar in the section dedicated to decay.

\end{proof}

\vspace{.1in}

\section{Bilinear estimates involving the {\mathversion{bold}$Y$} spaces}

\vspace{.1in}

In the previous section we have just seen that the theory of bilinear estimates cannot be completely closed in the $X^{s,\q}$ spaces. This is the reason for introducing a more refined structure to measure our solutions, namely the wave-packet one. We concluded that the interactions causing problems in the $X^{s,\q,1}$ theory are the low-high ones. This is why we need to complete Theorem \ref{tb1} with a result for this particular case. As we did there, we assume that $B$ is of type (\ref{bf}). 

\begin{t1} \label{tb2}

Assume we have $5i \leq j$. We have the bilinear estimates:

\beq \label{b7}
||B(u,v)||_{\mathcal{R}\mathcal{D}W^{s}_{j}} \leq i^{\frac{3}{2}} 2^{(1-s)i} ||u||_{\mathcal{R}\mathcal{D}Z^{s}_{i}} ||v||_{\mathcal{R}\mathcal{D}Z^{s}_{j}}
\eeq

The estimate remains valid if $B(u,v)$ is replaced by $B(\bar{u},v)$ or $B(u,\bar{v})$. 

\end{t1}

 In what follows we make few important remarks for the rest of this section. 

\begin{r3} 
We work under the hypothesis that $5i \leq j$. 
\end{r3}

The result in (\ref{s130}) shows that it is fine to use the $X^{s,\q,1}$ structure to measure the low frequency and part of the high frequency (both input and output) at distance greater than $2^{-i}$ from $P$. Thus we shall obtain estimates for:

\beq \label{t1}
X_{i}^{0,\q,1} \cdot Y_{j, \leq 2^{-i}} \rightarrow \mathcal{Y}_{j, \leq 2^{-i}} \ + \ X^{0,\q,1}_{j,\geq 2^{-i}}; \  X_{i}^{0,\q,1} \cdot X^{0,\q,1}_{j, \geq 2^{-i}} \rightarrow \mathcal{Y}_{j, \leq 2^{-i}} 
\eeq

 We also need the corresponding estimates when we involve conjugates of these spaces. The condition $5i \leq j$ implies that the the low frequency does not see the curvature of the parabola at the high frequency, in other words the parabola at high frequency is flat in these interactions. This is why the estimates for $B(\bar{u}_{i},v_{j})$ are similar to the ones for $B(u_{i},v_{j})$. 
 
 If we have to deal with $B(u_{i},\bar{v}_{j})$, a simple geometric argument shows that the interaction is localized at high frequency and in a region with $\tau \leq 0$. This makes these estimates weaker than the ones in (\ref{t1}).

\begin{r3}
Once we get one of the estimates in (\ref{t1}), we trivially get the corresponding ones with conjugate spaces. 
\end{r3}

Our spaces involve rotations, therefore:

\begin{r3}
We use the principles in (\ref{t2}) and (\ref{t3}) in dealing with rotations.  

\end{r3}

We have to involve and recover decay in these estimates. We prove:

$$||B(u,v)||_{\mathcal{R}W^{s}_{j}} \leq i^{\frac{3}{2}} 2^{(1-s)i} ||u||_{\mathcal{R}\mathcal{D}Z^{s}_{i}} ||v||_{\mathcal{R}Z^{s}_{j}}$$

\noindent
and the similar ones. In the end we obtain the estimates with decay on all terms by a similar argument as in Theorem \ref{td1}.  

\begin{r3}

We first prove the estimates without involving decay on the bilinear term and on the high frequency. But we do involve decay on the low frequency. 

\end{r3}

The structure of this section is the following:

- continue with a few definitions;

- provide estimates for interactions between $Y$ and $\mathcal{D}L^{2}$ - spaces;

- analyze the geometry of interactions;

- provide the bilinear estimates in Theorem \ref{tb2}.   

\vspace{.1in}

We record a change in the geometry, namely  we want to work with estimates in strips of width $2^{-i}$.  For this recall the definition of $\Xi^{i}_{n}$ and the corresponding $\phi^{i}_{\xi}$, see (\ref{s120}). 
 
We consider a system of functions $(\phi_{l}^{\tau})_{l \in \Z}$ to be smooth approximations of $\chi_{[l-\q,l+\q]}$ and with the standard partition of unity for $\R$ property. 
 
For each $\xi \in \Xi^{i}$ and $l \in \Z$ we define $g_{\xi,l}$ by  $\hat{g}_{\xi,l}=\phi^{i}_{\xi} \cdot \phi_{l}^{\tau} \cdot \hat{g}$. The support of $\hat{g}_{\xi, l}$ is approximately a tube centered at $(\xi,l)$ and of size $2^{-i} \times 2^{-i} \times 1$, the last one being in the $\tau$ direction. Since the distance of these tubes will play an important role, sometimes it would be convenient if we were able to work with $(g_{\xi, \xi^{2}+l})_{\xi \in \Xi^{i}, l \in \Z}$ instead. The only problem is that it is not guaranteed that $\xi^{2} \in \Z$ for all $\xi \in \Xi^{i}$. Of course we could change the way we cut in the $\tau$ direction, but this would complicate notations even more. We choose instead to ignore that $\xi^{2}$ may not be integer, and go on and use $g_{\xi, \xi^{2}+l}$. It will be obvious from the argument that this does not affect in any way the rigorousness of the proof. The last notation we introduce is $g_{\xi, \xi^{2} \pm l}=g_{\xi, \xi^{2} + l}+g_{\xi, \xi^{2} - l}$.

For $k \in \N$, we define $g_{i, k2^{-i}}$ to be the part of $g_{i}$ whose Fourier transform is localized in $\{ (\xi,\tau): |\tau - \xi^{2}| \in [k-1,k+1] \} \subset \{(\xi,\tau): d((\xi,\tau),P) \approx k 2^{-i} \}$. This inclusion should rather be seen in a strict way. Notice that if we used $d$ instead of $k2^{-i}$ then we obtain a different localization and we hope this will not create confusions. To be more suggestive about which way we go, 

\begin{r3} \label{re}
We choose to run $k$'s over a discrete scale of values and $d$'s over a dyadic scale.
\end{r3}

 This is the case when, as mentioned in the introduction, we may choose to localize in a linear way rather than a dyadic way with respect to the distance to $P$ when it is more convenient. Related to the above notations, we can easily define $S_{i,2^{-i}k}$, $X_{i,2^{-i}k}^{0,\q}$ and similar entities.
 
For $k \leq 2^{2i-2}$ we obtain a new decomposition of $g_{i, k 2^{-i}}$:

$$g_{i, k 2^{-i}} = \sum_{n} \sum_{\xi \in \Xi^{i}_{n}} g_{\xi, \xi^{2} \pm k}$$

Notice that the $\xi$'s $\in \Xi^{i}_{n}$ involved in the above summation have $|\xi| \approx 2^{i}$ since we deal with the part of $\hat{g}$ close to $P$.

For the part of $\hat{g}$ supported away from $P$ we come with a different decomposition. In this region it turns out that the important parameter is the distance to the $\tau$ axis. We have the decomposition: 

\beq \label{aux1}
g_{i, \geq 2^{i-2}} = \sum_{n} \sum_{\xi \in \Xi^{i}_{n}} \sum_{l \in I_{\xi}} g_{\xi,l}
\eeq

\noindent
where $I_{\xi}=\{l \in \Z: 2^{2i-2} \leq |l-\xi^{2}| \leq 2^{2i+2} \}$.

\vspace{.1in}

\subsection{Basic estimates}

\noindent
\vspace{.1in}

This section is concerned with providing results of type $Y \cdot \mathcal{D}L^{2} \rightarrow \mathcal{Y}$, $Y \cdot \mathcal{D} L^{2} \rightarrow L^{2}$ and $L^{2} \cdot \mathcal{D}L^{2} \rightarrow \mathcal{Y}$.

We make the convention that whenever we specify the sizes of a tube in the frequency space, the last size is the one in the $\tau$ direction.  

 We localize $\hat{g}$ on a scale $2^{-i} \times 2^{-i} \times 1$, hence the dual scale to localize in the physical space is $2^{i} \times 2^{i} \times 1$. Recall that the system of cubes $(Q^{m}_{i})_{m \in \Z^{2}}$ is a partition of $\R^{2}$ with the properties: $Q^{m}_{i}$ is centered at $(2^{i}m_{1},2^{i}m_{2})$ and has sizes $2^{i} \times 2^{i}$. Associated to this systems we build a partition $(Q_{i}^{m,l})_{(m,l) \in Z^{3}}$ of $\R^{3}$ defined by:

$$Q_{i}^{m,l}=\cup_{t \in [l,l+1]} Q^{m}_{i} \times \{ t \}= Q^{m}_{i} \times [l,l+1]$$

\begin{l1}

Let $g \in L^{2}$ such that $\hat{g}$ is supported in a tube of size $2^{-i} \times 2^{-i} \times 1$. We have the estimate:

\beq \label{c0}
\sum_{m} ||g||^{2}_{L^{\infty}(Q^{m,l}_{i})} \leq 2^{-2i} ||g||^{2}_{L^{2}}
\eeq

\end{l1}

\begin{proof}
The support of $\hat{g}$ is a tube with volume $2^{-2i}$ therefore we have:

\beq \label{c101}
||g||_{L^{\infty}(Q^{m,l}_{i}))} \leq 2^{-i} \sum_{(m',l') \in Z^{3}} C_{N} \langle (m,l)-(m',l') \rangle^{-N} ||g||_{L^{2}(Q^{m',l'}_{i})}
\eeq

If we chose $N \geq 4$, then we use Cauchy-Schwartz and estimate:

$$||g||^{2}_{L^{\infty}(Q^{m,l}_{i}))} \leq  2^{-2i} \sum_{(m',l') \in \Z^{3}} C^{2}_{N} \langle (m,l)-(m',l') \rangle^{-N} ||g||^{2}_{L^{2}(Q^{m',l'}_{i})} $$

We can perform the summation with respect to $(m,l)$:

$$\sum_{(m,l)} ||g||^{2}_{L^{\infty}(Q^{m,l}_{i}))} \leq  2^{-2i} \sum_{m,l} \sum_{m',l'}   \langle (m,l)-(m',l') \rangle^{-N} ||g||^{2}_{L^{2}(Q^{m',l'}_{i} )} \leq 2^{-2i} ||g||^{2}_{L^{2}}$$

In the last line we use again the fact that if $N \geq 4$, then we have:

$$\sum_{m,l}  \langle (m,l)-(m',l') \rangle^{-N} \leq C$$

This is enough to justify the claim.

\end{proof}

\begin{l1} \label{w1}

Let $g \in \mathcal{D}L^{2}$ such that $\hat{g}$ is supported at frequency $2^{i}$ in a tube of size $2^{-i} \times 2^{-i} \times 1$. We have the estimate:

\beq \label{c1} 
\langle m \rangle^{1+2\e}||g||^{2}_{L^{\infty}(Q^{m,l}_{i})} \leq 2^{-2i} C_{N} \sum_{m',l'}  \langle (m,l)-(m',l') \rangle^{-N} ||g||^{2}_{\mathcal{D}L^{2}(Q_{i}^{m',l'} )}
\eeq

\end{l1}

\begin{proof}

Making use of (\ref{c101}) we can continue with:

$$\langle m \rangle^{1+2\e}||g||^{2}_{L^{\infty}(Q_{i}^{m,l})} \leq $$

$$2^{-2i}  C^{2}_{N} \sum_{m',l'}  \langle m \rangle^{1+2\e} \langle (m',l')-(m,l) \rangle^{-N} ||g||^{2}_{L^{2}(Q^{m',l'}_{i})} \leq$$

$$2^{-2i} C^{2}_{N} \sum_{m',l'}  \langle m' \rangle^{1+2\e} \langle (m',l')-(m,l) \rangle^{-N+1+2\e} ||g||^{2}_{L^{2}(Q^{m',l'}_{i})} \leq$$

$$2^{-2i} C^{2}_{N}  \sum_{m',l'}  \langle (m',l')-(m,l) \rangle^{-N+\q+\e} ||g||^{2}_{\mathcal{D}L^{2}(Q^{m',l'}_{i} )} $$

This is enough to justify the claim.

\end{proof}

In this section we always work with $g=g_{\xi, k2^{-i}}$ or $g=g_{\xi,l}$. 

 We want $\hat{f}_{\eta,\leq 2^{-i}} * \hat{g}$ to be supported at distance less than $2^{-i}$ from $P$. $\hat{f}_{\eta,\leq 2^{-i}}$ and $\hat{f}_{\eta,\leq 2^{-i}} * \hat{g}$ are measured on family of tubes which are different; first family is associated to $\eta$ and the second one to $\eta+\xi$. We need to analyze some of the geometry of the these two families.

By fairly simple geometrical arguments we can conclude the following:
 
 - $T_{\eta+\xi}^{m,l} \cap T_{\eta}^{m',l'} \ne \emptyset \Rightarrow l=l'$
 
 - $T_{\eta+\xi}^{m,l} \cap T_{\eta}^{m',l} \ne \emptyset \Leftrightarrow |m-m'+t\xi| \leq 2$ for some $t \in [l,l+1]$

 - if $T_{\eta+\xi}^{m,l} \cap T_{\eta}^{m',l} \ne \emptyset$, then it is approximately a tube of length $2^{j-i}$
 
 In the last observation we think of the intersection as a subtube of either tubes, and by its length we mean the size of the subtube in the direction of the longest size of either original tubes. We introduce:

- $A^{m,l}=\{m': T_{\eta+\xi}^{m,l} \cap T_{\eta}^{m',l} \ne \emptyset \}$

- $B^{m,m',l}= \{m'': Q^{m'',l}_{i} \cap T_{\eta+\xi}^{m,l} \cap T_{\eta}^{m',l} \ne \emptyset \}$

- $\g_{m}^{l}(m') : A^{m,l} \rightarrow \N \setminus \{0\}$ defined by:

$$\g_{m}^{l}(m')=\{\sup 2^{i-j}||p||: (p,t) \in T_{\eta+\xi}^{m,l} \cap T_{\eta}^{m',l} \}$$

The idea behind the $\g_{m}$ function is that the points in $T_{\eta+\xi}^{m,l} \cap T_{\eta}^{m',l}$ (if nonempty) are of the form $(p,t) \in \R^{2} \times \R$ and the variation of $||p||$ in the intersection is at most  $2^{j-i}$ since the subtubes have length $\approx 2^{j-i}$. It is easy to check the following consequences:

- if $\g_{m}^{l}(m') \leq 2$ then $\forall (p,t) \in T_{\eta+\xi}^{m,l} \cap T_{\eta}^{m',l}$ we have $||p|| \leq 2^{j-i+2}$

- if $\g_{m}^{l}(m') \geq 3$ then $\forall (p,t) \in T_{\eta+\xi}^{m,l} \cap T_{\eta}^{m',l}$ we have $||p|| \in [2^{j-i}(\g_{m}^{l}(m')-1),2^{j-i}(\g_{m}^{l}(m')+1)]$.

We need a simple result about the function $\g_{m}^{l}(m')$:

\begin{l1} For each $l \in \Z$ and $m' \in \Z^{2}$ we have:

\beq \label{g1} 
\sum_{m: m' \in A^{m,l}} \g_{m}^{l}(m')^{-1-2\e} \leq C_{\e} 
\eeq

\end{l1}

\begin{proof}
Things should not be seen as too complicated in the statement above. We simply fix $m'$ (equivalent to fixing $T^{m',l}_{\eta}$), collect all $m$'s for which $T^{m,l}_{\eta+\xi}$ intersect $T^{m',l}_{\eta}$ and perform the summation above over this range.

The function $h(p,t)=||p||$ defined on $T^{m',l}_{\eta}$ attains a minimum at one point, let's call it $p_{0}$. Reminding that $T_{\eta+\xi}^{m,l} \cap T_{\eta}^{m',l}$, if not void, is approximately a subtube of  $T^{m',l}_{\eta}$ of length $2^{j-i}$, we can easily conclude that:

$$\g_{m}^{l}(m')^{2} \approx ||2^{i-j}p_{0}||^{2}+k^{2}+1$$

\noindent
for some $k \in \N$ with $k \leq 2^{i}$. More exactly, $T_{\eta+\xi}^{m,l} \cap T_{\eta}^{m',l}$ is the subtube of $T^{m',l}_{\eta}$ at distance $ \approx k 2^{j-i}$ from $p_{0}$, or in other words the points in the subtube of $T^{m',l}_{\eta}$ at distance $ \approx k 2^{j-i}$ from $p_{0}$. We may have at most two subtubes $T_{\eta+\xi}^{m,l} \cap T_{\eta}^{m',l}$ at distance $k 2^{i-j}$ from $p_{0}$, hence we may conclude that:

$$\sum_{m: m' \in A^{m,l}} \g_{m}^{l}(m')^{-1-2\e} \leq \sum_{k=1}^{2^{i}} (||2^{i-j}p_{0}||^{2}+k^{2}+1)^{-\q-\e} \leq C_{\e}$$

\end{proof}

\begin{l1}

In the same hypothesis as in Lemma \ref{w1} we have the estimate:

\beq \label{c2}
\sum_{m' \in A^{m,l}} \g_{m}^{l}(m')^{1+2\e} \left( \sum_{m'' \in B^{m,m',l}} ||g||_{L^{\infty}(Q_{i}^{m'',l} )}\right)^{2} \leq 2^{-2i}||g||^{2}_{\mathcal{D}L^{2}}
\eeq

\end{l1}

\begin{proof}

For those $m' \in A^{m,l}$ for which $\g_{m}^{l}(m') \leq 2$ (if there are any) we obviously get:

$$\sum_{m'' \in B^{m,m',l}} ||g||_{L^{\infty}(Q_{i}^{m'',l})} \leq  \left( \sum_{m'' \in B^{m,m',l}} \langle m'' \rangle^{1+2\e} ||g||^{2}_{L^{\infty}(Q_{i}^{m'',l})} \right)^{\q} $$

 There are at most $5$ possible $m' \in A^{m,l}$ for which $\g_{m}^{l}(m') \leq 2$.

For those $m' \in A^{m,l}$ for which $\g_{m}^{l}(m') \geq 3$ we proceed as follows. $T_{\eta+\xi}^{m,l} \cap T_{\eta}^{m',l}$ is a subtube of length $2^{j-i}$, hence the cardinality of $B^{m,m',l}$ is $ \approx 2^{j-2i}$. For each fixed $m' \in A^{m,l}$ we have:

$$\sum_{m'' \in B^{m,m',l}} ||g||_{L^{\infty}(Q_{i}^{m'',l})} \leq 2^{\frac{j-2i}{2}} \left( \sum_{m'' \in B^{m,m',l}} ||g||^{2}_{L^{\infty}(Q_{i}^{m'',l})} \right)^{\q} \leq$$

$$ 2^{\frac{j-2i}{2}} (2^{j-2i} \g_{m}^{l}(m'))^{-\q-\e} \left( \sum_{m'' \in B^{m,m',l}} \langle m'' \rangle^{1+2\e} ||g||^{2}_{L^{\infty}(Q_{i}^{m'',l})} \right)^{\q} \leq$$

$$ \g_{m}^{l}(m')^{-\q-\e} \left( \sum_{m'' \in B^{m,m',l}} \langle m'' \rangle^{1+2\e} ||g||^{2}_{L^{\infty}(Q_{i}^{m'',l})} \right)^{\q}$$

\noindent
where we make use of the fact that $j-2i \geq 0$ and $\e$ is positive. At this time we can perform the summation with respect to $m'$:

$$\sum_{m' \in A^{m,l}} \g_{m}^{l}(m')^{1+2\e} \left(  \sum_{m'' \in B^{m,m',l}} ||g||_{L^{\infty}(Q_{i}^{m'',l})} \right)^{2} \leq$$ 

$$\sum_{m' \in A^{m,l}} \sum_{m'' \in B^{m,m',l}} \langle m'' \rangle^{1+2\e} ||g||^{2}_{L^{\infty}(Q_{i}^{m'',l})}$$

The family of $(Q^{m'',l})_{m'' \in B^{m,m',l}, m' \in A^{m,l}}$ does not contain repeated cubes, therefore we can make use of (\ref{c1}) to get:

$$\sum_{m' \in A^{m,l}} \g_{m}^{l}(m')^{1+2\e} \left(  \sum_{m'' \in B^{m,m',l}} ||g||_{L^{\infty}(Q_{i}^{m'',l})} \right)^{2} \leq$$ 

$$2^{-2i} \sum_{m' \in A^{m,l}} \sum_{m'' \in B^{m,m',l}} ||g||^{2}_{\mathcal{D}L^{2}(Q_{i}^{m'',l})} \leq 2^{-2i}||g||^{2}_{\mathcal{D}L^{2}}$$

\end{proof}

\begin{p2} Let $f$ and $g$ be two functions with the following properties: $f=f_{\eta, \leq 2^{-i}} \in Y_{j}$, $|\eta| \approx 2^{j}$, $g \in \mathcal{D}L^{2}$, $\hat{g}$ is supported at frequency $2^{i}$ in a tube of size $2^{-i} \times 2^{-i} \times 1 (\xi \times \tau)$, then we have the estimates:

\beq \label{m8}
|| f \cdot  g||_{\mathcal{Y}_{j}} \leq 2^{-j}||f||_{Y_{j}} ||g||_{\mathcal{D}L^{2}}
\eeq

\beq \label{m9}
|| f \cdot  g||_{L^{2}} \leq  2^{-\frac{i+j}{2}} ||f||_{Y_{j}} ||g||_{L^{2}}
\eeq

\end{p2}

\begin{proof}

For a particular $m$, $m' \in A^{m,l}$ and $m'' \in B^{m,m',l}$ the intersection $Q^{m'',l}_{i} \cap T^{m,l}_{\eta+\xi}$ is included in a rectangular parallelepiped of sizes $2^{i} \times 2^{i} \times 2^{i-j}$ (last one in the $t$ direction). Therefore we have:

$$||f \cdot g||_{L_{t}^{1}L_{x}^{2}(T^{m,l}_{\eta+\xi} \cap T^{m',l}_{\eta})}=\sum_{m'' \in B^{m,m',l}} ||f \cdot g||_{L_{t}^{1}L_{x}^{2}(Q^{m'',l}_{i} \cap T^{m,l}_{\eta+\xi})} \leq$$

$$ 2^{i-j} \sum_{m'' \in B^{m,m',l}} ||f \cdot g||_{L_{t}^{\infty}L_{x}^{2}(Q^{m'',l}_{i} \cap T^{m,l}_{\eta+\xi})} \leq$$

$$ 2^{i-j} \sum_{m'' \in B^{m,m',l}} ||f ||_{L_{t}^{\infty}L_{x}^{2}(Q^{m'',l}_{i})} ||g||_{L^{\infty}(Q^{m'',l}_{i})} \leq $$

$$ 2^{i-j}  ||f ||_{L_{t}^{\infty}L_{x}^{2}(T^{m',l}_{\eta})} \sum_{m'' \in B^{m,m',l}} ||g||_{L^{\infty}(Q^{m'',l}_{i})} $$

We can go on and perform the summation with respect to $m' \in A^{m,l}$:

$$||f \cdot g ||_{L_{t}^{1}L_{x}^{2}(T^{m,l}_{\eta +\xi})} \leq$$ 

$$2^{i-j} \sum_{m' \in A^{m,l}}  ||f ||_{L_{t}^{\infty}L_{x}^{2}(T^{m',l}_{\eta})} \sum_{m'' \in B^{m,m',l}} ||g||_{L^{\infty}(Q^{m'',l}_{i})} \leq$$

$$2^{i-j} \left (\sum_{m' \in A^{m,l}} \g_{m}^{l}(m')^{-1-2\e}  ||f||^{2}_{L_{t}^{\infty}L_{x}^{2}(T^{m',l}_{\eta})} \right)^{\q} \cdot$$

$$ \left (\sum_{m' \in A^{m,l}} \g_{m}^{l}(m')^{1+2\e} \left( \sum_{m'' \in B^{m,m',l}} ||g||_{L^{\infty}(Q^{m'',l}_{i})} \right)^{2} \right)^{\q} \leq$$

$$ 2^{-j} \left (\sum_{m' \in A^{m,l}} \g_{m}^{l}(m')^{-1-2\e}  ||f||^{2}_{L_{t}^{\infty}L_{x}^{2}(T^{m',l}_{\eta})} \right)^{\q} \cdot ||g||_{\mathcal{D}L^{2}}$$

In the last line we have used the result in (\ref{c2}) . The norm in $\mathcal{Y}_{j}$ is an $l^{2}_{m,l}$ of the norms above:

$$||f \cdot g ||^{2}_{\mathcal{Y}_{j}} \approx \sum_{(m,l) \in \Z^{3}} ||f \cdot g ||^{2}_{L_{t}^{1}L_{x}^{2}(T^{m,l}_{\eta + \xi})} \leq$$

$$ 2^{-2j}||g||^{2}_{\mathcal{D}L^{2}} \sum_{m,l} \sum_{m' \in A^{m,l}} \g_{m}^{l}(m')^{-1-2\e} ||f||^{2}_{L_{t}^{\infty}L_{x}^{2}(T^{m',l}_{\eta})} \leq $$

$$2^{-2j} ||g||^{2}_{\mathcal{D}L^{2}} \sum_{m',l} ||f||^{2}_{L_{t}^{\infty}L_{x}^{2}(T^{m',l}_{\eta})} \sum_{m: m' \in A^{m,l}} \g_{m}^{l}(m')^{-1-2\e} \leq 2^{-2j}||g||^{2}_{\mathcal{D}L^{2}} ||f||^{2}_{Y_{j}}$$

In the last estimates we have made use of (\ref{g1}).

The $L^{2}$ estimates are much easier. For each $(m,l) \in \Z^{3}$ let's denote by $C^{m,l}=\{m' \in Z^{2}: Q^{m',l} \cap T^{m,l}_{\eta} \ne \emptyset \}$. Then we have:

$$||f \cdot g||^{2}_{L^{2}(T^{m,l}_{\eta})}=\sum_{m' \in C^{m,l}} ||f \cdot g||^{2}_{L^{2}(Q_{i}^{m',l} \cap T^{m,l}_{\eta})} \leq$$

$$ 2^{i-j} \sum_{m' \in C^{m,l}} ||f||^{2}_{L_{t}^{\infty}L^{2}_{x}(T^{m,l}_{\eta})} ||g||^{2}_{L^{\infty}(Q_{i}^{m',l})} \leq $$

$$2^{i-j} ||f||^{2}_{L_{t}^{\infty}L^{2}_{x}(T^{m,l}_{\eta})} \sum_{m' \in C^{m,l}} ||g||^{2}_{L^{\infty}(Q_{i}^{m',l})} \leq   2^{-i-j} ||f||^{2}_{L_{t}^{\infty}L^{2}_{x}(T^{m,l}_{\eta})} ||g||^{2}_{L^{2}}$$

In the last estimate we have used the result in (\ref{c0}).

We sum the above estimate with respect to $(m,l)$ over $Z^{3}$ to obtain (\ref{m9}).

\end{proof}

From (\ref{m9}) we obtain, by duality, the following result:

\begin{p2} Let $f \in L^{2}$ and and $g$ be two functions with the following properties: $f \in L^{2}$, $\hat{f}$ is supported at frequency $2^{j}$, $\hat{g}$ is supported at frequency $2^{i}$ in a tube of size $2^{-i} \times 2^{-i} \times 1 (\xi \times \tau)$, then we have the estimates:

\beq \label{m6}
|| f \cdot g||_{\mathcal{Y}_{j}} \leq 2^{-\frac{i+j}{2}} ||f||_{L^{2}} ||g||_{L^{2}}
\eeq

\end{p2}

The next Lemma is a geometrical one. We work with $f=f_{\eta, \leq 2^{-i}}$ and $g=g_{\xi^{0},l}$, $\xi^{0} \in \Xi^{i}, l \in \Z$ where $|\eta| \approx 2^{j}$, $|(\xi^{0},l)| \approx 2^{i}$ and recall that $5i \leq j$.

\begin{l1}
$\hat{f} * \hat{g}$ is supported in a region where $|\tau-\xi^{2}| \in 2^{-i}|\eta|[k-1,k+1]$ iff

\beq \label{d1}
|\cos{\a}| \in 2^{-i} |\xi^{0}|^{-1} [k-1,k+1]
\eeq

\noindent
where $\a$ is the angle between $\xi^{0}$ and $\eta$. 

\end{l1}

\begin{proof}

 $\hat{f}$ is supported in a region where $|\tau_{2}-\eta^{2}| \leq 2^{-i}|\eta|$, while $\hat{g}$ is supported in a region where $|\xi-\xi^{0}| \leq 2^{-i}$ and $|\tau_{1}-l| \leq \q$. A generic point in the support of $\hat{f}*\hat{g}$ is of type $(\xi_{1}+\xi_{2},\tau_{1}+\tau_{2})$ where $(\xi_{1},\tau_{1})$ is in the support of $\hat{f}$ and $(\xi_{2},\tau_{2})$ is in the support of $\hat{g}$. We want this point to satisfy $|\tau_{1}+\tau_{2}-(\xi_{1}+\xi_{2})^{2}| \in 2^{-i}|\eta|[k-1,k+1]$.

We have $|\tau_{1} -\xi_{1}^{2}| \leq 2^{2i} \leq 2^{-i}|\eta|$, $\D |\xi_{1}| \approx 2^{-i}$, $\D |\eta| \approx 1$, therefore the condition is equivalent to $|2\eta \cdot \xi^{0}| \in |\eta|2^{-i}[k-1,k+1]$. This implies (\ref{d1}).

\end{proof}

\vspace{.1in}

\subsection{Estimates: {\mathversion{bold}$ \mathcal{D}\mathcal{R}X^{0,\q,1}_{i} \cdot \mathcal{R}Y_{j, \leq 2^{-i}} \rightarrow \mathcal{R}\mathcal{Y}_{j, \leq 2^{-i}}$}}

\noindent
\vspace{.1in}

The first Proposition deals with the case when we have the low frequency input close to $P$.

\begin{p2} If $k \leq 2^{2i-2}$ we have:

\beq \label{a1}
|| v_{j, \leq 2^{-i}} \cdot  u_{i, k 2^{-i}}||_{\mathcal{R}\mathcal{Y}_{j, \leq 2^{-i}}} \leq 2^{-j} ||v_{j, \leq 2^{-i}}||_{\mathcal{R}Y_{j}} \cdot ||u_{i, k 2^{-i}}||_{\mathcal{D}\mathcal{R}L^{2}}
\eeq

\end{p2}

\begin{proof}

We first deal with the case $k=1$ and then use this as a model for the other $k$'s. We decompose:

\beq \label{dec1}
v_{j, \leq 2^{-i}}=\sum_{n_{1}=2^{j+i-1}}^{2^{j+i+1}} v_{n_{1}, \leq 2^{-i}} = \sum_{n_{1}=2^{j-1}}^{2^{j+1}} \sum_{\eta \in \Xi_{n_{1}}} v_{\eta, \leq 2^{-i}}
\eeq

\beq \label{dec2}
u_{i, \leq 2^{-i}}=\sum_{n_{2}= 2^{2i-2}}^{2^{2i+2}} u_{n_{2}, \leq 2^{-i}} =\sum_{n_{2}= 2^{2i-1}}^{2^{2i+1}} \sum_{\xi \in \Xi^{i}_{n_{2}}} u_{\xi, \leq 2^{-i}}
\eeq

Fix $n_{1}$ and $n_{2}$. Take $\xi \in \Xi^{i}_{n_{2}}$ and $\eta \in \Xi_{n_{1}}$ and denote by $\a$ their angle. We want $\hat{v}_{\eta, \leq 2^{-i}} * \hat{u}_{\xi, \leq 2^{-i}}$ to be supported at distance less than $2^{-i}$ from $P$. Using (\ref{d1}) we obtain the following condition on $\a$:

$$|\cos{\a}| \leq 2^{-2i}$$ 

This condition suggests splitting $\Xi_{n_{1}}$ and $\Xi^{i}_{n_{2}}$ in $2^{-2i}$ angular subsets, i.e. such that two elements in the same subset make an angle less then $2^{-2i}$ and two elements from different subsets make an angle greater than $2^{-2i}$.  

 $\Xi^{i}_{n_{2}}= \cup_{\xi \in \Xi^{i}_{n_{2}}} \{ \xi \}$ is exactly what we need by the definition of  $\Xi^{i}_{n_{2}}$. 

 For each $\xi \in \Xi^{i}_{n_{2}}$ there are $\approx 2^{j-2i}$ $\eta$'s which make an angle less than $2^{-2i}$ with $\xi$ and we denote by $A_{\xi}$ this set. It is obvious that if $\xi \ne \xi'$ then $A_{\xi}$ and $A_{\xi'}$ are disjoint and $\cup _{\xi} A_{\xi} = \Xi_{n_{1}}$. 

The last geometrical detail we have to clarify is the separation of the supports of $\hat{v}_{\eta, \leq 2^{-i}} * \hat{u}_{\xi, \leq 2^{-i}}$ as we vary $\xi$ and $\eta \in A_{\xi}$. 

The support of $\hat{u}_{\xi, \leq 2^{-i}}$ is a tube of sizes $2^{-i} \times 2^{-i} \times 1$ and the (long) axis points in the direction of $\tau$. The support of $\hat{v}_{\eta, \leq 2^{-i}}$ is a parallelepiped of sizes $2^{-i} \times 1 \times 2^{j}$ whose longest side is tangent to $P$. The key property is that we can translate the support of $\hat{u}_{\xi, \leq 2^{-i}}$ so that it is included in the support of $\hat{v}_{\eta, \leq 2^{-i}}$ (by simply translating the center of the first to the center of the second). Therefore the support of $\hat{v}_{\eta, \leq 2^{-i}} * \hat{u}_{\xi, \leq 2^{-i}}$ is a translate of the support of $\hat{v}_{\eta, \leq 2^{-i}}$ by the vector $(\xi, \xi^{2})$. If we keep $\xi$ fixed and take  $\eta \ne \eta'$  both in $A_{\xi}$, then the supports of $\hat{v}_{\eta, \leq 2^{-i}} * \hat{u}_{\xi, \leq 2^{-i}}$ and $\hat{v}_{\eta', \leq 2^{-i}} * \hat{u}_{\xi, \leq 2^{-i}}$ are disjoint. 

If we take $\xi \ne \xi'$ and $\eta \in A_{\xi}$, $\eta' \in A_{\xi'}$ then the supports of $\hat{v}_{\eta, \leq 2^{-i}} * \hat{u}_{\xi, \leq 2^{-i}}$ and $\hat{v}_{\eta', \leq 2^{-i}} * \hat{u}_{\xi', \leq 2^{-i}}$ are in different angular regions, therefore they are disjoint again.

 We can apply (\ref{m8}) to each pair $v_{\eta, \leq 2^{-i}}, u_{\xi, \leq 2^{-i}}$. Using the orthogonality with respect to $\eta \in A_{\xi}$ of the convolution we get:

$$|| \sum_{\eta \in A_{\xi}} v_{\eta, \leq 2^{-i}} \cdot  u_{\xi,\leq 2^{-i}}||^{2}_{\mathcal{Y}_{j, \leq 2^{-i}}} = \sum_{\eta \in A_{\xi}} ||v_{\eta, \leq 2^{-i}} \cdot  u_{\xi,\leq 2^{-i}}||^{2}_{\mathcal{Y}_{j, \leq 2^{-i}}} \leq$$

$$2^{-2j} \sum_{\eta \in A_{\xi}} ||v_{\eta, \leq 2^{-i}}||^{2}_{Y_{j}}  ||u_{\xi,\leq 2^{-i}}||^{2}_{\mathcal{D}L^{2}}$$

Using the orthogonality with respect to $\xi$ of the convolution we get:

$$||v_{n_{1}, \leq 2^{-i}} \cdot u_{n_{2}, \leq 2^{-i}}||^{2}_{\mathcal{Y}_{j, \leq 2^{-i}}} \approx \sum_{\xi} || \sum_{\eta \in A_{\xi}} v_{\eta, \leq 2^{-i}} \cdot  u_{\xi,\leq 2^{-i}}||^{2}_{\mathcal{Y}_{j, \leq 2^{-i}}} \leq$$

$$2^{-2j} \sum_{\xi} \sum_{\eta \in A_{\xi}} ||v_{\eta, \leq 2^{-i}}||^{2}_{Y_{j}}  ||u_{\xi,\leq 2^{-i}}||^{2}_{\mathcal{D}L^{2}}$$

We can perform the summation with respect to $\xi$ on each of the both terms above and use a $\sup$ on the other one:

$$||v_{n_{1}, \leq 2^{-i}} \cdot u_{n_{2}, \leq 2^{-i}}||^{2}_{\mathcal{Y}_{j, \leq 2^{-i}}} \leq 2^{-2j} \sup_{\xi \in \Xi^{i}_{n_{2}}} (\sum_{\eta \in A_{\xi}} ||v_{\eta, \leq 2^{-i}}||^{2}_{Y_{j}}) ||u_{n_{2}, \leq2^{-i}}||^{2}_{\mathcal{D}L^{2}}$$

\noindent
and

$$||v_{n_{1}, \leq 2^{-i}} \cdot u_{n_{2}, \leq 2^{-i}}||^{2}_{\mathcal{Y}_{j, \leq 2^{-i}}} \leq 2^{-2j} ||v_{j, \leq 2^{-i}}||^{2}_{Y_{j}} \sup_{\xi \in \Xi^{i}_{n_{2}}} ||u_{\xi, \leq 2^{-i}}||^{2}_{\mathcal{D}L^{2}}$$

In both cases we can bound the $\sup$ by using $\mathcal{R}$ on the corresponding term:

$$\sup_{\xi \in \Xi^{i}_{n_{2}}} (\sum_{\eta \in A_{\xi}} ||v_{\eta, \leq 2^{-i}}||^{2}_{Y_{j}}) \leq 2^{-2i} ||v_{n_{1}, \leq 2^{-i}}||^{2}_{\mathcal{R}Y_{j}}$$

$$\sup_{\xi \in \Xi^{i}_{n_{2}}} ||u_{\xi, \leq 2^{-i}}||^{2}_{\mathcal{D}L^{2}} \leq 2^{-2i} ||u_{n_{2},\leq 2^{-i}}||^{2}_{\mathcal{R}\mathcal{D}L^{2}}$$

Therefore we get the estimates:

$$||v_{n_{1}, \leq 2^{-i}} \cdot u_{n_{2}, \leq 2^{-i}}||_{\mathcal{Y}_{j, \leq 2^{-i}}} \leq 2^{-i-j} ||v_{n_{1}, \leq 2^{-i}}||_{\mathcal{R}Y_{j}} ||u_{n_{2},\leq 2^{-i}}||_{\mathcal{D}L^{2}}$$

$$||v_{n_{1}, \leq 2^{-i}} \cdot u_{n_{2}, \leq 2^{-i}}||_{\mathcal{Y}_{j, \leq 2^{-i}}} \leq 2^{-i-j} ||v_{n_{1}, \leq 2^{-i}}||_{Y_{j}} ||u_{n_{2},\leq 2^{-i}}||_{\mathcal{R}\mathcal{D}L^{2}}$$

In the end we have to sum with respect to $n_{1}$ and $n_{2}$. $\hat{v}_{n_{1}, \leq 2^{-i}} * \hat{u}_{n_{2}, \leq 2^{-i}}$ is supported in a region with $\tau \approx n^{2}_{1} +(n_{2}2^{-i})^{2}$ and $\D \tau \leq 3n_{1}+2^{-i+2}n_{2} \leq 4n_{1}$.  As a consequence, if we keep $n_{2}$ fixed and take $|n_{1}-n_{1}'| \geq 2$, then the supports of $\hat{v}_{n_{1}, \leq 2^{-i}} * \hat{u}_{n_{2}, \leq 2^{-i}}$ and $\hat{v}_{n'_{1}, \leq 2^{-i}} * \hat{u}_{n_{2}, \leq 2^{-i}}$ are disjoint. This implies:

$$||v_{j, \leq 2^{-i}} \cdot u_{n_{2}, \leq 2^{-i}}||^{2}_{\mathcal{Y}_{j, \leq 2^{-i}}} \approx \sum_{n_{1}} ||v_{n_{1}, \leq 2^{-i}} \cdot u_{n_{2}, \leq 2^{-i}}||^{2}_{\mathcal{Y}_{j, \leq 2^{-i}}} \leq $$

$$ 2^{-2i-2j} \sum_{n_{1}} ||v_{n_{1}, \leq 2^{-i}}||^{2}_{\mathcal{R}Y_{j}} ||u_{n_{2},\leq 2^{-i}}||^{2}_{\mathcal{D}L^{2}} \approx  2^{-2i-2j} ||v_{j, \leq 2^{-i}}||^{2}_{\mathcal{R}Y_{j}} ||u_{n_{2},\leq 2^{-i}}||^{2}_{\mathcal{D}L^{2}}$$

 We do not have orthogonality with respect to $n_{2}$, so we use the trivial estimate:

$$||v_{j, \leq 2^{-i}} \cdot u_{i, \leq 2^{-i}}||_{\mathcal{Y}_{j, \leq 2^{-i}}} \approx \sum_{n_{2} \leq 2^{2i+1}} ||v_{j, \leq 2^{-i}} \cdot u_{n_{2}, \leq 2^{-i}}||_{\mathcal{Y}_{j, \leq 2^{-i}}} \leq $$

$$2^{i} \left( \sum_{n_{2} \leq 2^{2i+1}} ||v_{j, \leq 2^{-i}} \cdot u_{n_{2}, \leq 2^{-i}}||^{2}_{\mathcal{Y}_{j, \leq 2^{-i}}} \right)^{\q} \leq 2^{-j}||v_{j, \leq 2^{-i}}||_{\mathcal{R}Y_{j}} ||u_{n_{2}, \leq 2^{-i}}||_{\mathcal{D}L^{2}}$$

In a similar way we can perform the estimate when we use the rotations on $v$ and obtain:  
$$||v_{j, \leq 2^{-i}} \cdot u_{i, \leq 2^{-i}}||_{\mathcal{Y}_{j, \leq 2^{-i}}} \leq 2^{-j} ||v_{j, \leq 2^{-i}}||_{Y_{j}} ||u_{n_{2}, \leq 2^{-i}}||_{\mathcal{R}\mathcal{D}L^{2}}$$

Making use of the principle stated in (\ref{t2}) we can derive (\ref{a1}) from these two last estimates.  

What changes if $k \ne 1$? We start in a similar manner, namely decompose $v_{j, \leq 2^{-i}}$ as in (\ref{dec1}) and

\beq \label{dec3}
 u_{i, k 2^{-i}}=\sum_{n_{2} \leq 2^{2i+1}} u_{n_{2}, k 2^{-i}} =\sum_{n_{2} = 2^{2i-1}}^{2^{2i+1}} \sum_{\xi \in \Xi^{i}_{n_{2}}} u_{\xi, \xi^{2} \pm k}
\eeq

We fix $n_{1}$ and $n_{2}$. We reduce the estimates at this level to the ones we have just proved by showing that the main geometrical elements are similar. 

Take $\xi \in \Xi^{i}_{n_{2}}$ and $\eta \in \Xi_{n_{1}}$ and denote by $\a$ their angle. We want $\hat{v}_{\eta, \leq 2^{-i}} * \hat{u}_{\xi, \xi^{2} \pm k}$ to be supported at distance less than $2^{-i}$ from $P$. Using (\ref{d1}) we get the following condition on $\a$:

$$|\cos{\a}| \leq 2^{-2i}$$ 

From this point we can use the same argument as in the case $k=1$.

\end{proof}

The next Proposition deals with the case when we have the low frequency input close to $\tau$ axis. 

\begin{p2}  We have 

\beq \label{a5}
|| v_{j, \leq 2^{-i}} \cdot  \sum_{\xi \in \Xi^{i}_{k}} u_{\xi, \geq 2^{i-2}}||_{\mathcal{R}\mathcal{Y}_{j, \leq 2^{-i}}} \leq 
\eeq

$$2^{i-j} k^{-\q}||v_{j, \leq 2^{-i}}||_{\mathcal{R}Y_{j}} \cdot ||\sum_{\xi \in \Xi^{i}_{k}} u_{\xi, \geq 2^{i-2}}||_{\mathcal{D}\mathcal{R}L^{2}}$$

\end{p2}

\begin{proof} We decompose $v_{j, \leq 2^{-i}}$ as in (\ref{dec1}) and using (\ref{aux1}):

$$ \sum_{\xi \in \Xi^{i}_{k}} u_{\xi, \geq 2^{i-2}}= \sum_{\xi \in \Xi^{i}_{k}} \sum_{l \in I_{\xi}} u_{\xi, l}$$

(\ref{d1}) gives us a necessary and sufficient condition for the support of $\hat{v}_{\eta, \leq 2^{-i}} * \hat{u}_{\xi,l}$ to be localized at distance less than $2^{-i}$ from $P$: $ |\cos{\a}| \leq k^{-1}$, where $\a$ is the angle between $\xi$ and $\eta$. We use the fact that in the support of $\hat{u}_{\xi,l}$ we have $|\xi| \approx k 2^{-i}$. 

This suggests to split $\Xi_{n_{1}}$ and $\Xi^{i}_{k2^{-i}}$ in angular sectors of size $\approx k^{-1}$.  $\Xi^{i}_{k2^{-i}}$ comes already with this splitting since that the angle between every two different $\xi \in \Xi^{i}_{k2^{-i}}$ is at least $k^{-1}$. We define $A_{\xi}$ to be the set of $\eta \in \Xi_{n_{1}}$ whose angle $\a$ with $\xi$ satisfies $|\cos{\a}| \leq k^{-1}$. 

We have the same geometrical setup as before. At the numerology level we record the following changes:
 
- the gain from spherical symmetry changes now to $k^{-\q}$ since this is the angular localization. 

- since we have to perform a summation with respect to $l$ (which we did not have to do before) we pick a factor of $2^{i}$. 

Other than that the argument is the same as before.

\end{proof}

The next result sums up the results in the two previous Propositions and provides us with the estimate we wanted in this section. 

\begin{p2} We have 

\beq \label{a7}
|| v_{j, \leq 2^{-i}} \cdot  u_{i}||_{\mathcal{R} \mathcal{Y}_{j, \leq 2^{-i}}} \leq i^{\q} 2^{-j}||v_{j, \leq 2^{-i}}||_{\mathcal{R}Y_{j}} \cdot ||u_{i}||_{\mathcal{D}\mathcal{R}X^{0,\q,1}}
\eeq

\end{p2}

\begin{proof} We decompose:

$$u_{i}= u_{i, \leq 2^{i-2}} + u_{i, \geq 2^{i-2}}=\sum_{k=1}^{2^{2i-2}} u_{i,k 2^{-i}} + \sum_{k} \sum_{\xi \in \Xi^{i}_{k }} \sum_{l \in I_{\xi}} u_{\xi, l}$$

Using the result in (\ref{d3001}) we have:  

$$||u_{i}||^{2}_{\mathcal{R}\mathcal{D}X^{0,\q}} \approx \sum_{k} k||u_{i,k 2^{-i}}||^{2}_{\mathcal{R}\mathcal{D}L^{2}} + 2^{-i} \sum_{k} \sum_{\xi \in \Xi^{i}_{k}} \sum_{l \in I_{\xi}} ||u_{\xi,l}||^{2}_{\mathcal{R}\mathcal{D}L^{2}} $$

We continue with:

$$|| v_{j, \leq 2^{-i}} \cdot  u_{i, \leq 2^{i-2}}||_{\mathcal{R}\mathcal{Y}_{j, \leq 2^{-i}}} \leq \sum_{k=1}^{2^{2i-2}} || v_{j, \leq 2^{-i}} \cdot  u_{i, k 2^{-i}}||_{\mathcal{R}\mathcal{Y}_{j, \leq 2^{-i}}}  \leq$$

$$2^{-j} \sum_{k=1}^{2^{2i-2}} ||v_{j, \leq 2^{-i}}||_{\mathcal{R}Y_{j}} \cdot ||u_{i, k 2^{-i}}||_{\mathcal{D}\mathcal{R}L^{2}} \leq $$

$$ 2^{-j} || v_{j, \leq 2^{-i}}||_{\mathcal{R}Y_{j}} (\sum_{k=1}^{2^{2i-2}} k^{-1})^{\q} \left( \sum_{k=1}^{2^{2i-2}} k ||u_{i,k 2^{-i}}||^{2}_{\mathcal{D}\mathcal{R}L^{2}}\right)^{\q} \leq $$

$$ i^{\q} 2^{-j} || v_{j, \leq 2^{-i}}||_{\mathcal{R}Y_{j}} ||u_{i, \leq 2^{i-2}}||_{\mathcal{D}\mathcal{R}X^{0,\q}}$$

For the second part:

$$|| v_{j, \leq 2^{-i}} \cdot  u_{i, \geq 2^{i-2}}||_{\mathcal{R}\mathcal{Y}_{j, \leq 2^{-i}}} \leq \sum_{k} || v_{j, \leq 2^{-i}} \cdot \sum_{\xi \in \Xi^{i}_{k }} \sum_{l \in I_{\xi}} u_{\xi,l}||_{\mathcal{R}\mathcal{Y}_{j, \leq 2^{-i}}}  \leq$$

$$2^{-j}\sum_{k} 2^{i} k^{-\q} ||v_{j, \leq 2^{-i}}||_{\mathcal{R}Y_{j}} \cdot ||\sum_{\xi \in \Xi^{i}_{k}} \sum_{l \in I_{\xi}} u_{\xi,l}||_{\mathcal{D}\mathcal{R}L^{2}} \leq $$

$$2^{i}2^{-j} || v_{j, \leq 2^{-i}}||_{\mathcal{R}Y_{j}} (\sum_{k} k^{-1})^{\q} \left( \sum_{k} ||\sum_{\xi \in \Xi^{i}_{k}} \sum_{l \in I_{\xi}} u_{\xi,l}||^{2}_{\mathcal{D}\mathcal{R}L^{2}}\right)^{\q} \leq$$

$$ i^{\q} 2^{-j} || v_{j, \leq 2^{-i}}||_{Y_{j}} ||u_{i, \geq 2^{i-2}}||_{\mathcal{D}\mathcal{R}X^{0,\q}}$$

In the end we sum up the two estimates and use the trivial fact that $||u_{i}||_{X^{0,\q}} \leq ||u_{i}||_{X^{0,\q,1}}$ to obtain (\ref{a7}).

\end{proof}

\vspace{.1in}

\subsection{Estimates: {\mathversion{bold}$ \mathcal{R} X^{0,\q,1}_{i} \cdot \mathcal{R}Y_{j, \leq 2^{-i}} \rightarrow \mathcal{R}X^{0,-\q,1}_{j, \geq 2^{-i}}$}}

\noindent
\vspace{.1in}

The first result deals with the case when we have the low frequency input close to $P$ and we measure the part of the output which is close to $P$.

\begin{p2} If $k_{2} \leq 2^{2i-2}$ and $k_{1} \leq 2^{2i-2}$ then we have 

\beq \label{a9}
|| v_{j, \leq 2^{-i}} \cdot  u_{i, 2^{-i} k_{2}}||_{\mathcal{R}X^{0,-\q}_{j, 2^{-i} k_{1}}} \leq 2^{-j} k_{1}^{-\q} ||v_{j, \leq 2^{-i}} ||_{\mathcal{R}Y_{j}} \cdot ||u_{i, 2^{-i} k_{2}}||_{\mathcal{R}L^{2}}
\eeq

\end{p2}

\begin{proof} We decompose $v_{j, \leq 2^{-i}}$ as in (\ref{dec1}) and  $u_{i, k 2^{-i}}$ as in (\ref{dec3}).

Fix $n_{1}$ and $n_{2}$. Take $\xi \in \Xi^{i}_{n_{2}}$ and $\eta \in \Xi_{n_{1}}$ and denote by $\a$ their angle. We want $\hat{v}_{\eta, 2^{-i}} * \hat{u}_{\xi, \xi^{2}+k_{2}}$ to be at distance $\approx 2^{-i} k_{1}$ from $P$. The condition we get from (\ref{d1}) is:

$$|\cos{\a}| \in 2^{-2i} [k_{1}-1, k_{1}+1]$$

We have that $k_{1} \leq 2^{2i-2}$, therefore $\cos{\a} \leq \q$. This is important because it implies that the solution of the above inclusion is an interval of size $\approx 2^{-2i}$. Therefore we go on and split $\Xi_{n_{1}}$ and $\Xi^{i}_{n_{2}}$ in $2^{-2i}$ angular subsets. We use the same kind of decomposition as in the proof of (\ref{a1}).

 $\Xi^{i}_{n_{2}}= \cup_{\xi \in \Xi^{i}_{n_{2}}} \{ \xi \}$ is exactly what we need by the definition of  $\Xi^{i}_{n_{2}}$. 

 For each $\xi \in \Xi^{i}_{n_{2}}$ there are $\approx 2^{j-2i}$ $\eta$'s whose angle $\a$ with $\xi$ satisfies $|\cos{\a}| \in 2^{-2i}[k_{1}-1, k_{1}+1]$ and we denote by $A_{\xi}$ this set. It is obvious that if $\xi \ne \xi'$ then $A_{\xi}$ and $A_{\xi'}$ are disjoint and $\cup _{\xi} A_{\xi} = \Xi_{n_{1}}$. 

From this point on we can copy verbatim the argument in used in proving (\ref{a1}), of course making use of (\ref{m9}). We obtain:

$$|| v_{j, \leq 2^{-i}} \cdot  u_{i, 2^{-i} k_{2}}||_{X^{0,-\q}_{j, 2^{-i} k_{1}}} \approx (2^{j}2^{-i} k_{1})^{-\q} ||S_{j, 2^{-i} k_{1}}( v_{j, \leq 2^{-i}} \cdot  u_{i, 2^{-i} k_{2}})||_{L^{2}} \leq $$

$$(2^{j}2^{-i} k_{1})^{-\q} 2^{-\frac{i+j}{2}}||v_{j, \leq 2^{-i}} ||_{Y_{j}} \cdot ||u_{i, 2^{-i} k_{2}}||_{\mathcal{R}L^{2}}=$$

$$ k_{1}^{-\q} 2^{-j}||v_{j, \leq 2^{-i}} ||_{Y_{j}} \cdot ||u_{i, 2^{-i} k_{2}}||_{\mathcal{R}L^{2}}$$ 

In a similar way we obtain the estimates with rotations on $Y_{j}$.

\end{proof}

The next result deals with the case when we have the low frequency input close to $P$ and we measure the part of the output which is away from $P$. 

\begin{p2} If $k_{2} \leq 2^{2i-2}$ and $d \geq 2^{2i-2}$ we have 

\beq \label{a50}
|| v_{j, \leq 2^{-i}} \cdot  u_{i, 2^{-i} k_{2}}||_{X^{0,-\q}_{j, 2^{-i} d}} \leq 2^{-j} ||v_{j, \leq 2^{-i}}||_{Y_{j}} \cdot ||u_{i, 2^{-i} k_{2}}||_{L^{2}}
\eeq

\end{p2}

\begin{proof}

What is the particularity of this case? When we localize  the support of $\hat{u}_{\eta, \leq 2^{-i}} * \hat{u}_{\xi,\xi^{2}+ k_{2}}$  at distance $\geq 2^{i-2}$ from $P$ we do not have anymore that the angle $\a$ between $\eta$ and $\xi$ satisfies $|\cos{\a}| \leq \q$. Therefore we cannot conclude that if $|\cos{\a}|$ is in an interval of size $2^{-2i}$ then so does $\a$. 

So we have to come up with a different way of organizing the interacting elements, the main reason being to bring some sort of orthogonality into play. For $\t=\frac{\pi}{2} l 2^{-2i}$ with $l \in \{0,1,...,2^{2i+4}-1\}$ we define 

$$\Xi^{i}_{\t}=\{\xi=(r, \t); r=n 2^{-i} \ \mbox{for} \ n \leq 2^{2i+1}  \}$$ 

We decompose $v_{j, \leq 2^{-i}}$ as in (\ref{dec1}) and:

$$ u_{i, k_{2} 2^{-i}}= \sum_{\t} \sum_{\xi \in \Xi^{i}_{\t}} u_{\xi, \xi^{2}+k_{2}}$$

Let's take $\xi \in \Xi^{i}_{\t}$ and $\eta \in \Xi_{n_{1}}$ and denote by $\a$ their angle. From (\ref{d1}) we obtain that if the support of $\hat{v}_{\eta, \leq 2^{-i}} * \hat{u}_{\xi, \xi^{2}+ k_{2}}$ is at distance $ \geq 2^{i-2}$ from $P$ then $|\cos{\a}| \geq 2^{-2} $.  

We fix $\t$. For fixed $\eta$, as we change $\xi \in \Xi^{i}_{\t}$ we actually change $|\xi|$ in increments of $2^{-i}$. The support of $\hat{u}_{\xi,\xi^{2}+ k_{2}}$ moves in a direction transversal to the support of $\hat{v}_{\eta, \leq 2^{-i}}$ in increments of $2^{-i}$, therefore the supports of $\hat{v}_{\eta, \leq 2^{-i}} * \hat{u}_{\xi,\xi^{2}+ k_{2} }$ are disjoint with respect to $\xi \in \Xi^{i}_{\t}$. 

Because the support of $\hat{v}_{\eta, \leq 2^{-i}} * \hat{u}_{\xi, \xi^{2}+k_{2}}$ is a translation of the support of $\hat{v}_{\eta, \leq 2^{-i}}$ in a transversal direction to $P$ and $\hat{v}_{\eta, \leq 2^{-i}}$ are supported on $P$, we get orthogonality with respect to $\eta$ too. 

Taking into account the above two remarks and (\ref{m9}) we can estimate:

$$||\sum_{\eta \in \Xi_{n_{1}}, \xi \in \Xi^{i}_{\t}} v_{\eta, \leq 2^{-i}} \cdot u_{\xi, \xi^{2}+k_{2}}||^{2}_{L^{2}} \leq \sum_{\eta \in \Xi_{n_{1}}, \xi \in \Xi^{i}_{\t}} ||v_{\eta, \leq 2^{-i}} \cdot u_{\xi,\xi^{2}+ k_{2}}||^{2}_{L^{2}} \leq $$

$$2^{-i-j} \sum_{\eta \in \Xi_{n_{1}}, \xi \in \Xi^{i}_{\t}} ||v_{\eta, \leq 2^{-i}}||^{2}_{Y_{j}} \cdot|| u_{\xi,\xi^{2}+ k_{2}}||^{2}_{L^{2}} \approx$$

$$2^{-i-j} (\sum_{\eta \in \Xi_{n_{1}}} ||v_{\eta, \leq 2^{-i}}||^{2}_{Y_{j}}) \cdot (\sum_{\xi \in \Xi^{i}_{\t}}|| u_{\xi,\xi^{2}+ k_{2}}||^{2}_{L^{2}}) =$$

$$2^{-i-j} ||v_{n_{1}, \leq 2^{-i}}||^{2}_{Y_{j}} \cdot (\sum_{\xi \in \Xi^{i}_{\t}}|| u_{\xi,\xi^{2}+ k_{2}}||^{2}_{L^{2}})$$

We do not have orthogonality of the interaction with respect to $\t$, hence:

$$||v_{n_{1}, \leq 2^{-i}} \cdot u_{i, k_{2} 2^{-i}}||_{L^{2}} \leq \sum_{\t} ||\sum_{\eta \in \Xi_{n_{1}}, \xi \in \Xi^{i}_{\t}} v_{\eta, \leq 2^{-i}} \cdot u_{\xi,\xi^{2}+ k_{2}}||_{L^{2}} \leq$$

$$2^{i} \left( \sum_{\t} ||\sum_{\eta \in \Xi_{n_{1}}, \xi \in \Xi^{i}_{\t}} v_{\eta, \leq 2^{-i}} \cdot u_{\xi,\xi^{2}+ k_{2}}||^{2}_{L^{2}}\right)^{\q} \leq$$

$$2^{\frac{i-j}{2}} ||v_{n_{1}, \leq 2^{-i}}||_{Y_{j}} \left( \sum_{\t} \sum_{\xi \in \Xi^{i}_{\t}}|| u_{\xi,\xi^{2}+ k_{2}}||^{2}_{L^{2}}  \right)^{\q} \approx 2^{\frac{i-j}{2}} ||v_{n_{1}, \leq 2^{-i}}||_{Y_{j}} || u_{i, k_{2} 2^{-i}}||_{L^{2}} $$ 

The summation with respect to $n_{1}$ has been already discussed (see proof of (\ref{a1})):

$$||v_{j, \leq 2^{-i}} \cdot u_{i, k_{2} 2^{-i}}||_{L^{2}} \leq 2^{\frac{i-j}{2}} ||v_{j, \leq 2^{-i}}||_{Y_{j}} || u_{i, k_{2} 2^{-i}}||_{L^{2}} $$

Therefore:

$$||v_{j, \leq 2^{-i}} \cdot u_{i, k_{2} 2^{-i}}||_{X^{0,-\q}_{j, d 2^{-i}}} \approx 2^{-\frac{i}{2}-\frac{j}{2}}||v_{j, \leq 2^{-i}} \cdot u_{i, k_{2} 2^{-i}}||_{L^{2}} \leq $$

$$2^{-j}||v_{j, \leq 2^{-i}}||_{Y_{j}} || u_{i, k_{2} 2^{-i}}||_{L^{2}}$$

\end{proof}

The next result deals with the case when we have the low frequency input close to $\tau$ axis and we measure parts of the output which are close to $P$. 

\begin{p2} If $2k_{1} \leq k_{2} \leq 2^{2i-2}$ we have 

\beq \label{a11}
|| v_{j, \leq 2^{-i}} \cdot  \sum_{\xi \in \Xi^{i}_{k_{2}}} \sum_{l \in I_{\xi}} u_{\xi,l}||_{\mathcal{R}X^{0,-\q}_{j, k_{1} 2^{-i}}} \leq 
\eeq

$$2^{i-j} (k_{1} k_{2})^{-\q} ||v_{j, \leq 2^{-i}}||_{\mathcal{R}Y_{j}} \cdot ||\sum_{\xi \in \Xi^{i}_{k_{2}}} \sum_{l \in I_{\xi}} u_{\xi,l}||_{\mathcal{R}L^{2}}$$

\end{p2}

\begin{proof}
From (\ref{d1}) we get the necessary and sufficient condition that the support of $\hat{v}_{\eta, \leq 2^{-i}} * \hat{u}_{\xi, l}$ is localized at distance $ \approx 2^{-i}k_{1}$ from $P$:

$$|\cos{\a}| \in k_{2}^{-1} [k_{1}-1,k_{1}+1]$$

\noindent
where, as usual $\a$ is the angle between $\xi$ and $\eta$ and we use the fact that in the support of $\hat{u}_{\xi,l}$ we have $|\xi| \approx 2^{-i} k_{2}$. Since $2k_{1} \leq k_{2}$ we get $|\cos{\a}| \leq \q$.  

This suggest a splitting of $\Xi^{i}_{n_{1}}$ and $\Xi^{i}_{k_{2}}$ in angular sectors of size $\approx k_{2}^{-1}$. From this point on, the steps are exactly as in the proof of (\ref{a9}). 

At the numerology level we record the following changes:

- the gain from spherical symmetry changes now to $k_{2}^{-\q}$. 

- the summation with respect to $l$ brings an additional factor of $2^{i}$.

Other than that the argument is the same as before.

\end{proof}

We need to complete the result in the previous Proposition by analyzing the cases left out. In what follows $d \in \{1,2, 2^{2}, ..., 2^{3i+2}\}$ and we remark that the outcome $\hat{v}_{j, \leq 2^{-i}} * \hat{u}_{i}$ cannot be supported at distance higher than $2^{2i+2}$ from $P$.

\begin{p2}  We have 

\beq \label{a13}
|| v_{j, \leq 2^{-i}} \cdot  \sum_{k \leq 4 d}  \sum_{\xi \in \Xi_{k}} \sum_{l \in I_{\xi}} u_{\xi,l}||_{X^{0,-\q}_{j, d 2^{-i}}} \leq 
\eeq

$$2^{-j} ||v_{j, \leq 2^{-i}}||_{Y_{j}} \cdot ||\sum_{k \leq 4 d}  \sum_{\xi \in \Xi_{k}} \sum_{l \in I_{\xi}} u_{\xi,l}||_{X^{0,\q}}$$

\end{p2}

\begin{proof} We decompose $v_{j, \leq 2^{-i}}$ as in (\ref{dec1}). We fix $n_{1}$ and group the $\eta$'s in $2^{j-2i}$  subsets which have orthogonal interaction with $\sum_{k \leq 4 d}  \sum_{\xi \in \Xi_{k}} \sum_{l \in I_{\xi}} u_{\xi,l}$. Write $\Xi_{n_{1}}=\cup_{m=0}^{2^{j-2i}} A_{m}$ where $A_{m}=\{\eta \in \Xi_{n_{1}}: \arg \eta \in [2^{2i-j}m, 2^{2i-j}(m+1)] \}$. The size of such a block in the tangential direction to $P$ is $\approx 2^{2i}$ which is greater than the size of the support of  $\sum_{k \leq 4 d}  \sum_{\xi \in \Xi_{k}} \sum_{l \in I_{\xi}} u_{\xi,l}$. 

The basic estimate we use is the one in (\ref{m9}) which says:

$$||v_{\eta, \leq 2^{-i}} \cdot u_{\xi,l}||_{L^{2}} \leq 2^{-\frac{i+j}{2}} ||v_{\eta, \leq 2^{-i}} ||_{Y_{j}} ||u_{\xi,l}||_{L^{2}}$$

The rest is just a careful examination of possible orthogonalities we can take advantage while summing up this estimate. We fix $m$ and note that the supports of $\hat{v}_{\eta, \leq 2^{-i}} * \hat{u}_{\xi, l}$ are disjoint with respect to $\eta \in A_{m}$, hence:

$$||\sum_{\eta \in A_{m}} v_{\eta, \leq 2^{-i}} \cdot u_{\xi,l}||_{L^{2}} \leq 2^{-\frac{i+j}{2}} ||\sum_{\eta \in A_{m}} v_{\eta, \leq 2^{-i}} ||_{Y_{j}} ||u_{\xi,l}||_{L^{2}}$$

Then we pick an arbitrary point in the support of $\sum_{\eta \in A_{m}} \hat{v}_{\eta, \leq 2^{-i}}$ and take the normal vector to $P$ at this point and denote by  $n_{m}$. We fix $l$ and notice that the interaction $\sum_{\eta \in A_{m}} \hat{v}_{\eta, \leq 2^{-i}} * \hat{u}_{\xi,l}$ is almost orthogonal with respect to $\xi$ as we move $\xi$ in the direction of $n_{m}$. In the end we have to sum with respect to $\xi$ in the orthogonal direction to $n_{m}$ and with respect to $l$ and there is no more orthogonality we can exploit. Therefore we pick factors of $d^{\q}$ and $2^{i}$ respectively and obtain:

$$||\sum_{\eta \in A_{m}} v_{\eta, \leq 2^{-i}} \cdot \sum_{k \leq 4 d}  \sum_{\xi \in \Xi_{k}} \sum_{l \in I_{\xi}} u_{\xi,l} u_{\xi,l}||_{L^{2}} \leq $$

$$2^{\frac{i-j}{2}} d^{\q} ||\sum_{\eta \in A_{m}} v_{\eta, \leq 2^{-i}} ||_{Y_{j}} ||\sum_{k \leq 4 d}  \sum_{\xi \in \Xi_{k}} \sum_{l \in I_{\xi}} u_{\xi,l} u_{\xi,l}||_{L^{2}}$$

Using the orthogonality with respect to $m$ of the interaction as explained at the beginning helps us to perform the last summation with respect to $m$ from where we obtain the (\ref{a13}) by passing to $X^{s, \pm \q}$.

\end{proof}

At the end of this section we sum up all the estimates. 

\begin{p2}

We have:

\beq \label{a20}
||v_{j, \leq 2^{-i}} \cdot u_{i}||_{\mathcal{R}X_{j,d 2^{-i}}^{0,-\q}} \leq i^{\q} 2^{-j}||v_{j, \leq 2^{-i}}||_{\mathcal{R}Y_{j}} ||u_{i}||_{\mathcal{R}X^{0,\q}}
\eeq

\end{p2}

\begin{proof}
For fixed $d$ we decompose:

$$u_{i}= u_{i, \leq 2^{i-2}} + u_{i, \geq 2^{i-2}}=\sum_{k=1}^{2^{2i-2}} u_{i,k 2^{-i}} + $$

$$ \sum_{k \geq 4d} \sum_{\xi \in \Xi_{k}} \sum_{l \in I_{\xi}} u_{\xi,l} + \sum_{k \leq 4 d}  \sum_{\xi \in \Xi_{k}} \sum_{l \in I_{\xi}} u_{\xi,l}$$

\noindent
where the factors have disjoint support in the frequency space and: 

$$||u_{i}||^{2}_{\mathcal{R}X^{0,\q}} \approx \sum_{k} k ||u_{i,k 2^{-i}}||^{2}_{\mathcal{R}L^{2}} + $$

$$2^{i} \sum_{k \geq 4 d}  \sum_{\xi \in \Xi_{k}} \sum_{l \in I_{\xi}} ||u_{\xi,l}||^{2}_{\mathcal{R}L^{2}} +  2^{i} \sum_{k \leq 4 d}  \sum_{\xi \in \Xi_{k}} \sum_{l \in I_{\xi}} ||u_{\xi,l}||^{2}_{\mathcal{R}L^{2}}$$

If $d \leq 2^{2i-2}$ we make use of (\ref{a9}) to derive:

$$|| v_{j, \leq 2^{-i}} \cdot  u_{i, \leq 2^{i-2}}||^{2}_{\mathcal{R}X^{0,-\q}_{j,d 2^{-i}}} \approx \sum_{k_{1}=2^{-1}d}^{2d} || v_{j, \leq 2^{-i}} \cdot  u_{i, \leq 2^{i-2}}||^{2}_{\mathcal{R}X^{0,-\q}_{j,2^{-i}k_{1}}} \leq$$

$$ \sum_{k_{1}=2^{-1}d}^{2d} \left( \sum_{k_{2}=1}^{2^{2i-2}} || v_{j, \leq 2^{-i}} \cdot  u_{i, k_{2} 2^{-i}}||_{\mathcal{R}X^{0,-\q}_{j,2^{-i}k_{1}}} \right)^{2} \leq$$

$$ \sum_{k_{1}=2^{-1}d}^{2d} \left( \sum_{k_{2}=1}^{2^{2i-2}} 2^{-j} k_{1}^{-\q} ||v_{j, \leq 2^{-i}}||_{\mathcal{R}Y_{j}} \cdot ||u_{i, k_{2} 2^{-i}}||_{\mathcal{R}L^{2}} \right)^{2} \leq $$

$$ \sum_{k_{1}=2^{-1}d}^{2d} \left( \sum_{k_{2}=1}^{2^{2i-2}} 2^{-j} k_{1}^{-\q}  k_{2}^{-\q}||v_{j, \leq 2^{-i}}||_{\mathcal{R}Y_{j}} \cdot ||u_{i, k_{2} 2^{-i}}||_{\mathcal{R}X^{0,\q}} \right)^{2} \leq $$

$$2^{-2j}  (\sum_{k_{2}=1}^{2^{2i-2}} k_{2}^{-1})^{\q} (\sum_{k_{1}=2^{-1}d}^{2d} k_{1}^{-1} ) || v_{j, \leq 2^{-i}}||^{2}_{\mathcal{R}Y_{j}} ||u_{i, \leq 2^{i-2}}||_{\mathcal{R}X^{0,\q}}^{2} \leq$$

$$ i 2^{-2j} || v_{j, \leq 2^{-i}}||^{2}_{\mathcal{R}Y_{j}} ||u_{i, \leq 2^{i-2}}||_{\mathcal{R}X^{0,\q}}^{2}$$

If $d \geq 2^{2i-2}$ we use (\ref{a50}) to derive in a similar way:

$$|| v_{j, \leq 2^{-i}} \cdot  u_{i, \leq 2^{i-2}}||_{X^{0,-\q}_{j,d 2^{-i}}}  \leq i^{\q} 2^{-j} ||v_{j, \leq 2^{-i}}||_{Y_{j}} ||u_{i, \leq 2^{i-2}}||_{X^{0,\q}}$$

For the second part we use (\ref{a11}) to obtain:

$$|| v_{j, \leq 2^{-i}} \cdot  \sum_{k_{2} \geq 4d}  \sum_{\xi \in \Xi_{k_{2}}} \sum_{l \in I_{\xi}} u_{\xi,l}||^{2}_{\mathcal{R}X^{0,-\q}_{j, 2^{-i}d}} \leq $$

$$\sum_{k_{1}=2^{-1}d}^{2d} || v_{j, \leq 2^{-i}} \cdot \sum_{k_{2} \geq 4d}  \sum_{\xi \in \Xi_{k_{2}}} \sum_{l \in I_{\xi}} u_{\xi,l} ||^{2}_{\mathcal{R}X^{0,-\q}_{j, 2^{-i}k_{1}}} \leq$$

$$\sum_{k_{1}=2^{-1}d}^{2d} \left( \sum_{k_{2} \geq 4d} || v_{j, \leq 2^{-i}} \cdot \sum_{\xi \in \Xi_{k_{2}}} \sum_{l \in I_{\xi}} u_{\xi,l}||_{\mathcal{R}X^{0,-\q}_{j, 2^{-i}k_{1}}} \right)^{2}  \leq$$

$$\sum_{k_{1}=2^{-1}d}^{2d} \left( \sum_{k_{2} \geq 4d} 2^{i-j} k_{1}^{-\q} k_{2}^{-\q} || v_{j, \leq 2^{-i}}||_{\mathcal{R}Y_{j}} \cdot ||\sum_{\xi \in \Xi_{k_{2}}} \sum_{l \in I_{\xi}} u_{\xi,l}||_{\mathcal{R}L^{2}} \right)^{2}  \approx$$

$$\sum_{k_{1}=2^{-1}d}^{2d} \left( \sum_{k_{2} \geq 4d} 2^{-j} k_{1}^{-\q} k_{2}^{-\q} || v_{j, \leq 2^{-i}}||_{\mathcal{R}Y_{j}} \cdot ||\sum_{\xi \in \Xi_{k_{2}}} \sum_{l \in I_{\xi}} u_{\xi,l}||_{\mathcal{R}X^{0,\q}} \right)^{2}  \leq$$

$$i 2^{-2j} || v_{j, \leq 2^{-i}}||^{2}_{\mathcal{R}Y_{j}} \cdot ||\sum_{k_{2} \geq 4d}  \sum_{\xi \in \Xi_{k_{2}}} \sum_{l \in I_{\xi}} u_{\xi,l}||^{2}_{\mathcal{R}X^{0,\q}}$$

Finally, for the third part (\ref{a13}) is exactly what we need. Summing up all these estimates gives us the result claimed.

\end{proof}

\begin{p2}

We have:

\beq \label{a21}
||v_{j, \leq 2^{-i}} \cdot u_{i}||_{\mathcal{R}X_{j, \geq 2^{-i}}^{0,-\q,1}} \leq i^{\frac{3}{2}} 2^{-j}||v_{j, \leq 2^{-i}}||_{\mathcal{R}Y_{j}} ||u_{i}||_{\mathcal{R}X^{0,\q,1}}
\eeq

\end{p2}

\begin{proof}

$$S_{j,\geq 2^{-i}} (v_{j, \leq 2^{-i}} \cdot u_{i}) = \sum_{d} S_{j, d2^{-i}} (v_{j, \leq 2^{-i}} \cdot u_{i})$$

Because $5i \leq j$, we have that $\hat{u}_{i}$ cannot move the support of $\hat{v}_{j, \leq 2^{-i}}$ with respect to the distance to $P$ with more than $2^{2i+2}$. Therefore the set of values for $d$ in the above sum is $\{1, 2, 2^{2}, ..., 2^{3i+2}\}$ which implies it contains $\approx 3i$ values.  We apply (\ref{a20}) for each $d$, perform the summation with respect to $d$ and get the claim in (\ref{a21}). 

\end{proof}

\subsection{Estimates: {\mathversion{bold}$ \mathcal{R}X^{0,\q,1}_{i} \cdot \mathcal{R}X^{0,\q,1}_{j, \geq 2^{-i}} \rightarrow \mathcal{R}\mathcal{Y}_{j, \leq 2^{-i}}$}}

\noindent
\vspace{.1in}

These estimates can be obtained by duality from the ones in $ \mathcal{R} X^{0,\q,1}_{i} \cdot \mathcal{R}Y_{j, \leq 2^{-i}} \rightarrow \mathcal{R}X^{0,-\q,1}_{j, \geq 2^{-i}}$. We state the main result we need:

\begin{p2}  We have the estimate:

\beq \label{a25}
||v_{j, \geq 2^{-i}} \cdot u_{i}||_{\mathcal{R}\mathcal{Y}_{j, \leq 2^{-i}}} \leq i^{\frac{3}{2}} 2^{-j} ||v_{j, \geq 2^{-i}}||_{\mathcal{R}X^{0,\q,1}} \cdot ||u_{i}||_{\mathcal{R}X^{0,\q,1}}
\eeq

\end{p2}

\subsection{Bilinear estimates on dyadic regions}

\begin{proof}[Proof of Theorem \ref{tb2}]

$$B(u_{i},v_{j}) = S_{j, \leq 2^{-i}}B(u_{i},v_{j}) +  S_{j, \geq 2^{-i}}B(u_{i},v_{j}) =$$

$$ S_{j, \leq 2^{-i}}B(u_{i},v_{j, \leq 2^{-i}}) + S_{j, \leq 2^{-i}}B(u_{i},v_{j, \geq 2^{-i}}) + $$

$$ S_{j, \geq 2^{-i}}B(u_{i},v_{j, \leq 2^{-i}}) + S_{j, \geq 2^{-i}}B(u_{i},v_{j, \geq 2^{-i}})$$

For the first term we make use of (\ref{a7}) and get:

$$|| S_{j, \leq 2^{-i}}B(u_{i},v_{j, \leq 2^{-i}})||_{\mathcal{R}\mathcal{Y}^{s}} \approx 2^{sj} || S_{j, \leq 2^{-i}}B(u_{i},v_{j, \leq 2^{-i}})||_{\mathcal{R}\mathcal{Y}} \leq $$

$$i^{\q} 2^{(s-1)j} ||\nabla v_{j, \leq 2^{-i}}||_{\mathcal{R}Y} ||\nabla u_{i}||_{\mathcal{D}\mathcal{R}X^{0,\q}} \approx $$

$$i^{\q} 2^{i} 2^{sj} ||v_{j, \leq 2^{-i}}||_{\mathcal{R}Y} ||u_{i}||_{\mathcal{D}\mathcal{R}X^{0,\q}} \approx $$

$$i^{\q} 2^{(1-s)i} ||v_{j, \leq 2^{-i}}||_{\mathcal{R}Y^{s}} ||u_{i}||_{\mathcal{D}\mathcal{R}X^{s,\q}} \leq i 2^{(1-s)i} ||v_{j, \leq 2^{-i}}||_{\mathcal{R}Y^{s}} ||u_{i}||_{\mathcal{D}\mathcal{R}Z^{s}}$$

For the second term we make use of (\ref{a25}) to get:

$$||S_{j, \leq 2^{-i}}B(u_{i},v_{j, \geq 2^{-i}})||_{\mathcal{R}\mathcal{Y}^{s}} \approx 2^{sj} ||S_{j, \geq 2^{-i}}B(u_{i},v_{j, \geq 2^{-i}})||_{\mathcal{R}\mathcal{Y}} \leq $$

$$i^{\q}  2^{(s-1)j} ||\nabla v_{j, \geq 2^{-i}}||_{\mathcal{R}X^{0,\q}} \cdot ||\nabla u_{i}||_{\mathcal{D}\mathcal{R}X^{0,\q}} \approx $$

$$i^{\q}  2^{(1-s)i} ||v_{j, \geq 2^{-i}}||_{\mathcal{R}X^{s,\q}} \cdot || u_{i}||_{\mathcal{D}\mathcal{R}X^{s,\q}} \leq $$

$$i^{\frac{3}{2}}  2^{(1-s)i} ||v_{j, \geq 2^{-i}}||_{\mathcal{R}Z^{s}} \cdot || u_{i}||_{\mathcal{D}\mathcal{R}Z^{s}} $$

For the third term we use of (\ref{a21}) to obtain:

$$||S_{j, \geq 2^{-i}}B(u_{i},v_{j, \leq 2^{-i}})||_{\mathcal{R}X^{s,-\q,1}} \approx 2^{sj} ||S_{j, \geq 2^{-i}}B(u_{i},v_{j, \leq 2^{-i}})||_{\mathcal{R}X^{0,-\q,1}} \leq $$

$$i^{\frac{3}{2}} 2^{(s-1)j} ||\nabla v_{j, \leq 2^{-i}}||_{\mathcal{R}Y} ||\nabla u_{i}||_{\mathcal{D}\mathcal{R}X^{0,\q}} \approx $$

$$i^{\frac{3}{2}} 2^{(1-s)i} ||v_{j, \leq 2^{-i}}||_{\mathcal{R}Y^{s}} ||u_{i}||_{\mathcal{D}\mathcal{R}X^{s,\q}} \leq i^{2} 2^{(1-s)i} ||v_{j, \leq 2^{-i}}||_{\mathcal{R}Y^{s}} ||u_{i}||_{\mathcal{D}\mathcal{R}Z^{s}} $$

The fourth term had been handled in Theorem \ref{tb1}. By adding all the estimates we obtain:

$$||B(u,v)||_{\mathcal{R}W^{s}_{j}} \leq i^{\frac{3}{2}} 2^{(1-s)i} ||u||_{\mathcal{R}\mathcal{D}Z^{s}_{i}} ||v||_{\mathcal{R}Z^{s}_{j}}$$

In the end we can recover the decay via an argument similar to the one in (\ref{b30}) in the case when the outcome is localized at high frequency. One would notice that over there we had to recover decay of type $d_{j}$ and all we used is that the high frequency comes with that decay. This is why we were allowed to make use of the decay property of the low frequency throughout the argument so far.

\end{proof}

\section{Bilinear estimates}

\begin{proof}[Proof of Theorem \ref{bg}]

We decompose our functions in dyadic pieces to get:

$$B(u,v) = \sum_{k} \sum_{i,j} S_{k} B(u_{i},v_{j}) = \sum_{k}\sum_{i \leq j}S_{k} B(u_{i},v_{j}) +\sum_{k} \sum_{i > j}S_{k} B(u_{i},v_{j})$$

We work out the first term, the second one being similar. We know that if $i < j-1$ then $B(u_{i},v_{j})$ is supported at frequency $2^{j}$ or $2^{j \pm 1}$. Therefore we have:

$$\sum_{i \leq j}S_{k} B(u_{i},v_{j}) = \sum_{i \leq k-2}S_{k} B(u_{i}, v_{k-1} + v_{k}+v_{k+1}) +$$

$$ \sum_{k \leq j}S_{k} B(u_{j-1}+u_{j}+u_{j+1}, v_{j-1}+v_{j}+v_{j+1})$$

From the bilinear estimates on dyadic pieces, see (\ref{b30}) and (\ref{b7}) we get two sets on inequalities:

$$||\sum_{i \leq k-2}S_{k} B(u_{i}, v_{k-1} + v_{k}+v_{k+1})||_{\mathcal{R}\mathcal{D}W^{s}} \leq $$

$$ C_{N} \sum_{i \leq k-2} i^{2} 2^{(1-s)i} ||u_{i}||_{\mathcal{R}\mathcal{D}Z^{s}} ||v_{k-1}+v_{k}+v_{k+1}||_{\mathcal{R}\mathcal{D}Z^{s}} \leq$$

$$C_{N} ||u||_{\mathcal{R}\mathcal{D}Z^{s}} ||v_{k-1}+v_{k}+v_{k+1}||_{\mathcal{R}\mathcal{D}Z^{s}}$$

and

$$||\sum_{k \leq j}S_{k} B(u_{j-1}+u_{j}+u_{j+1}, v_{j-1}+v_{j}+v_{j+1})||_{\mathcal{R}\mathcal{D}W^{s}} \leq $$

$$\sum_{k \leq j} j^{2} 2^{(1-s)j} 2^{(\q+\e-s)(j-k)} ||u_{j-1}+u_{j}+u_{j+1}||_{\mathcal{R}\mathcal{D}Z^{s}} ||v_{j-1}+v_{j}+v_{j+1}||_{\mathcal{R}\mathcal{D}Z^{s}} \leq$$

$$\sum_{k \leq j} 2^{(\q+\e-s)(j-k)} ||u_{j-1}+u_{j}+u_{j+1}||_{\mathcal{R}\mathcal{D}Z^{s}} ||v_{j-1}+v_{j}+v_{j+1}||_{\mathcal{R}\mathcal{D}Z^{s}} \leq $$

$$C_{\e,s} \left( \sum_{k \leq j} 2^{(1+2\e-2s)(j-k)} ||u_{j-1}+u_{j}+u_{j+1}||^{2}_{\mathcal{R}\mathcal{D}Z^{s}} ||v_{j-1}+v_{j}+v_{j+1}||^{2}_{\mathcal{R}\mathcal{D}Z^{s}} \right)^{\q}$$

In the last line we used Cauchy-Schwartz and the estimate:

$$\sum_{k \leq j} 2^{(1+2\e-2s)(j-k)} \leq C^{2}_{\e,s}$$

\noindent
which is valid as long as $1 < s$ and $\e < \q$ since these imply that $1+2\e-2s > 0$. 

In the end we sum up with respect to $k$:

$$||\sum_{k}\sum_{i \leq j}S_{k} B(u_{i},v_{j})||^{2}_{\mathcal{R}\mathcal{D}W^{s}} = \sum_{k} ||\sum_{i \leq j}S_{k} B(u_{i},v_{j})||^{2}_{\mathcal{R}\mathcal{D}W^{s}} \leq $$

$$\sum_{k} ||u||^{2}_{\mathcal{R}\mathcal{D}Z^{s}} ||v_{k-1}+v_{k}+v_{k+1}||^{2}_{\mathcal{R}\mathcal{D}Z^{s}} + $$

$$C^{2}_{\e,s} \sum_{k} \sum_{k \leq j} 2^{(1+2\e-2s)(j-k)} ||u_{j-1}+u_{j}+u_{j+1}||^{2}_{\mathcal{R}\mathcal{D}Z^{s}} ||v_{j-1}+v_{j}+v_{j+1}||^{2}_{\mathcal{R}\mathcal{D}Z^{s}} \leq $$

$$ C^{4}_{\e,s} ||u||^{2}_{\mathcal{R}\mathcal{D}Z^{s}} ||v||^{2}_{\mathcal{R}\mathcal{D}Z^{s}}$$

The term in the third line of the previous estimate was handled by an change in the order of summation: 

$$ \sum_{k} \sum_{k \leq j} 2^{(1+2\e-2s)(j-k)} ||u_{j-1}+u_{j}+u_{j+1}||^{2}_{\mathcal{R}\mathcal{D}Z^{s}} ||v_{j-1}+v_{j}+v_{j+1}||^{2}_{\mathcal{R}\mathcal{D}Z^{s}} =$$

$$ \sum_{j} \sum_{k \leq j} 2^{(1+2\e-2s)(j-k)} ||u_{j-1}+u_{j}+u_{j+1}||^{2}_{\mathcal{R}\mathcal{D}Z^{s}} ||v_{j-1}+v_{j}+v_{j+1}||^{2}_{\mathcal{R}\mathcal{D}Z^{s}}  \leq $$

$$ C^{2}_{\e,s} \sum_{j} ||u_{j-1}+u_{j}+u_{j+1}||^{2}_{\mathcal{R}\mathcal{D}Z^{s}} ||v_{j-1}+v_{j}+v_{j+1}||^{2}_{\mathcal{R}\mathcal{D}Z^{s}}  \leq $$

$$ C^{2}_{\e,s} ||u||^{2}_{\mathcal{R}\mathcal{D}Z^{s}} ||v||^{2}_{\mathcal{R}\mathcal{D}Z^{s}}$$

The estimate for $B(u,\bar{v})$ is obtained in a similar way. 

\end{proof}

\end{document}